\documentclass[11pt]{amsart}	
\usepackage[utf8]{inputenc}
\usepackage[T1]{fontenc}
\usepackage{amsmath,amsfonts,amssymb,dsfont}
\usepackage{graphicx}
\usepackage[english]{babel}
\usepackage{listings}
\usepackage{xspace}
\usepackage[colorinlistoftodos]{todonotes}
\usepackage{pgf,pgfkeys,pgfmath}
\usepackage{upgreek}

\usepackage{fullpage}
\usepackage[backref=page]{hyperref}
\hypersetup{
	plainpages=false,
	unicode=false,          
	pdftoolbar=true,        
	pdfmenubar=true,        
	pdffitwindow=true,      
	pdftitle={My title},    
	pdfauthor={Author},     
	pdfsubject={Subject},   
	pdfcreator={Creator},   
	pdfproducer={Producer}, 
	pdfkeywords={keywords}, 
	pdfnewwindow=true,      
	colorlinks=true,       
	linkcolor=blue,          
	citecolor=blue,        
	filecolor=blue,      
	urlcolor=magenta,          
	urlbordercolor=0 1 1
}

\usepackage{amsthm}

\usepackage{lmodern}

\usepackage{color}
\usepackage{xcolor}
  
\usepackage{pgf, tikz}

\usepackage{fancybox} 

\def\R{\mathbb{R}}
\def\N{\mathbb{N}}

\def\K{\mathbb{K}}

\def\H{\mathcal{H}}

\def\D{\mathcal{D}}
\def\Re{\mathcal{R}}

\def\S{\mathcal{S}}

\def\O{\mathcal{O}}

\DeclareMathOperator*{\argmin}{argmin}

\usepackage{mathtools}
\DeclarePairedDelimiter{\abs}{\lvert}{\rvert}
\makeatletter
\let\oldabs\abs
\def\abs{\@ifstar{\oldabs}{\oldabs*}}
\makeatother

\newtheorem{thm}{Theorem}[section]

\newtheorem{rem}[thm]{Remark}

\def\dt{\delta_t}

\usepackage{palatino}

\title{Learning phase field mean curvature flows with neural networks} 

\author{Elie Bretin}
\address{Univ Lyon, INSA de Lyon, CNRS UMR 5208, Institut Camille Jordan\\ 20 avenue Albert Einstein, F-69621 Villeurbanne, France\\ elie.bretin@insa-lyon.fr}
\author{Roland Denis}
\address{Univ Lyon, CNRS UMR 5208, Institut Camille Jordan\\43 boulevard du 11 novembre
1918, F-69622 Villeurbanne, France\\denis@math.univ-lyon1.fr}
\author{Simon Masnou}
\address{Univ Lyon, Universit\'e Claude Bernard Lyon 1, CNRS UMR 5208, Institut Camille Jordan \\43 boulevard du 11 novembre
1918, F-69622 Villeurbanne, France\\masnou@math.univ-lyon1.fr}
\author{Garry Terii}
\address{Univ Lyon, Universit\'e Claude Bernard Lyon 1, CNRS UMR 5208, Institut Camille Jordan\\43 boulevard du 11 novembre
1918, F-69622 Villeurbanne, France\\ terii@math.univ-lyon1.fr}
\makeatletter
\@namedef{subjclassname@2020}{\textup{2020} Mathematics Subject Classification}
\makeatother

\begin{document}
\maketitle

\begin{abstract}
We introduce in this paper new and very effective numerical methods based on neural networks for the approximation of the mean curvature flow of either oriented or non-orientable surfaces. To learn the correct interface evolution law, our neural networks are trained on phase field representations of exact evolving interfaces. The structures of the networks draw inspiration from splitting schemes used for the discretization of the Allen-Cahn equation. But when the latter approximate the mean curvature motion of oriented interfaces only, the approach we propose extends very naturally to the non-orientable case. Through a variety of examples, we show that our networks, trained only on flows of smooth and simplistic interfaces, generalize very well to more complex interfaces, either oriented or non-orientable, and possibly with singularities. Furthermore, they can  be coupled easily with additional constraints which opens the way to various applications illustrating the flexibility and effectiveness of our approach: mean curvature flows with volume constraint, multiphase mean curvature flows, numerical approximation of Steiner trees or minimal surfaces.
\end{abstract}


\section{Introduction}

Many applications in physics, biology, mechanics, engineering, or image processing involve interfaces 
whose shape evolves to decrease a particular surface energy.  \textcolor{black}{A very common example is the area energy that explains for instance the shapes of soap bubbles, bee honeycombs, the interface between two fluids, some crystalline materials, some communication networks, etc.} In image processing, the area energy is used to quantify regularity, 
for instance in the celebrated TV (total variation) model. It is well-known that the $L^2$-gradient flow of the area energy,
i.e. the flow which decreases the energy in the direction of steepest descent with respect to the $L^2$ metric, is the mean curvature flow. \\

The mean curvature flow is classically defined for smooth, embedded $(N-1)$-surfaces in $\R^N$ without \textcolor{black}{a}  boundary: each point of the surface moves with the velocity vector equal to the (vector) mean curvature, see~\cite{BellettiniBook}. Such a flow is well defined until the onset of singularities. Various definitions of mean curvature flows have been proposed to handle singularities as well, but also to handle non orientable sets (typically planar networks with triple points satisfying the Herring's condition~\cite{MR3967812}), higher codimensional sets, or sets with boundaries, see the many references in~\cite{BellettiniBook}.\\

We propose in this paper to learn with neural networks a numerical approximation of a phase field representation of either the mean curvature flow of an oriented set, or the mean curvature flow of a possibly non orientable interface. Interestingly, as will be seen later, we train our neural networks on few examples of smooth sets flowing by mean curvature but, once trained, some networks can handle consistently sets with singularities as well.
\\

There is a vast literature on the numerical approximation of mean curvature flows, and the methods roughly divide in four categories (some of them are exhaustively reviewed and compared in~\cite{review_interface}):

\begin{enumerate}
\item Parametric methods \cite{Deckelnick2005,Barrett2008_2} are based on explicit parameterizations of smooth surfaces. The numerical approximation of the parametric mean curvature flow is quite simple in 2D and the approach can be extended to non-orientable surfaces since there
is no need for the surface to be the boundary of a domain nor to separate its interior from its exterior. The numerical approximation is however more difficult in dimensions higher than $2$ for the method can hardly handle topological changes. The processing of singularities is difficult even in dimension $2$, see for instance the recent work~\cite{SMAI-JCM_2021__7__27_0}.
  
\item The level set method was introduced by \textcolor{black}{ Osher and  Sethian}~\cite{OsherSethian} for interface geometric evolution problems, see also~\cite{Osherbook1,Osherbook2,Evans_spruck,chen_giga_goto}.
The main idea is to represent implicitly the interface as the zero level-set of an auxiliary 
function $\varphi$ (typically the signed distance function associated with the domain enclosed by the interface) and the evolution is described 
through a Hamilton-Jacobi equation satisfied by $\varphi$.
The level set approach provides a convenient formalism to represent the mean curvature
flow in any dimension. In a strong contrast with the parametric approach, it can handle topological changes and it can be defined rigorously beyond singularities using the theory of viscosity solutions for the Hamilton-Jacobi equation. There are however several difficult issues regarding the numerical approximation of the level set approach. First, the Hamilton-Jacobi equation is nonlinear and highly degenerate, thus difficult to approximate numerically and, secondly, delicate methods are needed to preserve some needed properties of the level set function $\varphi$. Moreover, the method is basically designed to represent the evolution of the boundary of an oriented domain and, to the best of our knowledge, there is no level set method that can handle the mean curvature flow of non-orientable interfaces.  
\item Convolution/thresholding type algorithms~\cite{BenceMerrimanOsher,Ishhi_pires_souganidis,Ruuth_efficient} involve a time-discrete scheme alternating the convolution with a suitable kernel of the characteristic  function at time $t_n$ of the domain enclosed by the interface, followed by a thresholding step to define the set at time $t_n+dt$. The asymptotic limit of such a scheme coincides with the smooth mean curvature flow. \textcolor{black}{By reinterpreting convolution/thresholding type algorithms as gradient flows, convergence and stability results have been proved in~\cite{MR3333842,MR3556529} in more general contexts, e.g., for multiphase systems and for anisotropic energies.}

\item Phase field  approaches~\cite{Modica1977,Chen1992a} \textcolor{black}{constitute} the fourth category of methods for the numerical approximation of the mean curvature flow. In these approaches, the sharp interface between is approximated by a smooth transition, the interfacial area is approximated by a smooth energy depending
on the smooth transition, and the gradient flow of this energy appears to be a relatively simple reaction-diffusion system.
Phase field approaches are widely used in physics since the seminal works of van der Waals' on liquid-vapor interfaces (1893), of Ginzburg \& Landau on superconductivity (1950), and of Cahn \& Hilliard on binary alloys (1958). 
However, most phase field methods are designed to approximate the mean curvature flow of the boundary of a domain but cannot handle the evolution of non-orientable interfaces.
\end{enumerate}

\medskip
Our work starts with the following question: can we design and train neural networks to approximate the mean curvature flow of either oriented or non-orientable interfaces? Our strategy is to draw inspiration from phase field approaches and their numerical approximations.

\medskip
\paragraph{{\bf The phase field approach to the mean curvature flow of domain boundaries}}

\medskip
A time-dependent smooth domain $\Omega(t)\subset\R^d$ evolves under the classical mean curvature flow if its inner normal velocity satisfies $V_n(t) = H(t)$,  where $H(t)$ denotes the scalar mean curvature of the boundary $ \partial\Omega(t)$ (with the convention that the scalar mean curvature on the boundary of a convex domain is positive).
This evolution coincides with the $L^2$-gradient flow of the perimeter of $\Omega(t)$    
$$P(\Omega(t)) = \int_{\partial \Omega(t)} d\H^{d-1},$$
where  $\H^{d-1}$ denotes the $(d-1)$-dimensional Hausdorff measure. In the phase field approach, the perimeter functional is approximated (up to a multiplicative constant) in the sense of $\Gamma$-convergence \cite{Modica1977,Chen1992a} by 
the Cahn-Hilliard energy $P_{\varepsilon}$ defined for every smooth function $u$ by 
$$ P_{\varepsilon}(u) = \int_{\R^d} (\varepsilon \frac{|\nabla u|^2}{2} + \frac{1}{\varepsilon} W(u))dx,$$
where $\varepsilon$ is a real parameter which quantifies the accuracy of the approximation, and $W$ is a double-well potential, typically $W(s) = \frac{1}{2} s^2(1-s)^2$. It follows from the $\Gamma$-convergence as $\varepsilon\to 0$ of $P_{\varepsilon}$ to $c_WP$ ($c_W$ is a constant depending only on $W$) that when $u$ is a smooth approximation of the characteristic function of $\Omega(t)$, $ P_{\varepsilon}(u)$ is close to  $c_WP(\Omega(t))$.

The  $L^2$-gradient flow of the Cahn-Hilliard energy $P_{\varepsilon}$ leads to the celebrated Allen-Cahn equation \cite{AmbrosioNotes2000} which reads as, up to a rescaling:

\begin{equation}\label{eq:Allen-Cahn}
\partial_t u  = \Delta u - \frac{1}{\varepsilon^2}W'(u).
\end{equation}

The existence and uniqueness of a solution to the Allen-Cahn equation and the fact that it satisfies a comparison principle are well-known properties, see for instance~\textcolor{black}{{\cite[Thm 32]{AmbrosioNotes2000}}}.\\

With the choice $W(s) = \frac{1}{2} s^2(1-s)^2$, an evolving set associated naturally with the Allen-Cahn equation is
$$
\Omega_{\varepsilon}(t) = \left\{u_\varepsilon(\cdot,t) \geq \frac{1}{2}\right\},
$$
where $u_{\varepsilon}$ denotes the solution to \eqref{eq:Allen-Cahn}  with the well-prepared initial data
\begin{equation}\label{eq:Allen-Cahn-init}
u_{\varepsilon}(x,0) =  q\left(\frac{d(x,\Omega(0))}{\varepsilon}\right).
\end{equation}

 \begin{rem}
 \textcolor{black}{ 
 Note that the gradient of the solution $u_{\varepsilon}$ to the Allen-Cahn equation may cancel on the boundary of the set 
$\Omega_{\varepsilon}$ and thus may lead to a non-smooth interface. This is typically the case in dimension $3$ where singularities may appear in finite time. 
 }\end{rem}

Here, $d(\cdot, \Omega(0))$ denotes the signed distance 
function to $\Omega(0)$ with the convention that $d(\cdot, \Omega(0))<0$ in $\Omega(0)$ and $q:\R\to [0,1]$ is a so-called \textit{optimal profile} which minimizes the parameter-free one-dimensional Allen-Cahn energy under some constraints:

$$
    q = \argmin_{p} \left\{ \int_{\R} (\frac{|p'(s)|^2}2+{W(p(s))}) ds , ~ p\in C^{0,1}(\R), ~ p(-\infty) = 1, ~ p(0)=\frac 1 2,~ p(+\infty) = 0\right\}.
    $$
Note that in the particular case $W(s) = \frac{1}{2} s^2(1-s)^2$, one has $q(t)=\frac 1 2(1-\tanh(\frac t 2))$.
The initial condition $u_{\varepsilon}(x,0)=q\left(\frac{d(x,\Omega(0))}{\varepsilon}\right)$ is considered as well-prepared initial data because $q\left(\frac{d(x,\Omega(0))}{\varepsilon}\right)$ is almost energetically optimal: with a suitable gluing $\tilde q$ of $q$ to $1$ on the left and $0$ on the right, one gets that, as $\varepsilon\to 0^+$, 
$$P_\varepsilon\left(\tilde q\left(\frac{d(x,\Omega(0))}{\varepsilon}\right)\right)\longrightarrow c_WP(\Omega(0)).$$

A formal asymptotic expansion of the solution $u_{\varepsilon}$ to \eqref{eq:Allen-Cahn}-\eqref{eq:Allen-Cahn-init} near the associated interface $\partial\Omega_{\varepsilon}(t)=\partial \left\{u_\varepsilon(\cdot,t) \geq \frac{1}{2}\right\}$ gives (see~\cite{BellettiniBook})

\begin{equation}\label{dev_asympto}
    u_{\varepsilon} (x, t) = q\left(\frac{d(x,\Omega_{\varepsilon}(t))}{\varepsilon}\right) + \O(\varepsilon^2),
\end{equation}
which shows that $u_\varepsilon$ remains energetically quasi-optimal with second-order accuracy. Furthermore, the velocity $V_\varepsilon$ of the boundary $\partial\Omega_{\varepsilon}(t)$ satisfies
$$
V_{\varepsilon}(t) = H_\varepsilon(t) + \O(\varepsilon^2),
$$
where $H_\varepsilon(t)$ denotes the scalar mean curvature on $\partial\Omega_{\varepsilon}(t)$, which suggests that the Allen-Cahn equation approximates a mean curvature flow with an error of order $\varepsilon^2$~\cite{BellettiniBook}.

Rigorous proofs of the convergence to the smooth mean curvature flow for short times (in particular before the onset of singularities)
have been presented in \cite{Chen1992a, DeMottoniSchatzman1995, bellettini1995quasi} with a quasi-optimal error on the Hausdorff distance between
$\Omega(t)$ and $\Omega_{\varepsilon}(t)$ given by
$$ \text{dist}_{\H}(\Omega(t),\Omega_{\varepsilon}(t)) \leq C \epsilon^2 |\log(\varepsilon)|^2,$$
where the constant $C$ depends on the regularity of $\Omega(t)$. \\

These convergence results, combined with the good suitability of the Allen-Cahn equation for numerical approximation, make the phase field approach a very effective method to approximate the mean curvature flow. This is true however only for the mean curvature motion of domain boundaries, i.e. codimension $1$ orientable interfaces without boundary. 
There do exist some phase field energies to approximate the perimeter of non-orientable interfaces, as for instance in the Ambrosio-Torterelli functional \cite{MR1075076}, but these energies need to be coupled with additional terms and cannot be used to approximate directly the mean curvature motion of the interfaces. \\

 \medskip
 
\paragraph{{\bf Neural networks and phase field representation}}

We introduce in this paper neural networks that learn, at least approximately, how to move a set by mean curvature. Our networks are trained on a collection of sets whose motion by mean curvature is, preferentially, known exactly.  However, these networks are not designed to work with exact sets, but rather with an implicit, smooth representation of them. A key aspect of our approach is the choice of this representation. A first option (see Table~\ref{table:flow}, left column) is to choose the representation provided by the Allen-Cahn phase field approach, i.e., given a set $\Omega(0)$, we compute the solution $u_\varepsilon(\cdot,t)$ to the Allen-Cahn equation~\eqref{eq:Allen-Cahn} with initial data $u_\varepsilon(\cdot,0)=q\left(\frac{d(\cdot, \Omega(0))}{\varepsilon}\right)$. Recall that with the particular choice $W(s)= \frac{1}{2} s^2(1-s)^2$, the $1/2$--isosurface of $u_\varepsilon(\cdot,0)$ is {\it exactly} $\partial\Omega(0)$. However, using the Allen-Cahn phase field approach introduces a bias: if $t\mapsto\Omega(t)$ denotes the motion by mean curvature starting from $\Omega(0)$ (Table~\ref{table:flow}, middle column), the 
$1/2$--isosurface of $u_\varepsilon(\cdot,t)$ is no more than a good approximation of $\partial\Omega(t)$, in general it is different (and the same holds for other isosurfaces of $u_\varepsilon(\cdot,t)$). Instead, we use another phase field representation which introduces no bias, see Table~\ref{table:flow}, right column: for all $t$ (before the onset of singularities), the $1/2$--isosurface at time $t$ of $q\left(\frac{d(\cdot, \Omega(t))}{\varepsilon}\right)$ is {\it exactly} $\partial\Omega(t)$ (still assuming that $W(s)= \frac{1}{2} s^2(1-s)^2$ but the argument can obviously be adapted for other choices of $W$). With such a choice, we ensure that our neural networks will be trained on exact implicit representations of the sets moving by mean curvature.

\begin{table}[htbp]
    \centering
    \begin{tikzpicture}
    
    \node[text width=4cm] (MCF_0) at (5.6, 6.4) {\Ovalbox{$\Omega(0)$}};
    \node[text width=4cm] (MCF) at (4.5, 5.1) {Mean curvature flow};
    \node[text width=4cm] (MCF_t) at (5.6, 3.6) {\Ovalbox{$\Omega(t)$}};
    \node[text width=4cm] (MCF) at (4.5, 2.3) {Mean curvature flow};
    \node[text width=4cm] (MCF_dt) at (5.3, 0.8) {\Ovalbox{$\Omega(t + \dt)$}};
    
    \node[text width=4cm] (PF_0) at (-0.6, 6.4) {\Ovalbox{$u_{\varepsilon}(\cdot, 0)=q\left(\frac{d(\cdot, \Omega(0))}{\varepsilon}\right)$}};
    \node[text width=3.5cm] (PF) at (-0.4, 5.1) {Allen-Cahn flow};
    \node[text width=4cm] (PF_t) at (0.3, 3.6) {\Ovalbox{$u_{\varepsilon}(\cdot, t)$}};
    \node[text width=3.5cm] (PF) at (-0.4, 2.3) {Allen-Cahn flow};
    \node[text width=5cm] (PF_dt) at (0.3, 0.8) {\Ovalbox{$u_{\varepsilon}(\cdot, t+\dt)$}};
    
    \node[text width=4cm] (PPFR_0) at (8.8, 6.4) {\Ovalbox{$q\left(\frac{d(\cdot, \Omega(0))}{\varepsilon}\right)$}};
    \node[text width=4cm] (PPFR_t) at (8.8, 3.6) {\Ovalbox{$q\left(\frac{d(\cdot, \Omega(t))}{\varepsilon}\right)$}};
    \node[text width=4cm] (PPFR_dt) at (8.8, 0.8) {\Ovalbox{$q\left(\frac{d(\cdot, \Omega(t+\dt))}{\varepsilon}\right)$}};
 
    \draw[stealth-] (1.5, 6.4) -- (MCF_0);
	\draw[-stealth] (4.8, 6.4) -- (PPFR_0);	
	\draw[-stealth] (4.8, 3.6) -- (PPFR_t);	
	\draw[-stealth] (5.2, 0.8) -- (PPFR_dt);
	
	\draw[dashed] (4.1, 5.9) -- (4.1, 5.3); 
	\draw[dashed,-stealth] (4.1, 4.9) -- (4.1, 4.1); 
	\draw[] (4.1, 3.1) -- (4.1, 2.5); 
	\draw[-stealth] (4.1, 2.1) -- (4.1, 1.3); 
	\draw[] (-1.1, 3.1) -- (-1.1, 2.5); 
	\draw[-stealth] (-1.1, 2) -- (-1.1, 1.3); 
	\draw[dashed] (-1.1, 5.9) -- (-1.1, 5.3); 
	\draw[dashed,-stealth] (-1.1, 4.8) -- (-1.1, 4.1); 

    \end{tikzpicture}
    \caption{{\bf Approximate vs exact phase field representations of mean curvature flows}. Middle column: a smooth motion by mean curvature $t\mapsto\Omega(t)$. Left column: the solution to the Allen-Cahn equation starting from the energetically quasi-optimal and exact phase field representation $q\left(\frac{d(\cdot, \Omega(0))}{\varepsilon}\right)$. Right column: exact phase field representations $q\left(\frac{d(\cdot, \Omega(t))}{\varepsilon}\right)$ of the $\Omega(t)$'s.}
    \label{table:flow}
\end{table}

\begin{figure}[htbp]
    \centering
    \includegraphics[width=0.30\textwidth,height=0.26\textwidth]{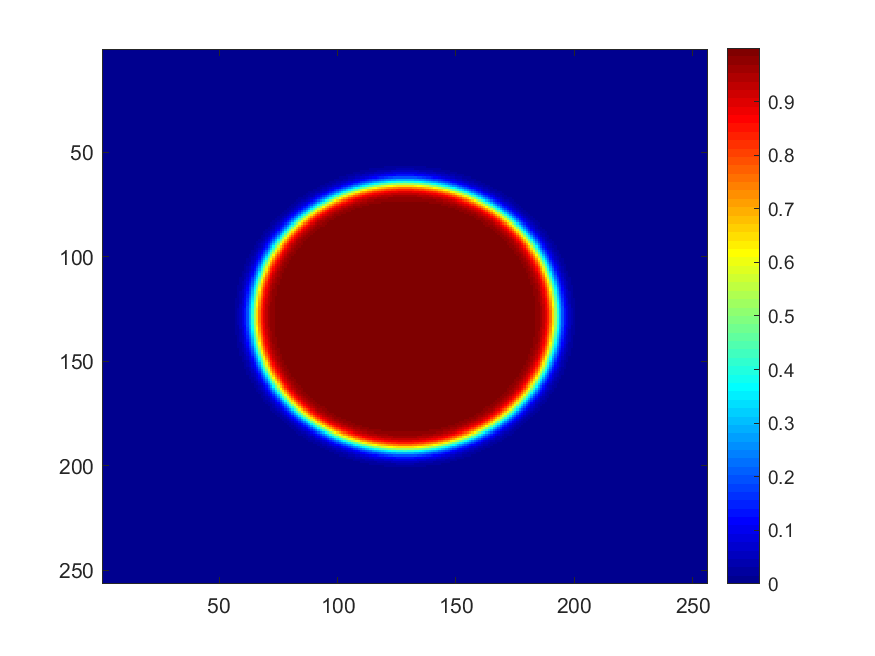}
    \includegraphics[width=0.30\textwidth,height=0.26\textwidth]{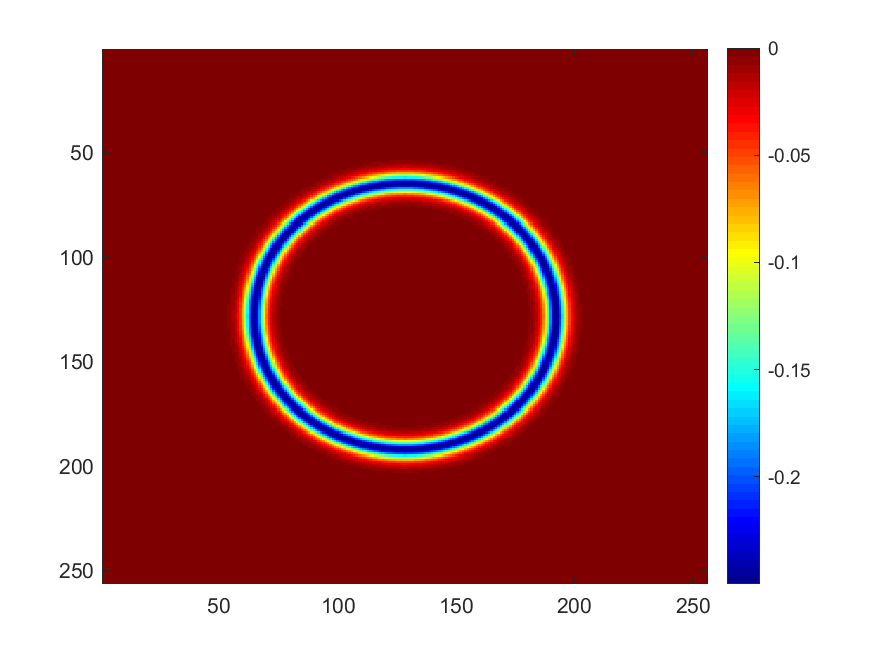}
    \caption{Two phase field representations of the same circle $C$ using the profile $q(t)=\frac 1 2(1-\tanh(\frac t 2))$ and its derivative: on the left, $q\left(\frac{d(\cdot, D}{\varepsilon}\right)$ with $D$ the disk enclosed by $C$. On the right,  $q'\left(\frac{\operatorname{dist}(\cdot, C)}{\varepsilon}\right)$. The first phase field representation carries an orientation information, the second one does not.}
    \label{fig:Phase-field rep}
\end{figure}

It is obviously more accurate to train our networks with the above choice of exact phase field representations, rather than with the Allen-Cahn phase field approximate solutions. Beyond accuracy, there is another key advantage of this choice: the possibility to address situations which are beyond the capacity of the classical Allen-Cahn phase field approach. Indeed, instead of working with the representation $q\left(\frac{d(\cdot, \Omega(t))}{\varepsilon}\right)$ which is well suited for domain boundaries, other implicit representations can be used which do not require any interface orientation. A typical example is the phase field $q'\left(\frac{d(\cdot, \Omega(t))}{\varepsilon}\right)$ where the optimal profile $q$ has been simply replaced by its derivative $q'$. In the case $W(s) = \frac{1}{2} s^2(1-s)^2$, one has $q(s)=\frac 1 2(1-\tanh(\frac s 2))$ and its derivative $q'(s)=\frac 1 4(\tanh^2(\frac s 2)-1)$ is even so  $q'\left(\frac{d(\cdot, \Omega(t))}{\varepsilon}\right)$ is symmetric on both sides of $\partial\Omega(t)$, therefore
$$ q'\left(\frac{d(\cdot, \Omega(t))}{\varepsilon}\right)=q'\left(\frac{\operatorname{dist}(\cdot, \partial\Omega(t))}{\varepsilon}\right)$$
where $\operatorname{dist}$ denotes the classical distance function. Note that $q'\left(\frac{d(\cdot, \Omega(t))}{\varepsilon}\right)$ is another {\it exact} phase field representation: its $(-\frac 1 4)$--isosurface coincides with $\partial\Omega(t)$. However, in contrast with $q\left(\frac{d(\cdot, \Omega(t))}{\varepsilon}\right)$, the phase field is identical on both sides of $\partial\Omega(t)$. To illustrate this idea of another phase field representation which does not carry any orientation information, figure~\ref{fig:Phase-field rep} shows the two phase field representations of a circle obtained using either $q(t)=\frac 1 2(1-\tanh(\frac t 2))$ or $q'$. Obviously, the approach can be adapted for other choices of the profile $q$.

~\\ 
Our contribution in this paper is a new class of neural networks aimed to learn, at least approximately, the flow by mean curvature of either oriented or non orientable sets. Since the classical motion by mean curvature of a smooth set can be defined with nonlinear partial differential equations, our contribution falls in the category of learning methods for PDEs. Let us give a  brief overview of the literature on neural network-based numerical methods for approximating the solutions to PDEs. 

\begin{itemize} 
\item A first category of approaches relies on the fact that very general functions can be approximated with neural networks~\cite{cybenko,KidgerLyons} so it is natural to seek the solution to a PDE as a neural network~ \cite{guo2016convolutional, zhu2018bayesian, adler2017solving, bhatnagar2019prediction, khoo2021solving,MR3847747}. Such a method is accurate and useful for specific problems, but not convenient: for each new initial condition, a new neural network must be trained.  \textcolor{black}{This approach was used for instance in \cite{MR3847747} to approximate numerically in high dimension the solutions to nonlinear PDEs with the particular example of the Allen-Cahn equation. It is not well suited to our purposes as it requires at each time to train a new network to produce an approximation of the solution at time $t_{n+1}$ from the solution at time $t_n$. In addition, the method requires that the specific structure of the equation is known, which is limited when the underlying equation is only partially known.}
Another approach considers a neural network as an operator between Euclidean spaces of same dimension depending on the discretization of the PDE \cite{raissi2019physics, bar2019unsupervised, smith2020eikonet, pan2020physics, yu2018deep}. This approach relies on the discretization and requires to modify the architecture of the network when the discrete resolution or the discretization are changed; 
    \item There are strong connections between neural networks and numerical schemes \cite{lu2018beyond, alt2021connections,doi:10.1137/20M1377199,pock2}. Most neural networks can actually be seen as operators acting between infinite-dimensional spaces (typically spaces of functions): for instance, for a time-dependent PDE, the forward propagation of an associated neural network can be viewed as the flow associated to the PDE when a time-step $\dt$ is fixed. Consequently, the neural network is trained only once: for each new initial condition, the solution is obtained by applying the neural network to the initial condition. This reduces significantly the computational cost in comparison with the first category of approaches mentioned above. This second category of approaches is mesh-free and fully data-driven, i.e., the training procedure, as well as the neural network, do not require any knowledge of the underlying PDE, only the knowledge of particular solutions to the PDE. This can be very advantageous when very limited information is given about the PDE, which is often the case in physics.  Recently, several works have followed this idea, see \cite{lu2019deeponet, bhattacharya2020model, nelsen2021random, anandkumar2020neural, patel2021physics, li2020fourier}.  In \cite{lu2019deeponet} the authors propose a network architecture based on a theorem of  approximation by neural networks of infinite dimensional operators. Very recently, the authors of \cite{li2020fourier} developed a new neural network based on the Fourier transform. The latter has the advantage of being very simple to implement and computationally very cheap thanks to the Fast Fourier Transform.    
    
    \item Lastly, there are also methods based on a stochastic approach which uses the links between PDEs and stochastic processes (see \cite{blechschmidt2021three} and references therein for more details). 
     
\end{itemize}

The approach proposed in this paper falls in the second category of approaches, i.e., our neural networks approximate the action of semi-group operators for which only few information is known, they are fully data-driven and enough training data is available to get an accurate approximation.\\

\paragraph{\textbf{Outline of the paper}} The paper is organized as follows: we first present in Section \ref{sec:Methodology} the strategy we adopt for the construction of our numerical schemes based on neural networks. In particular, we introduce the different semigroups involved in our study and detail the whole training protocol starting with the choice of the training data and metrics. The different neural network architectures are described in Section \ref{sec:Methodology}, the architectures being inspired by some discretization schemes for the Allen-Cahn equation as in  \cite{Eyre1998, zaitzeff2021high}. We focus on two particular networks denoted as $\S^{\text{NN}}_{\theta,1}$ and $\S^{\text{NN}}_{\theta,2}$, and we address in Section \ref{sec:validation} their capacity to approximate the mean curvature flows of either oriented or non-orientable surfaces. In the oriented case, we provide some numerical experiments and we give quantitative error estimates to highlight the accuracy and the stability of the numerical schemes associated with our networks. In the non orientable case, both networks seem to learn the flow but a first validation shows that $\S^{\text{NN}}_{\theta,1}$ fails to maintain the evolution of a circle while $\S^{\text{NN}}_{\theta,2}$ succeeds perfectly. To test the reliability of $\S^{\text{NN}}_{\theta,2}$ to approximate the mean curvature flow, evolutions starting from different initial sets are shown in Section \ref{sec:validation}. We observe in particular that $\S^{\text{NN}}_{\theta,2}$ is able to handle correctly non-orientable surfaces, even with singularities, and the Herring's condition \cite{bretin2018multiphase} seems to be respected at points with triple junction. This is somewhat surprising because only evolving smooth sets are used to train our networks (e.g., circles in $2D$, spheres in $3D$). Finally, we propose various applications in Section \ref{sec:application} (multiphase mean curvature flows, Steiner trees, minimal surfaces) to bring out the versatility of our approach and to show that the schemes derived from our trained networks are sufficiently stable to be coupled with additional constraints such as volume conservation or inclusion constraints, and sufficiently stable to be extended to the multiphase case.

\section{Neural networks and phase field mean curvature flows}
\label{sec:Methodology}

A first possible way, which is not the one we will opt for, to associate neural networks and phase field mean curvature flows is to compute approximations to the solutions of the Allen-Cahn equation \eqref{eq:Allen-Cahn}  by training a neural network $\S^{\mathrm{NN}}_{\theta}$ (depending on a parameter vector $\theta$)
to reproduce the action of the Allen-Cahn semigroup $\S^{\text{AC}}_{\dt, \varepsilon}$ defined by
 $$  \S_{\delta_t,\varepsilon}^{\text{AC}}[u_{\varepsilon}(\cdot,t)] =  u_{\varepsilon}(\cdot,t + \delta_t),$$ where  $u_\varepsilon$ is solution to the Allen-Cahn equation \eqref{eq:Allen-Cahn}  \textcolor{black}{on a set $Q$ with  periodic boundary conditions and with constant parameters $\delta_t$ and $\varepsilon$.} \\

 From this neural network, 
 it is possible to derive the simple numerical scheme 
  $$ u^{n+1} = \S^{\mathrm{NN}}_{\theta}[u^{n}],$$
 where the iterate $u^n$ is expected to be a good approximation of $u_{\varepsilon}$ at time $t_n = n \delta$. However, 
 as explained previously, learning to compute  solutions 
 to the Allen-Cahn equation is not necessarily of great interest since extremely simple and robust numerical schemes for 
 computing these solutions already exist. 
 
 The idea we propose in this paper is rather to train a network $\S^{\mathrm{NN}}_{\theta}$ to approximate 
 the semigroup $\S^{q}_{\delta_t,\varepsilon}$ defined  by 
 $$ \S^{q}_{\delta_t,\varepsilon}[v_{\varepsilon}(\cdot,t)] =  v_{\varepsilon}(\cdot,t + \delta_t),$$ 
 where $v_{\varepsilon} =  q\left(\frac{d(x,\Omega(t)}{\varepsilon} \right)$ is an exact phase field representation of an exact mean curvature flow $t\mapsto\Omega(t)$.  This second approach is by nature more accurate than the above one since the Allen-Cahn equation is only an approximation to the mean curvature flow whereas $v_{\varepsilon}$ encodes exactly the flow. However, the convergence results of the Allen-Cahn equation to the mean curvature flow show that the two semigroups $\S^{\text{AC}}_{\delta_t, \varepsilon}$ and $\S^{q}_{\delta_t, \varepsilon}$ are very \textcolor{black}{close}. And since very efficient numerical schemes exist for the Allen-Cahn equation, it makes sense to take inspiration from their structures to design efficient networks.

In a second step, we will adapt our strategy to approximate the mean curvature flow of possibly non-orientable sets
by simply replacing the phase field profile $q$. More precisely, 
we will introduce a network $\S^{\mathrm{NN}}_{\theta}$ to learn an
approximation of the semigroup $\S^{q'}_{\delta_t,\varepsilon}$ defined by
$$ \S^{q'}_{\delta_t,\varepsilon}[w_{\varepsilon}(\cdot,t)] =  w_{\varepsilon}(\cdot,t + \delta_t),$$
where $w_{\varepsilon} =  q'\left(\frac{ \text{dist}(x,\Gamma(t))}{\varepsilon} \right)$  
with $t \mapsto \Gamma(t)$ the mean curvature flow of a possible non orientable set.

\subsection{Training database and loss function }~\\
 Before providing details about the structure of our neural networks, let us describe the data on which they will be trained and the training energy (the so-called {\it loss function} in the literature of neural networks). As before, we shall denote with $t\mapsto\Omega(t)$ the mean curvature flow associated with an initial smooth open set $\Omega(0)$, and with $t\mapsto\Gamma(t)$ the mean curvature flow associated with an initial, possibly non-orientable set $\Gamma(0)$.
 
 Recall that our idea is to obtain an approximation of the operator  $\S^{\varphi}_{\delta_t,\varepsilon}$,  where $\varphi = q$ or $\varphi = q'$,  by training
 a neural network $S^{\mathrm{NN}}_{\theta}$ on a suitable dataset.
  We expect that the numerical scheme $$  u^{n+1} = \S^{\mathrm{NN}}_{\theta}[u^{n}]   $$ 
 coupled  with the initial data
 $$
 u^{0}(x) = \begin{cases}
                q\left(\frac{d(x,\Omega(0))}{\varepsilon}\right) & \text{ if } \varphi = q, \\
              q'\left(\frac{\text{dist}(x,\Gamma(0))}{\varepsilon}\right) & \text{ if } \varphi =q', 
                 \end{cases}
$$
  will be a good approximation of the mean curvature flows starting from either $\Omega(0)$ or $\Gamma(0)$. \\
 
 The training of the network $S^{\mathrm{NN}}_{\theta}$ consists in a gradient descent for the parameter vector $\theta$ with respect to a loss energy defined from the exact mean curvature flows of various circles in 2D, of various spheres in 3D, etc. 
 
 Let us describe more precisely how the dataset and 
 the associated loss energy are constructed in dimension $d=2$, the construction being similar in higher dimension. 
 Recall that a circle of radius $R_0$ evolving under mean curvature flow remains 
 a circle of radius 
 $R(t) =  \sqrt{R_0^2-2t}$ which decreases until the extinction time $T_{R_0} = \frac{R_0^2}{2}$.  
 We select a finite family of radii $R_i$ for $i \in \{ 1,\cdots,N_{\text{train}}\}$ and we define the 
 training dataset as follows:

$$\left \{ \left(X_i,Y_i\right) \right \}_{ i \in \{1 \cdots N_{\text{train}} \}}= \left(\varphi_{R_i}, \varphi_{\sqrt{R_i^2 - 2 \delta_t}} \right)_{ i \in \{1 \cdots N_{\text{train}} \}}$$
where  
$$
\varphi_{R}(x) = \begin{cases}
                q\left(\frac{\text{d}(x,B_R)}{\varepsilon}\right) & \text{ if } \varphi = q,\\
              q'\left(\frac{\text{dist}(x,\partial B_R)}{\varepsilon}\right) & \text{ if } \varphi =q'. 
                 \end{cases}
$$
with $B_{R}$ the ball of radius $R$ centred at $0$.

We introduce a first loss function:
$$ J_1(\theta) =  \frac{1}{ N_{\text{train}}} \sum_{i=1}^{  N_{\text{train}} } \| \S^{\mathrm{NN}}_{\theta}[X_i] - Y_i  \|^2 =   \frac{1}{ N_{\text{train}}} \sum_{i=1}^{  N_{\text{train}} } 
\int_Q \left( \S^{\mathrm{NN}}_{\theta}[\varphi_{R_i} ] -  \varphi_{\sqrt{R_i^2 - 2 \delta_t}} \right)^2 dx.$$

To stabilize the training of the network $S^{\mathrm{NN}}_{\theta}$, one can opt for a multipoint version which  consists in introducing the enriched data set
$$\left \{ \left(X_i,Y_{i,j}\right) \right \}_{\{ i \in \{1, \cdots, N_{\text{train}}\}, \, j \in \{1, \cdots, k   \} \}}= \left(\varphi_{R_i}, \varphi_{\sqrt{R_i^2 - 2 j \delta_t}} \right)_{ \{i \in \{1, \cdots, N_{\text{train}}\}, \, j \in \{1, \cdots, k   \} \}}$$
and in minimizing the loss functional $J_k$ defined by
$$ J_k(\theta) =  \frac{1}{ N_{\text{train}}} \sum_{i=1}^{  N_{\text{train}} } \sum_{j=1}^{k} \| (\S^{\mathrm{NN}}_{\theta})^{(j)}[X_i] - Y_{i,j}  \|^2 =   \frac{1}{ N_{\text{train}}} \sum_{i=1}^{  N_{\text{train}} } \sum_{j=1}^{k}
\int_Q \left( (\S^{\mathrm{NN}}_{\theta})^{(j)}[\varphi_{R_i} ] -  \varphi_{\sqrt{R_i^2 - 2 j \delta_t}} \right)^2 dx,$$
where  $(\S^{\mathrm{NN}}_{\theta})^{(j)}$ is the iterated $j$ times composition of $S^{\mathrm{NN}}_{\theta}$.\\

In all our experiments, the computational domain is $[0, 1]^d$ and 2D/3D training data consists of an average of $N_{\text{train}} = 100$ circles/spheres with a range of
radii ranging from $0.05$ to $0.45$ to capture enough information about the curvatures.

\subsection{From the numerical approximation of the Allen-Cahn semigroup to the structure of neural networks} We introduce in this section efficient neural network structures to approximate the semigroups $\S_{\dt, \varepsilon}^q$ and $\S_{\dt, \varepsilon}^{q'}$.
The design of a network architecture is a very delicate question because there is no single generic choice 
of network that can accurately approximate an operator. The approach developed in this paper consists in deriving
the architecture of our networks from splitting methods, which are often used to solve numerically evolution equations having a gradient flow
structure. Doing so, we draw on the rich knowledge of numerical analysis to guide us in designing new and more efficient networks. \\

To derive the structure of our networks, let us start by recalling the principle of splitting schemes \cite{LeeLee2014} to approximate the solutions of the Allen-Cahn equation defined in the domain $Q = [0, 1]^d$ with periodic boundary conditions (for a recent review of numerical methods for phase field approximation of various geometric flows see \cite{DuFeng2020}). As the two semigroups $S^{\text{AC}}_{\delta_t,\varepsilon}$ 
and $\S_{\dt, \varepsilon}^{q}$ are  closely related in the case of smooth mean curvature motion, 
our expectation is that imitating splitting schemes will lead us to very effective  
networks to approximate first $\S_{\dt, \varepsilon}^{q}$,  but also  $\S_{\dt, \varepsilon}^{q'}$.

Given a time step $\dt$, we construct an approximation sequence $(u^n)_{n\geq 0}$ of the solution $u_\varepsilon$ of \eqref{eq:Allen-Cahn} at time $n\dt$ using various splitting approaches. \\ 

\noindent {\bf {First-order neural network $\S^{\text{NN}}_{\theta,1}$}}\\
The first splitting method is the semi-implicit approach where the sequence $(u^n)$ is defined recursively from
$$  \frac{u^{n+1} - u^{n}}{\delta_t} =  \Delta u^{n+1} - \frac{1}{\varepsilon^2} W'(u^n),$$
i.e., 
\begin{equation}\label{scheme:semi-implicit}
 u^{n+1} = \left( I_d - \delta_t \Delta  \right)^{-1} \left( u^n - \frac{\delta_t}{\varepsilon^2} W'(u^{n})\right).
\end{equation}
More precisely, this numerical scheme can be written as a combination of a convolution kernel $K_1$ and an activation function $\rho_1$:
$$   u^{n+1} =   K_1 * \rho_1(u^{n}),$$
with $\rho_1(s) = s - \frac{\delta_t}{\varepsilon^2} W'(s)$ and
$$
K_1(x) = \mathcal{F}\left[\xi \mapsto \frac{1}{ 1+ \delta_t 4 \pi^2 |\xi|^2}\right](x),
 $$
 where $\mathcal{F}$ denotes the Fourier transform. This scheme is stable as soon as $\displaystyle\dt < \sup_{s \in [0, 1]} \abs{W''(s)} \varepsilon^2$. Moreover, the operator $ S^{\text{AC}}_{\delta_t,\varepsilon,1}: u \mapsto  K_1 * \rho_1(u)$ that encodes the
semi-implicit scheme can be viewed as an approximation of order $1$ of the Allen-Cahn semigroup 
$S^{\text{AC}}_{\delta_t,\varepsilon}$. This method has also the advantage of decoupling the action of 
the diffusion operator and the reaction operator, therefore each operator can be handled independently. \\

  The previous stability constraint can be avoided using a convex-concave splitting of the Cahn-Hilliard energy. 
  Following the idea introduced by Eyre \cite{Eyre1998}, the functional $P_\varepsilon = E_1 + E_2$ is 
  decomposed as the sum of a convex energy $E_1$ and a concave energy $E_2$ defined by
$$  E_1(u) =   \int_Q \left(\varepsilon \frac{|\nabla|^2}{2} +  \frac{\alpha}{\varepsilon}   \frac{u^2}{2} \right) dx  
 \quad\text{ and } \quad E_2(u) = \frac{1}{\varepsilon} \int_Q \left(  W(u) - \alpha\frac{u^2}{2} \right) dx,$$
with $\alpha$ a sufficiently large stabilization parameter. Treating the convex energy implicitly and the concave energy explicitly yields the scheme
$$
u^{n+1} = u^n - \frac{\dt}{\varepsilon} \left( \nabla E_1(u^{n+1}) + \nabla E_2(u^n)\right),
$$
i.e., 
\begin{equation}\label{scheme:convex-concave}
u^{n+1} = \left( I_d - \delta_t \left(\Delta - \frac{\alpha}{\varepsilon^2} I_d\right)   \right)^{-1}\left( u^n - \frac{\delta_t}{\varepsilon^2} (W'(u^{n}) - \alpha u^n)\right) = K_2 * \rho_2(u^{n}).
\end{equation}
Again, this numerical scheme is of the form
$$u^{n+1} = K_2 * \rho_2(u^{n})$$ 
where the convolution operator $K_2$ is now given by 
 $$
 K_2(x) = \mathcal{F}\left[ \xi \mapsto \frac{1}{ 1+ \delta_t \left(4 \pi^2 |\xi|^2 + \frac{\alpha}{\varepsilon^2}\right)}\right](x)
 $$ 
 and the activation function $\rho_2$ satisfies $\rho_2(s) = s - \frac{\delta_t}{\varepsilon^2}\left(W'(s) - \alpha s\right)$.\\
 
 The advantage of this scheme is to be unconditionally stable -- in the sense that  it decreases the Cahn-Hilliard energy -- as soon as the stabilization parameter $\alpha$ satisfies $\alpha > \sup_{s\in [0, 1]} \abs{W''(s)} = 1$ (to be complete, note that choosing a large value for $\alpha$ has also a bad influence on the dynamics of the flow).  It is also an approximation of order $1$ to the Allen-Cahn semigroup $\S_{\dt, \varepsilon}^{\mathrm{AC}}$.\\

\begin{figure}[htbp] 
    \centering
    \includegraphics[width=0.35\textwidth]{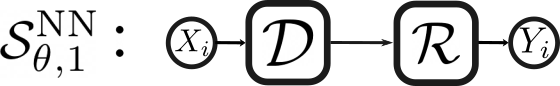}
    \caption{Definition of $\S^{\text{NN}}_{\theta,1}$ where $\D$ represents a convolution neuron (or diffusion neuron in reference to diffusion operators), and $\Re$ represents a multilayer perceptron (called reaction network in reference to reaction operators). The parameter $\theta$ is the vector of learning parameters, it contains all weights of the $\D$ and $\Re$ networks.}
    \label{fig:DRNN}
\end{figure}

 This stabilization effect leads us to believe that other such schemes could be obtained by directly
 learning the diffusion kernel $K$ and the activation function $\rho$ from the training data. 
The extension to networks is then straightforward since the latter operations can be interpreted as networks. 
Indeed, the schemes \eqref{scheme:semi-implicit} and \eqref{scheme:convex-concave}  are reminiscent of the structure
of convolutional neural networks (CNN)~\cite{lecun1989handwritten, mallat2016understanding}:
a convolution operation coupled with a nonlinear activation function, where the nonlinearity
is given by $\rho_1$ (or $\rho_2$ for \eqref{scheme:convex-concave}) and the learning parameters are the parameters of the kernel
associated with the convolution. A first idea of neural network would be to use a CNN with $\rho_1$ (for instance) as 
nonlinearity but it is restrictive because of the choice of $\rho_1$ that depends on the parameters $\dt$ and $\varepsilon$. 
We propose instead a neural network  $\S^{\text{NN}}_{\theta,1}$ constructed as the composition of a pure convolution neuron and a multilayer perceptron (MLP) \cite{rosenblatt1958perceptron} that will act as a nonlinearity. 
The network  $\S^{\text{NN}}_{\theta,1}$ we propose can be represented by the diagram pictured in figure~\ref{fig:DRNN}.

The choice of representing the reaction network by a MLP comes from the fact that the action of the
reaction operator is totally encoded by a 1D function. According to the fundamental neural network approximation theorem \cite{hornik1991approximation, hornik1993some, scarselli1998universal} which states that any function can be approximated by a two-layer MLP, it is quite natural to consider a MLP to represent the action of the reaction operator. However, we also choose to use a multilayer perceptron (MLP) as a non-linearity in order to keep some flexibility and not to depend on the context in which we work. The characteristics of the non-linearity are therefore entirely determined during the training stage.\\

\begin{rem}
Note that in contrast to the splitting schemes \eqref{scheme:convex-concave} and \eqref{scheme:semi-implicit}, we have chosen to define the network $\S^{\text{NN}}_{\theta,1}$ by starting with a diffusion neuron $\D$ followed by a non-linearity $\Re$ to have a structure similar to the usual convolutional network. The other architecture where we start with $\Re$ followed by $\D$ is also very legitimate and gives training results which are similar to those obtained with $\S^{\text{NN}}_{\theta,1}$. However, better numerical results are observed in the test phase when we use $\S^{\text{NN}}_{\theta,1}$ so we shall study only this network.
\end{rem} 

 \noindent{\bf {Higher-order neural networks $\S^{\text{NN}}_{\theta,2}$ and $\S^{\text{NN}}_{\theta,3}$}}\\ 
 The network $\S^{\text{NN}}_{\theta,1}$ provides us with the \textbf{basic block} to design more complex and  
 deeper networks as is often the case in deep learning. One can mention, for instance, deep convolutional neural networks 
(where the basic block is the convolution operation) and deep residual neural networks (with the residual block). Many different
networks can derive from $\S^{\text{NN}}_{\theta,1}$. 
 In particular, in order to obtain more complex and effective networks, we now try to adapt the structure of these networks 
 from the higher order semi-implicit schemes developed in \cite{zaitzeff2021high} in the case of a gradient flow structure.  
 These schemes are structured as follows: starting from $U_0 = u^n$, we set $u^{n+1} = U_M$ where $U_M$ is defined recursively from $U_0$ as 
  
  \begin{equation}\label{scheme:higher_order_semi_implicit}
  U_m = U_0 - \frac{\dt}{\varepsilon} \left(\sum_{i=0}^m \gamma_{i,m} \nabla E_1(U_i) +  \sum_{i=0}^{m-1} \tilde{\gamma}_{i,m} \nabla E_2(U_i)\right) \textrm{ for } m = 1, \dots, M
  \end{equation}
  where the parameters $\gamma_{i,m}$ and $\tilde{\gamma}_{i,m} $ are predefined to ensure the numerical scheme to be highly accurate, see \cite{zaitzeff2021high} for more information about the choice of the coefficients $\gamma_{i,m}$ and $\tilde{\gamma}_{i,m}$.
  As can be observed in this scheme \eqref{scheme:higher_order_semi_implicit}, the diffusion and reaction operators are applied in chain but
  also in parallel in order to keep the information of each $U_i$ for $i=1, \ldots, M-1$.
  
\begin{figure}[htbp]
    \centering
    \includegraphics[width=0.6\textwidth]{{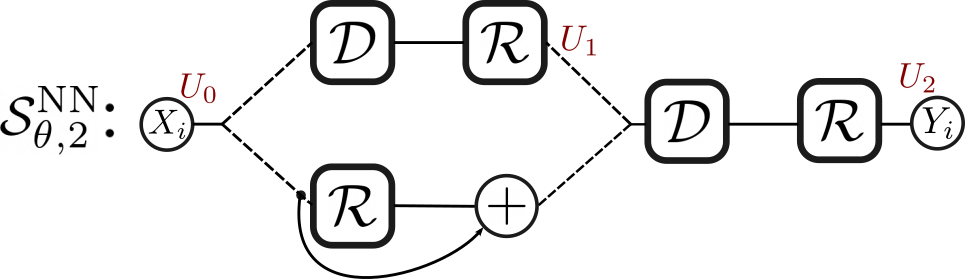}}
    \caption{Representation of the 2nd order neural network DR $\S^{\text{NN}}_{\theta,2}$. We adopt the same formalism as 
    for MLPs by assimilating the $\D$ and $\Re$ networks to neurons. Here all blocks $\D$ and $\Re$ are different, independent networks.
    The dotted lines mean that the output of the source neuron is multiplied by a weight. 
    The curved line corresponds to a skip connection and indicates that the network $\Re$ is equipped with a residual structure.}
    \label{fig:2Nd_order_NN}
\end{figure}

  The second-order network $\S^{\text{NN}}_{\theta,2}$ which is shown in~figure~\ref{fig:2Nd_order_NN} is inspired by the scheme \eqref{scheme:higher_order_semi_implicit} with $M=2$, which reads
   $$ U_2 = K_2*[ \rho_3(U_1) + [U_0 + \rho_2(U_1)],$$
   with $U_1 = K_1*[\rho_1(U_0)]$, and  $\gamma_{0,1} = \gamma_{0,2} = \gamma_{0,1} = 0$. \textcolor{black}{The curved line in ~figure~\ref{fig:2Nd_order_NN}, ended by the symbol $\bigoplus$, indicates that the neuron/network $\Re$ is equipped with a residual structure in the sense that each entry $u$ to the network $\Re$ is added to the output $\Re(u)$, i.e. the output of the bottom part of the architecture is $u + \Re(u)$. The curved line is actually a skip connection, a technique frequently used in machine learning which consists in skipping one or several layers of a neural network to improve the convergence of the training process by diminishing the effects of vanishing gradients.}

\begin{rem}
 The residual structure is also a consequence of different experiments we have done where networks with 
 this type of structure had better learning scores and the learning process was faster. In the particular case 
 of the network $\S^{\text{NN}}_{\theta,2}$  we can make the following conjecture: in the classical 
 schemes \eqref{scheme:semi-implicit} and \eqref{scheme:convex-concave} reaction operators are related to the double-well 
 potential but also to the phase field profile. In the case of the network $\S^{\text{NN}}_{\theta,2}$, the residual structure allows 
 to maintain this profile after each iteration. Indeed, several recent works \cite{greff2016highway} highlight the particularity of
 residual networks to stay close to the inputs by adding the identity. 
\end{rem}

\begin{figure}[htbp]
    \centering
    \includegraphics[width=0.8\textwidth]{{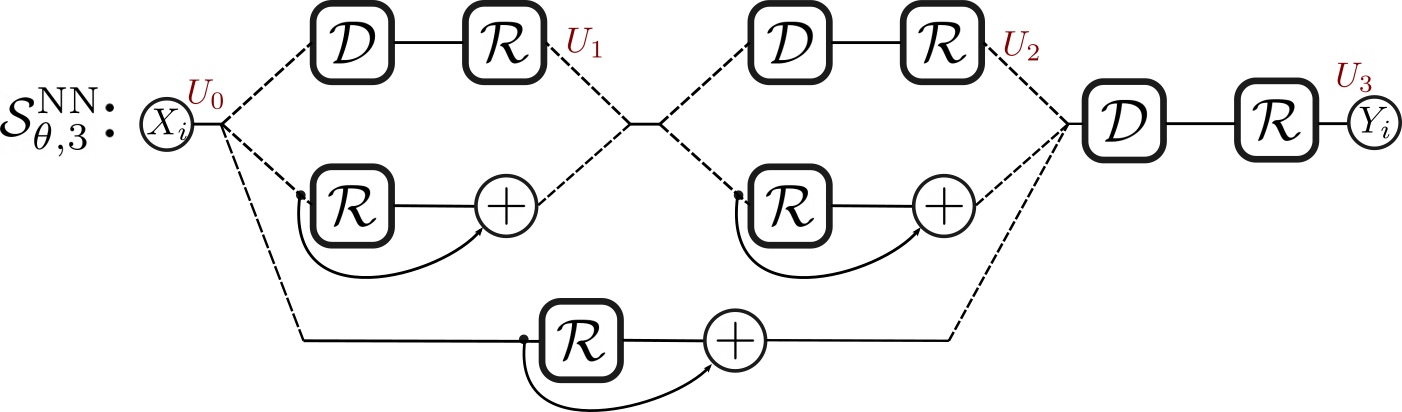}}
    \caption{Representation of the 3rd order neural network DR $\S^{\text{NN}}_{\theta,3}$. The same idea used for the structure of the network $\S^{\text{NN}}_{\theta,2}$ is applied to design $\S^{\text{NN}}_{\theta,3}$.}
    \label{fig:third_order_NN}
\end{figure}

We can go further in the complexity and the depth of our networks by taking inspiration again from the structure of the
scheme \eqref{scheme:higher_order_semi_implicit} with $M = 3$. For instance, following the principle of the previous network,
we obtain the architecture shown in figure \ref{fig:third_order_NN}.

In this paper, we will limit ourselves to the first two networks $\S^{\text{NN}}_{\theta,1}$ and $\S^{\text{NN}}_{\theta,2}$ for 
the reason that they already give very satisfactory results and it appeared more relevant to us to focus on these first two structures.

\subsection{Discretization and specification of our neural networks}

In practice, we will consider some approximations of phase field mean curvature flow  defined  on a square-box $Q = [0, 1]^d$
using Cartesian grid with $N$ nodes in each direction and using periodic boundary conditions.  We consider 
the approximation parameter $\varepsilon= 2/N$ and the fixed timestep $\delta_t= \varepsilon^2$. 
Our neural networks are applied to discrete images $U^{n}_{i,j}$ which correspond to a sampling
of the function $u^{n}(x_{i,j})$ at the points $x_{i,j} = ((i-1)/N,(j-1)/N)$. Our networks depend  on the specific 
choice of the parameters   $(\delta_x, \varepsilon, \delta_t)$  which have been fixed before the training procedure.
In particular, in all the numerical experiments presented below, we set  $\delta_x = 1/N =  1/2^8$,  $\varepsilon = 2 \delta_x$
and  $\delta_t = \varepsilon^2$.

\begin{rem}
Our approach does not correspond exactly to a method based on neural operators
since the definition of the convolutional kernel depends on the discretization in space. It will therefore not be 
possible to use the trained networks for other choices of discretization parameters. Nevertheless, 
it is possible to use the trained networks on different sizes of computational box as long as the value 
of the parameter $\delta_x$ remains unchanged. For instance, 
we can consider a square-box $Q = [0, 2]^d$ discretized with $N = 2^9$ nodes in each direction. In this sense, our method 
remains resolution-invariant thanks to multi-scale techniques.
\end{rem}

\noindent {\bf {Diffusional network based on a discrete kernel convolution.}}~\\
As $\D$ corresponds to a convolution operation, one needs to specify its hyper-parameters, especially the kernel size $N_K$ which is 
related to the domain discretization and the timestep $\dt$. In practice, we consider
a discrete square kernel $\K$ of size $N_K = 17$ where the kernel convolution reads as
$$ \D[u^{N}]_{\bf k} =  (\K*u^N)_{\boldsymbol k} =  \sum_{ {\boldsymbol \ell}\in L_N  }  \K_{\boldsymbol \ell } \; u^{N}_{\boldsymbol k - \boldsymbol \ell}   $$
 with $  L_N =  [-(N_K-1)/2,(N_K-1)/2 ]^d$. Here, we assume that the padding extension of $u^N$ is periodic.
 Moreover, the convolution product is computed in practice using the Fast Fourier Transform.  
 Here,  recall that the Fourier $\boldsymbol K$-approximation of a function $u$ defined in a box 
$Q = [0,L]^d$ is given by
$$u^{\boldsymbol N}(x) = \sum_{{\boldsymbol k}\in K_N } c_{\boldsymbol k} e^{2i\pi{\boldsymbol \xi}_k\cdot x},$$
where  $K_N =  [-\frac{N}{2},\frac{N}{2}-1 ]^d$,   ${\boldsymbol k} = (k_1,\dots,k_d)$ and ${\boldsymbol \xi_k} = (k_1/L,\dots,k_d/L)$.
In this formula, the $c_{\boldsymbol k}$'s denote the $K^d$ first discrete Fourier coefficients of $u$. 
The inverse discrete Fourier transform leads to 
$u^{N}_{\boldsymbol k} =   \textrm{IFFT}[c_{\boldsymbol k}]$ 
where $u^{N}_{\boldsymbol k}$ denotes the value of $u$ at the points 
$x_{\boldsymbol k} = (k_1 h, \cdots, k_d h)$ and where $h = L/N$. Conversely,
$c_{\boldsymbol k}$ can be computed as the discrete Fourier transform of $u^K_{\boldsymbol k},$ {\em i.e.}, $c_{\boldsymbol k} =
\textrm{FFT}[u^N_{\boldsymbol k}].$\\

\begin{rem}
In practice, it can be interesting to use kernels with sufficiently small  size to facilitate their training but 
it is necessary to choose a size large enough to avoid anisotropy phenomena due
to the square shape of the kernels' domain.
\end{rem}

\begin{figure}[htbp]
    \centering
    \includegraphics[width=0.25\textwidth]{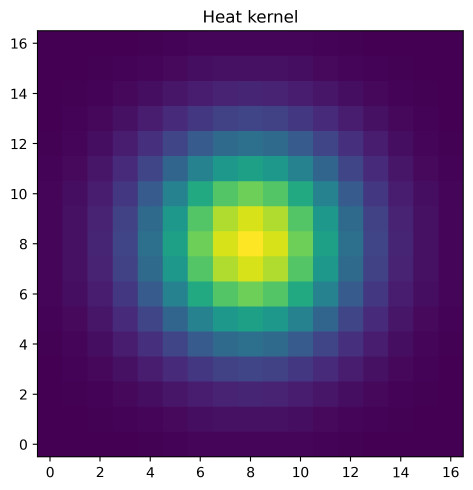}  $\quad \quad $
    \includegraphics[width=0.25\textwidth]{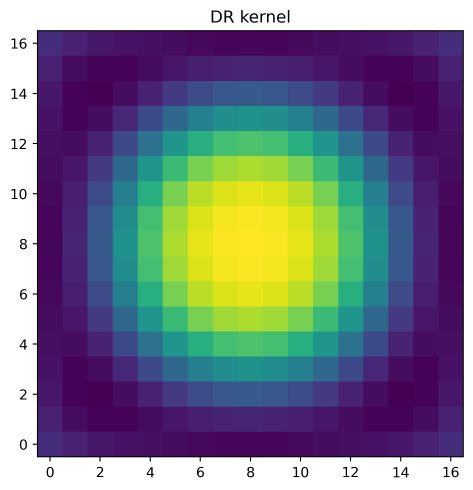}
    \caption{Discrete kernel $\K$. Left: sampling of the heat kernel; Right: example of learned kernel $\K$ for the network $\S^{\text{NN}}_{\theta,1}$ trained to approximate the semigroup $\S_{\varepsilon,\delta_t}^q$.}
    \label{fig:Kernel}
\end{figure}

\begin{rem}
As we have seen previously, the characteristic size of $\K$ depends strongly on the choice of the time step $\delta_t$.
In our case, using the time step $\delta_t = \epsilon^2 = (2/N)^2$, as shown in figure~\ref{fig:Kernel},
the support of the heat kernel $K_{\delta_t}$ is contained in the interval  $[-8/N,8/N]$  which suggests 
a good  approximation by a discrete kernel of size $N_K = 17$. On the other hand, a time step four times larger $\delta_t = 4 \epsilon^2$
would have required a kernel twice as large, i.e., $N_K = 33$.
\end{rem}

\begin{figure}[htbp]
    \centering
    \includegraphics[width=0.6\textwidth]{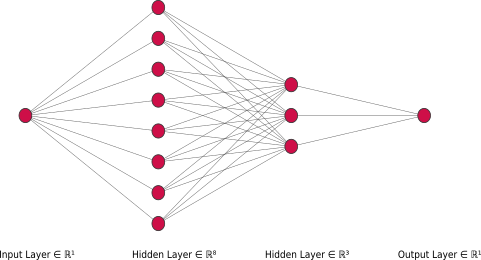} 
     \caption{Reaction network  based on a $1D$ multilayer perceptron  using two hidden layers with respectively $8$ and $3$ neurons.}
    \label{Rection_network}
\end{figure}

\noindent {\bf {Reaction network  based on a $1D$ multilayer perceptron}} ~\\
The reaction network is a multilayer perceptron acting as a function from $\R$ to $\R$.
In practice, the number of hidden layers is set to $2$ with $8$ neurons on the first layer and $3$ neurons on the second layer, as pictured in figure~\ref{Rection_network}. By a slight abuse of notation, we denote the same way this network from $\R$ to $\R$ and its action on an image, i.e. we consider that the reaction network applies independently to each component of $u^N$:
$$ \forall  k \in K_N, \qquad (\mathcal{R}[u^{N}])_{\bf k}=\mathcal{R}[u^{N}_{\bf k}]=  \rho{[2]}( W^{[3]} (\rho{[2]}( W^{[2]}  (  \rho{[1]}( W^{[1]}  ( u^N_{\bf k} )   +   b^{[1]}  )   )   +   b^{[2]})) + b^{[3]}),$$
\textcolor{black}{where  $W^{[1]},b^{[1]},W^{[2]},b^{[2]},W^{[3]},b^{[3]}$  denote the parameters of the network to be optimized and
$\rho{[1]},\rho{[2]},\rho{[3]}$ are the three activation functions. Thus, although the reaction neuron is a 1D multilayer perceptron, it extends to an operator that transforms an image into another image of the same size.}

We also select the appropriate activation function in the reaction network by testing different
choices of activation functions (ReLU, ELU, Sigmoid, SiLU, Tanh, Gaussian) in our numerical experiments. 
In the end, the Gaussian model seems to have the best properties by presenting more efficient learning rates and speeds.

\subsection{Neural network optimization, stopping criteria and Pytorch environment} ~\\
We now turn to the question of how the learning parameters of our networks are estimated.
The learning procedure consists in seeking an optimal vector parameter $\theta$ minimizing 
the energy $J_1$ (or $J_k,~k>1$ depending on the problem) using a mini-batch stochastic gradient descent algorithm 
with adaptive momentum \cite{kingma2014adam, bottou2010large, masters2018revisiting} on a set of training data.

\textcolor{black}{The reason we use a mini-batch approach} is that our problem is a high-dimensional non-convex optimization problem.
It is well known that optimizing on mini-batches rather than on the whole dataset reduces the computational cost and 
the memory cost, and it is commonly used to escape saddle points of the non-convex energy $J_1$ (see \cite{ge2015escaping} on this subject).

More precisely, the training strategy is as follows:

\begin{description}
\item[Step 1] \textbf{Define the neural network } $\S_\theta$ with the training vector parameter $\theta$.
\item[Step 2] \textbf{Generate data using mini-batch:} Randomly sample a mini-batch of size $B$ of training labeled
couples $(X_{i_k}, Y_{i_k}),~k = 1, \cdots, B$ over the training dataset $\{(X_i, Y_i)\}_{i=1,\cdots, N_{\mathrm{train}}}$. 
This means that the training dataset is organized into $N_{\mathrm{train}}/B$ mini-batches and the parameter $\theta$ 
will be updated for each of these batches. The training \textbf{batch size} $B$ is set to $B=10$ over all experiments. \textcolor{black}{The dataset is reshuffled after every epoch, so that the collection of mini-batches is different at every epoch.}
\item[Step 3] \textbf{Compute the loss} of the mini-batch 
$J_1(\theta) =\displaystyle \frac{1}{B} \sum_{k=1}^B \left\|\S_{\theta}[X_{i_k}] - Y_{i_k}\right\|^2$.
\item[Step 4] \textbf{Compute the gradient} of the previous loss using \textbf{back-propagation}. 
\item[Step 5] \textbf{Update the learning parameter} $\theta$ using the stochastic gradient descent algorithm with adaptive  
momentum (\textbf{Adam} \cite{kingma2014adam} in our case) and using with the gradient computed at the previous step. 
\item[Step 6] Repeat the steps $2$ to $5$ until all mini-batches have been used. Once all the mini batches have 
been used, an \textbf{epoch} is said to have passed. In practice, the epoch is set to $\mathrm{epoch} = 20000$ for all experiments.
\item[Step 7] Repeat steps $2$ to $6$ until a \textbf{stopping criterion} is satisfied.
\end{description}

\begin{rem}
 The number of learning parameters of each network is rather small with 
$N_\theta = 336$ for the first-order network $\S^{\text{NN}}_{\theta,1}$ and
 $N_{\theta} = 724$ for the second-order network $ \S^{\text{NN}}_{\theta,2}$.
\end{rem}

\begin{rem}
Very roughly speaking, stochastic gradient descent algorithms are parametrized by the learning rate, i.e. the step parameter just 
next to the gradient. The choice of the learning rate is a challenge. Indeed, a too small value may lead to a long 
training process that could get stuck, while a value that is too large may lead to a selection of sub-optimal 
learning parameters or an unstable learning process. In practice, schedulers are often used to gradually adjust the 
learning rate during training. In our case we use the Adam algorithm in step 5 which has the advantage of adjusting 
the learning rate according to the obtained value of the gradient in step 4.  
\end{rem}

\begin{rem} \textcolor{black}{Regarding the numerical evaluation of the loss function} 
$$ \textcolor{black}{J_1(\theta) =  \frac{1}{ N_{\text{train}}} \sum_{i=1}^{  N_{\text{train}} } \| \S^{\mathrm{NN}}_{\theta}[X_i] - Y_i  \|^2 =   \frac{1}{ N_{\text{train}}} \sum_{i=1}^{  N_{\text{train}} } 
\int_Q \left( \S^{\mathrm{NN}}_{\theta}[\varphi_{R_i} ] -  \varphi_{\sqrt{R_i^2 - 2 \delta_t}} \right)^2 dx,}$$
\textcolor{black}{we observe that all quadrature methods have the same order in periodic boundary conditions, and we opt for the following discretization:}
$$\textcolor{black}{ 
J^{N}_1(\theta) =   \frac{1}{ N_{\text{train}}}  \delta_x^d \sum_{{\bf k} \in K_N}   \left(  \S^{\mathrm{NN}}_{\theta}[\varphi^{N}_{R_i} ]_{\bf k} 
-  [\varphi^{N}_{\sqrt{R_i^2 - 2 \delta_t}}]_{\bf k} \right)^2. }
$$
\end{rem}

\paragraph{\textcolor{black}{\textbf{Checkpointing and selecting criterion}}}~\\
 It may happen that after the training process the learning parameters are biased according to the training data. 
 To prevent this and to avoid under- or over-fitting, it is recommended to add an intermediate step during the training 
 procedure. This is called the validation step which consists in evaluating the model on new data (the so-called validation data)
 and then measuring the  error made on these data using one or more metrics. In our case, the validation step is applied at 
 regular intervals every 100 epochs. It is important to note that the parameters of the network are not updated during this procedure.

 This step is used as a stopping criterion with the following validation metric which measures the estimation error for the volume of a given sphere evolving by mean curvature. In some sense, this metric is also used to enforce the \textcolor{black}{robustness of the neural network to predict the correct output over several iterations}. \textcolor{black}{In practice, the networks parameters which are selected in the end among the parameters saved at each checkpoint are those for which
 the validation metric gives the best score. }

More precisely, the validation metric is given by 

$$
E_l(\theta) = \sum_{i=1}^l \sum_{n = 0}^{n_{\mathrm{max}}^i} \abs{\int_Q u_i^n - \int_Q \varphi_{\sqrt{R_i^2 - 2n\dt}}}^2,
$$

where the sequence $(u_i^n)$ is defined iteratively by 

 $$ \begin{cases}
     u_i^{n+1} &=  \S_{\theta}[u_i^{n}],  \\
      u_i^{0} &= \varphi_{R_i},
    \end{cases}
 $$
and $n_{\mathrm{max}}^i = \max \{n\in \N,~ n\dt < T_{R_i}\}$ with $T_{R_i} = \frac{R_i^2}{2}$ the extinction time of 
the mean curvature flow of the initial sphere of radius $R_i$,  $i=1,\ldots, l$.

 Note that this particular choice of metric stems from the following observation: the evolution of the volume of a sphere with initial radius $R_0$ evolving under mean curvature flow is explicitly 
given by $V(t) = \pi R(t)^2 =  \pi \left(R_0^2 - 2t \right)$.

Moreover, if $\varphi = q$,  as $\varphi_R \approx 1_{B_R}$ one can show that 

$$
\int_Q \varphi_{R(t)} = \mathrm{Vol}(B_{R(t)}) + \O(\varepsilon^2) = V(t) + \O(\varepsilon^2).
$$

In addition, it can also be shown with a similar argument that for $\varphi =-\frac{1}{\varepsilon} q'$,

$$
\int_Q \varphi_{R(t)} = \H^{d-1}(\partial B_{R(t)}) + \O(\varepsilon^2).
$$

In this way, one can see $t\mapsto\int_Q \varphi_{R(t)}$ as a good measure of the evolution of the volume (or the perimeter
when $\varphi =-\frac{1}{\varepsilon} q'$) of a sphere with initial radius $R_0$ that evolves by mean curvature flow. In this sense, the energy $E_{l}$ allows to preserve some stability over time (at least on circles or spheres). \\

All our networks, datasets, visualization tools, etc., have been developed in a Python module 
dedicated to machine learning
for the approximation of interfaces that evolve according to a geometric law (as the mean curvature flow) and are represented by a phase field. 
This module uses the machine learning framework PyTorch \cite{NEURIPS2019_9015}
and the high-level interface PyTorch Lightning \cite{falcon2019pytorch} as the machine learning back-end. For the sake of reproducibility, source codes will be made available at \href{https://github.com/PhaseFieldICJ}{https://github.com/PhaseFieldICJ}.

Training computation has been made on a computational server hosted at the Camille Jordan Institute
and using a general-purpose graphical processing unit Nvidia V100S.
Typical training (20000 epochs) used in the following numerical experiments requires less than 30 minutes of computation without the validation steps.

\section{Validation}
\label{sec:validation}

In this first numerical section, we test the ability of two neural networks to learn either the mean curvature
flow $t\mapsto\Omega(t)$ or oriented domains or the mean curvature flow $t\mapsto \Gamma(t)$ of a possible non-orientable sets. These two neural networks are:

\begin{itemize}
\item The 1st-order DR network $\S^{\text{NN}}_{\theta,1}$ which combines successively a diffusion neuron and a reaction network. 
\item The 2nd-order network $\S^{\text{NN}}_{\theta,2}$, which uses two diffusion neurons and three reaction networks, 
according to the architecture illustrated on figure~\ref{fig:2Nd_order_NN}.
\end{itemize}

More specifically, in each case we will propose a methodology to learn the action of the semigroups  $\S_{\delta_t,\varepsilon}^q$ and  $S_{\delta_t,\varepsilon}^{q'}$
where the model parameters are set to $\varepsilon = 2 \delta_x$ and $\delta_t = \varepsilon^2$. \\

First, the training of the two networks $\S^{\text{NN}}_{\theta,1}$ and $\S^{\text{NN}}_{\theta,2}$ allows us to obtain sufficiently accurate and effective approximations of the semigroup $\S_{\delta_t,\varepsilon,q}$ to get an approximation of the mean curvature flow of oriented domains whose quality is at least
comparable with what would be obtained by solving the Allen-Cahn equation. Going further, a numerical error analysis in the particular case of the evolution of a circle shows that both networks make it possible to reduce the phase field approximation error related to the convergence of the Allen-Cahn equation to the motion by mean curvature in $O(\varepsilon^2)$.\\

In the case of the mean curvature flow of possibly non-orientable sets, the first order DR network $\S^{\text{NN}}_{\theta,1}$ fails to learn a sufficiently accurate approximation of the semigroup $\S_{\delta_t,\varepsilon}^{q'}$ and reveals that the basic block architecture DR is not well adapted to this specific setting. On the other hand, the structure of the second order DR network seems much more suitable and provides an interesting  approximation of the semigroup $\S_{\delta_t,\varepsilon}^{q'}$ allowing to obtain, at least qualitatively, approximate evolutions of interface 
under mean curvature flow even in the case of a non-orientable interface $\Gamma$ with multiple junctions.

\subsection{Oriented mean curvature flow $t\mapsto \partial \Omega(t)$  and approximation of $S_{\delta_t,\varepsilon}^q$}\label{subsect:Oriented_MCF}
\begin{figure}[htbp]
    \centering
    \includegraphics[width=1\textwidth]{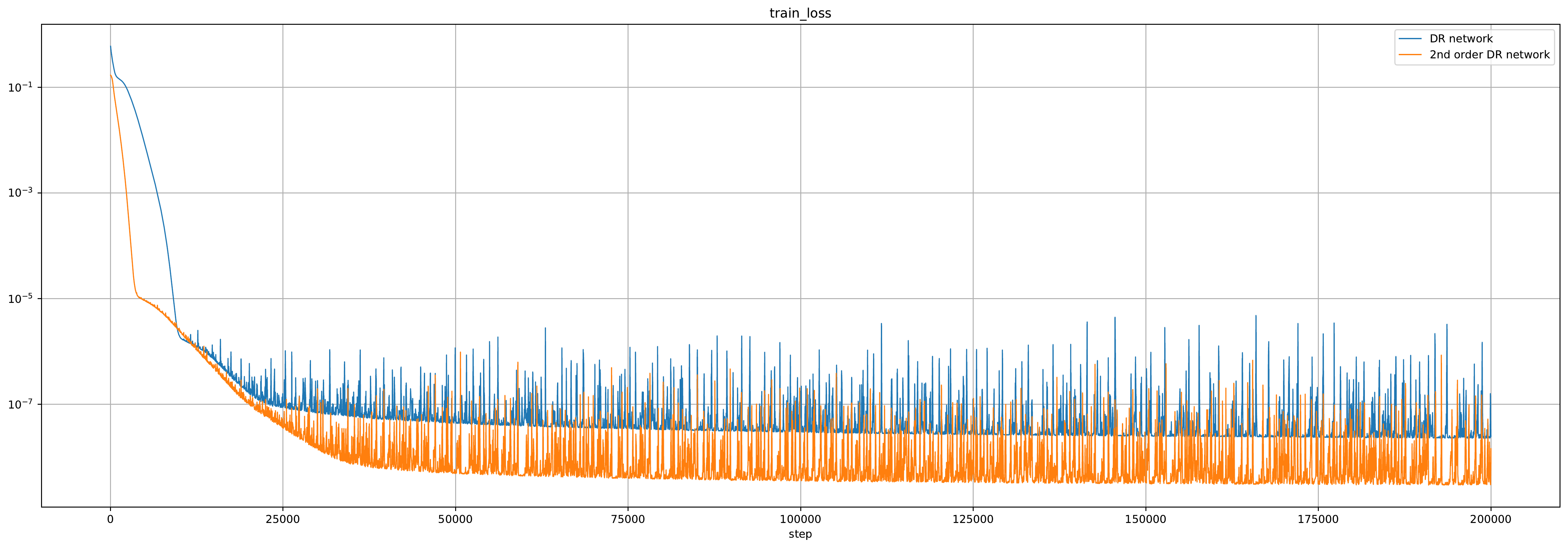}
    \caption{Learning procedure via optimization of the error $J_1$ over training data. Training losses of $\S^{\text{NN}}_{\theta,1}$ and $\S^{\text{NN}}_{\theta,2}$ plotted in blue and orange, respectively.}
    \label{fig:Learning_process_oriented}
\end{figure}

~\\
We now present some preliminary results about the approximation of the semigroup $S_{\delta_t=\varepsilon^2,\varepsilon}^q$ based on the training of both neural networks $\S^{\text{NN}}_{\theta,1}$ and $\S^{\text{NN}}_{\theta,2}$ on circles evolving by mean curvature flow and using exact oriented phase field representation.  We first plot on figure~\ref{fig:Learning_process_oriented} the evolution of the training loss energy 
throughout the  optimization process, in blue for the network $\S^{\text{NN}}_{\theta,1}$ 
and in orange for the second one $\S^{\text{NN}}_{\theta,2}$. Clearly, we observe that we manage to learn an approximation of $S_{\delta_t=\varepsilon^2,\varepsilon}^q$ with both networks. However, the second order DR network 
$^{\text{NN}}_{\theta,2}$ reaches more quickly acceptable errors on the learning loss,  
which suggests that the structure of this network is better suited  to our problem.

We now want to check that both trained networks achieve accurate approximations of
the mean curvature  flow.  To do so, we test the two associated numerical schemes 
$u^{n+1} = \S^{\text{NN}}_{\theta,\alpha}[u^{n}]$ with $\alpha \in \{1,2 \}$, coupled with  an initial condition $u^{0}$ given
by $u^{0}= q(d(\Omega_0,\cdot)/\varepsilon)$ where $\Omega_0$ corresponds to a disk of radius $R_0 = 0.3$. \\

\begin{figure}[htbp]
    \centering
    \includegraphics[width=0.2\textwidth]{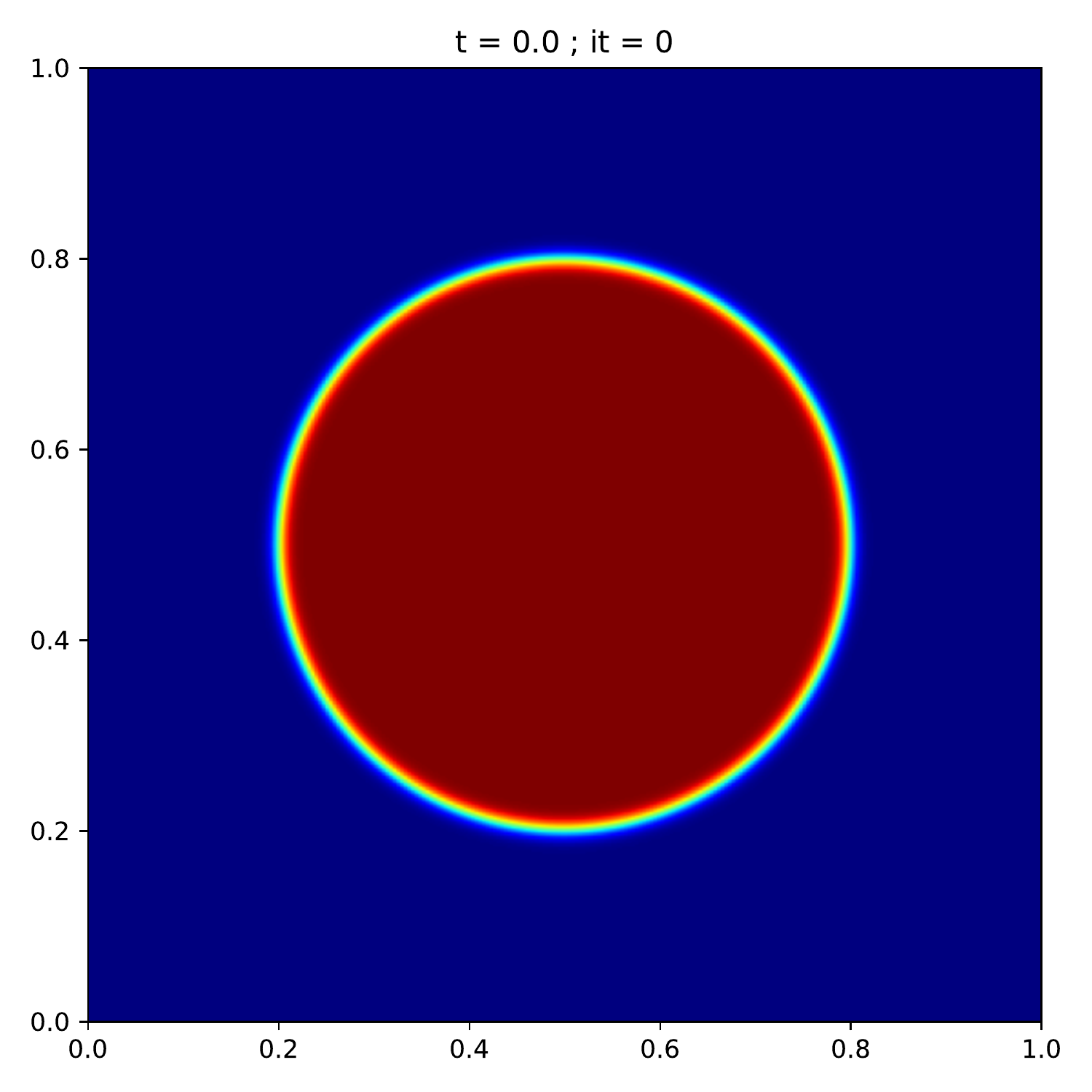}
    \includegraphics[width=0.2\textwidth]{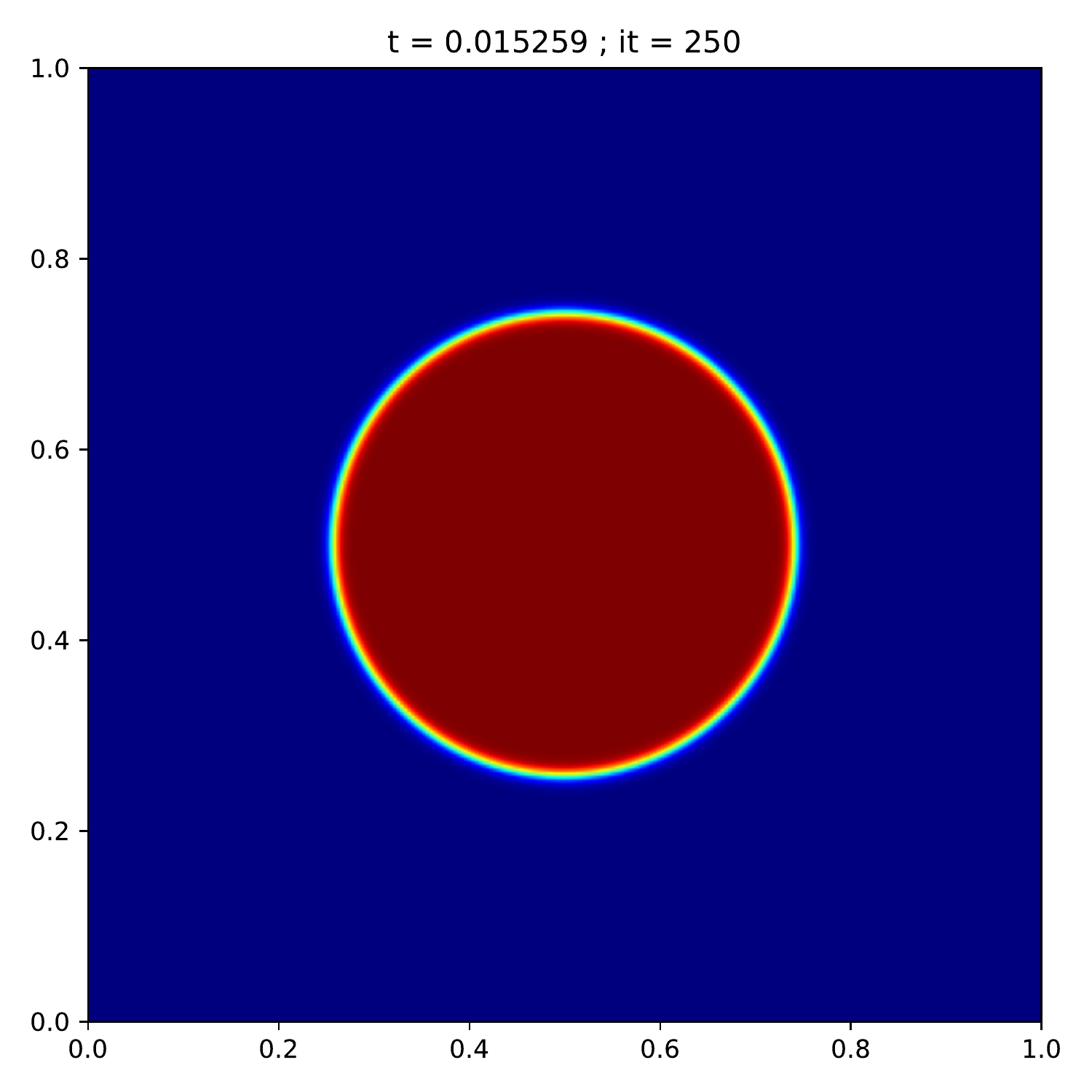}
    \includegraphics[width=0.2\textwidth]{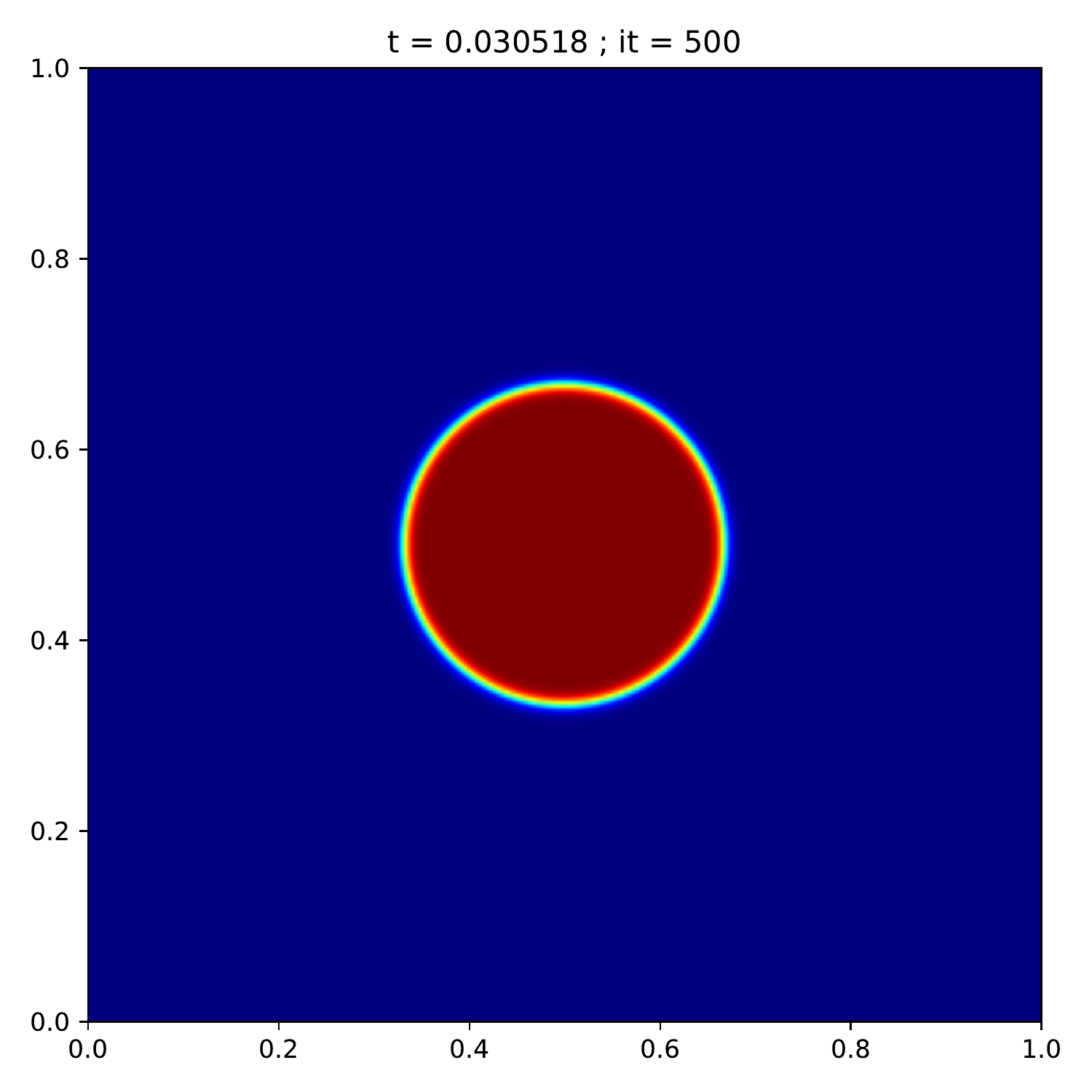}
    \includegraphics[width=0.2\textwidth]{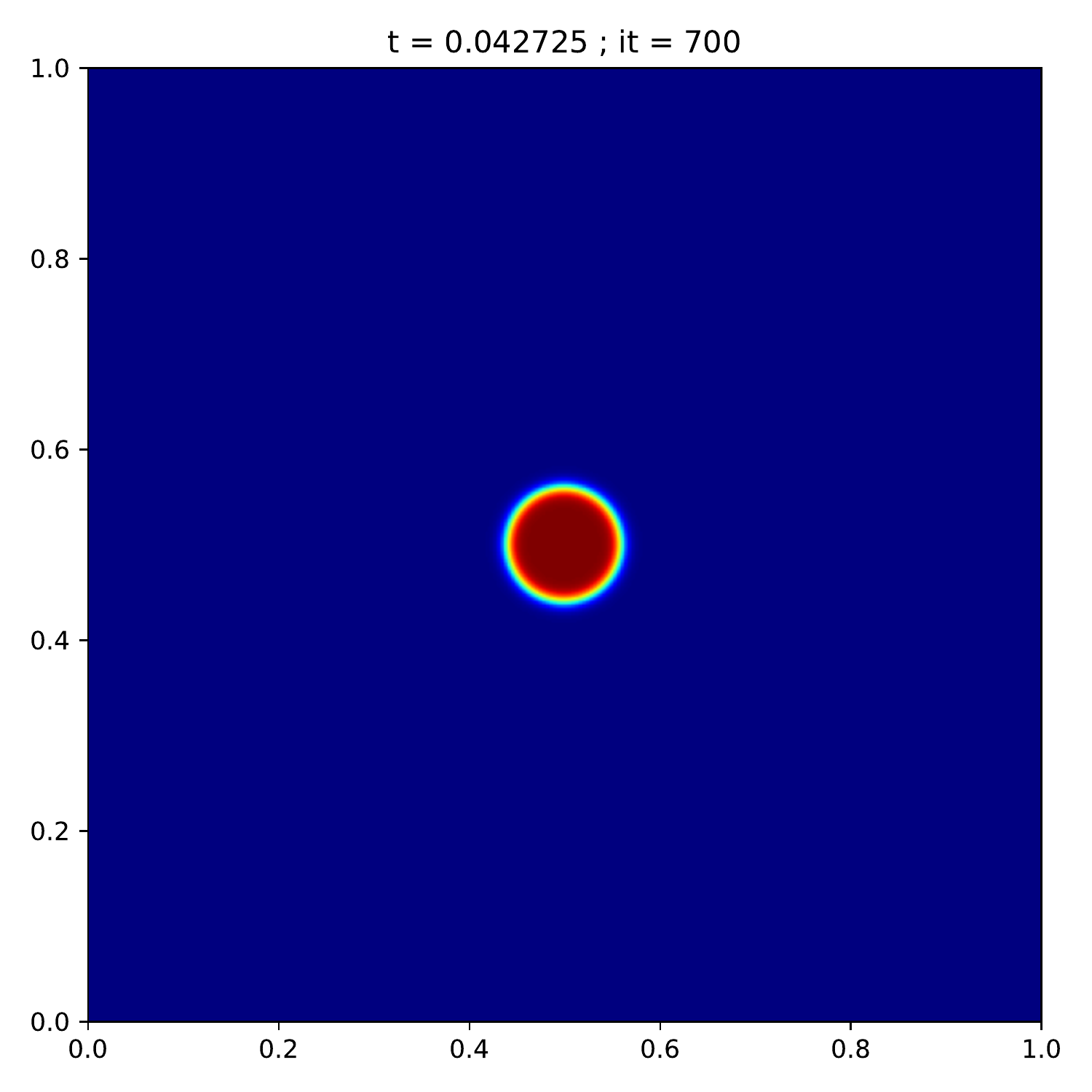} \\
      \includegraphics[width=0.2\textwidth]{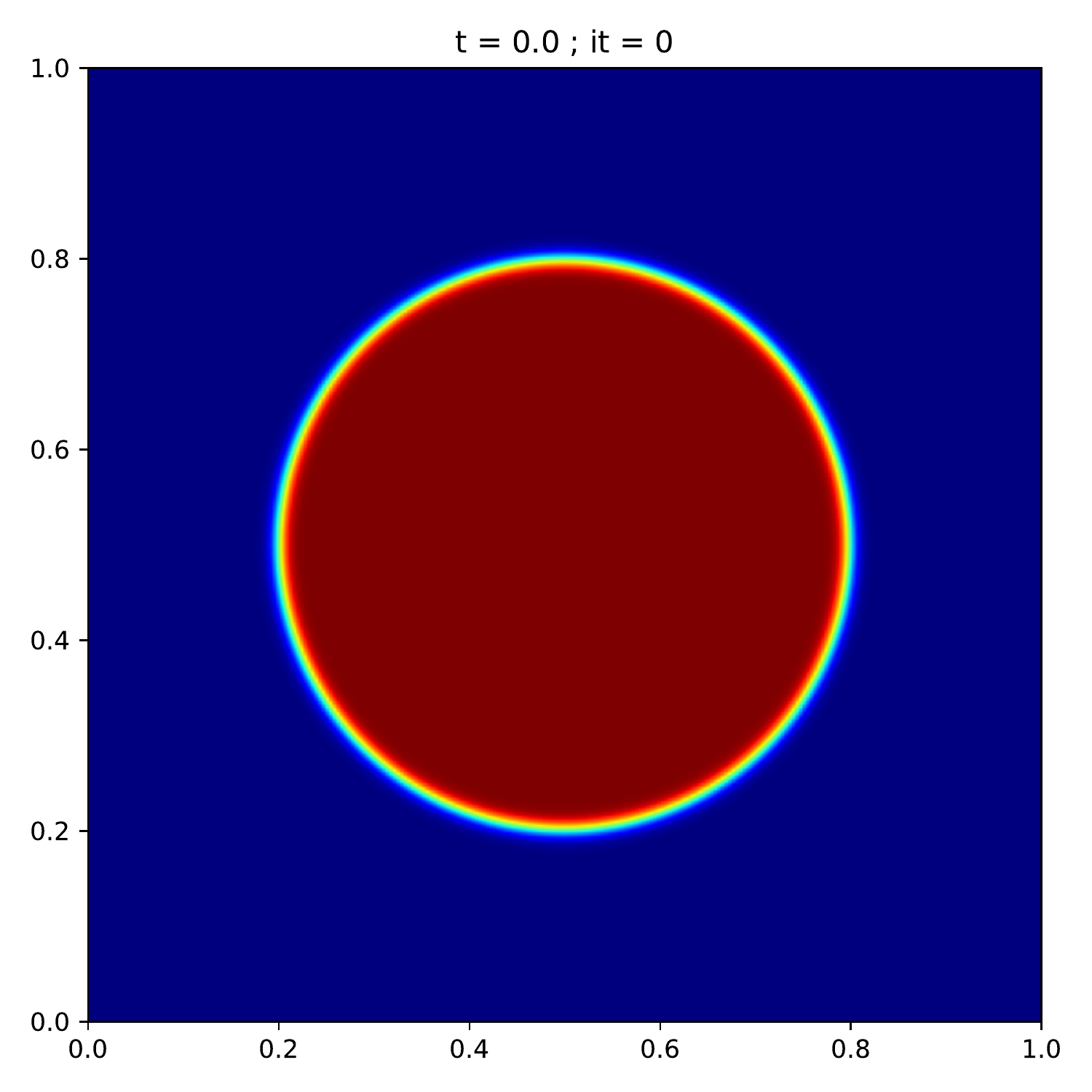}
         \includegraphics[width=0.2\textwidth]{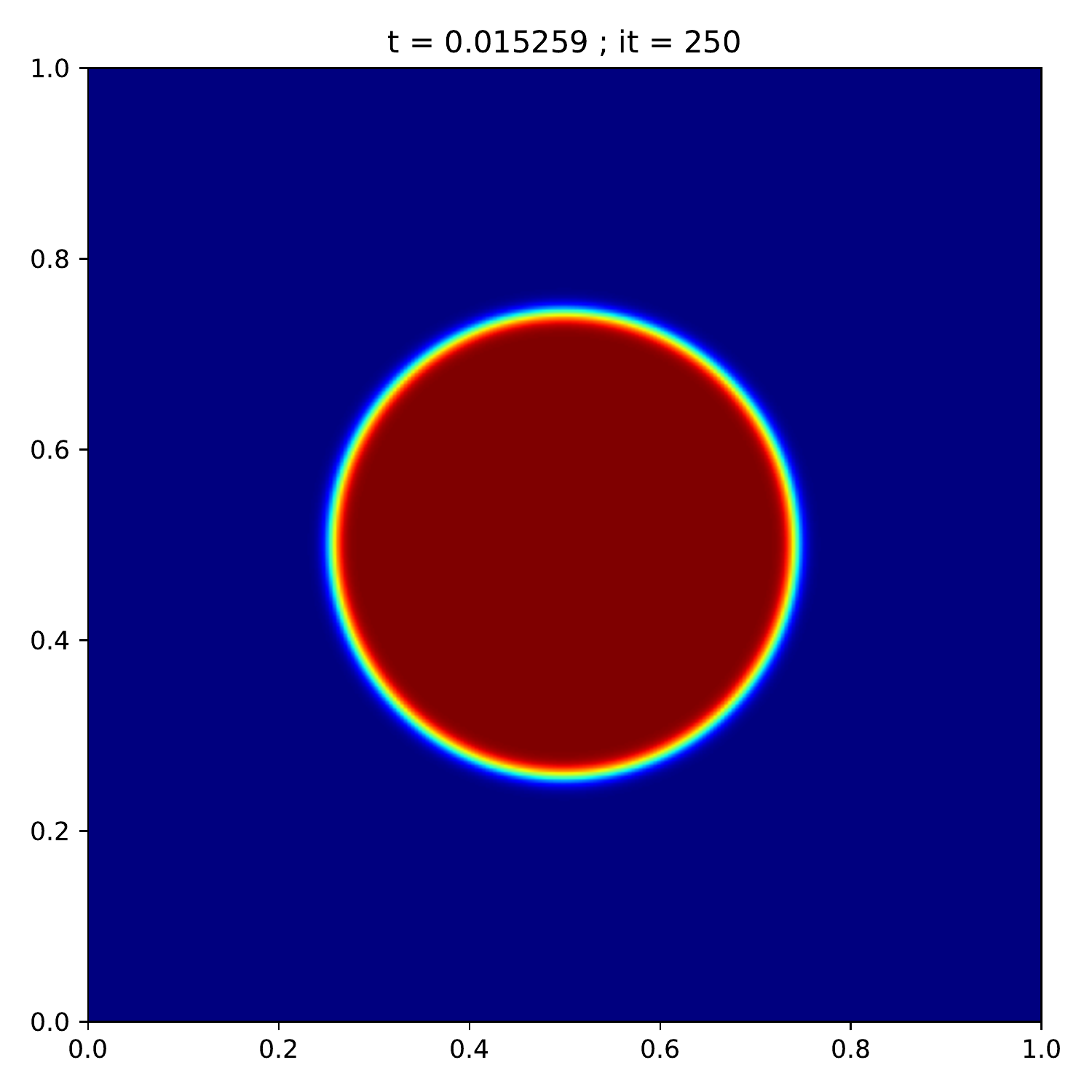}
     \includegraphics[width=0.2\textwidth]{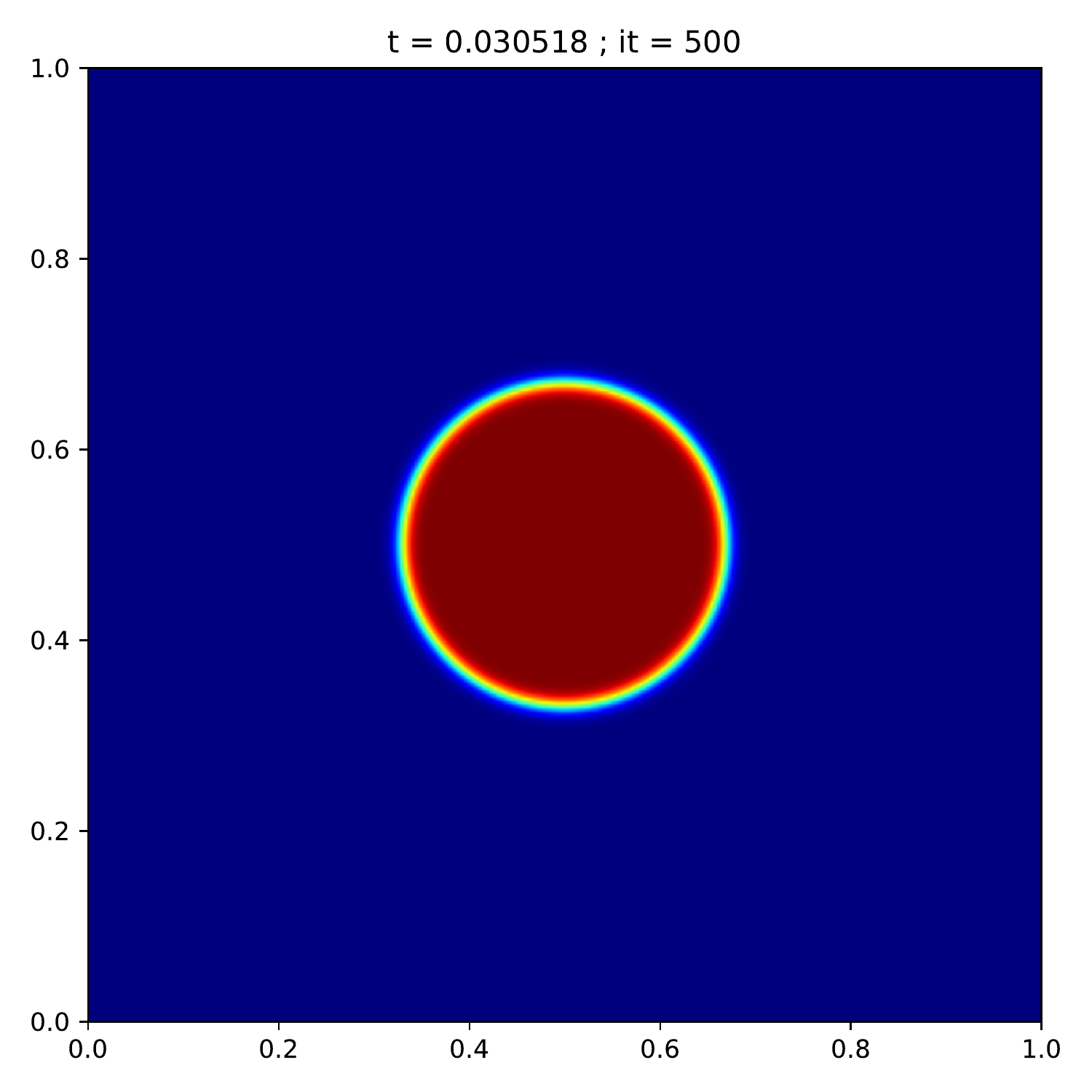}
          \includegraphics[width=0.2\textwidth]{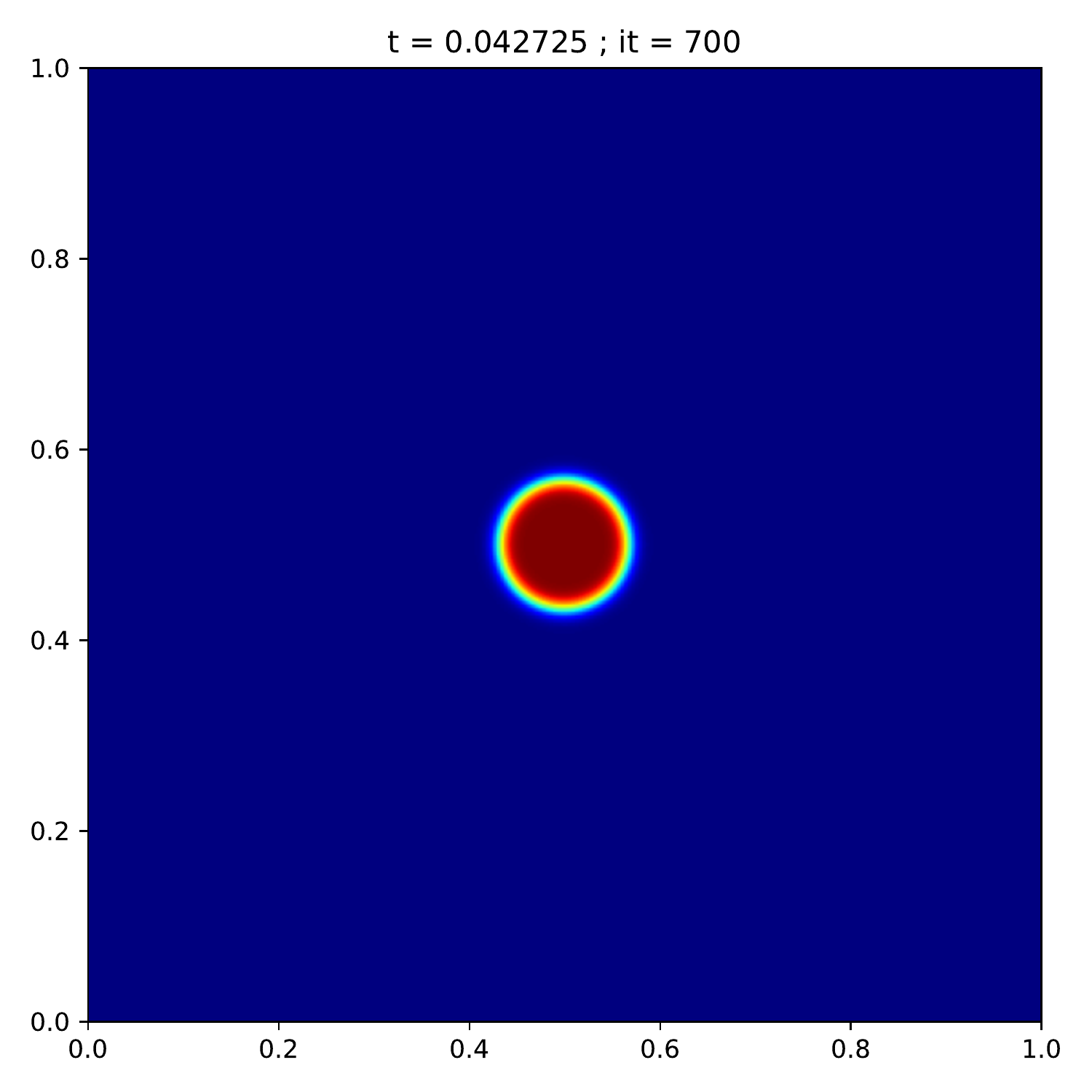} \\
      \includegraphics[width=0.2\textwidth]{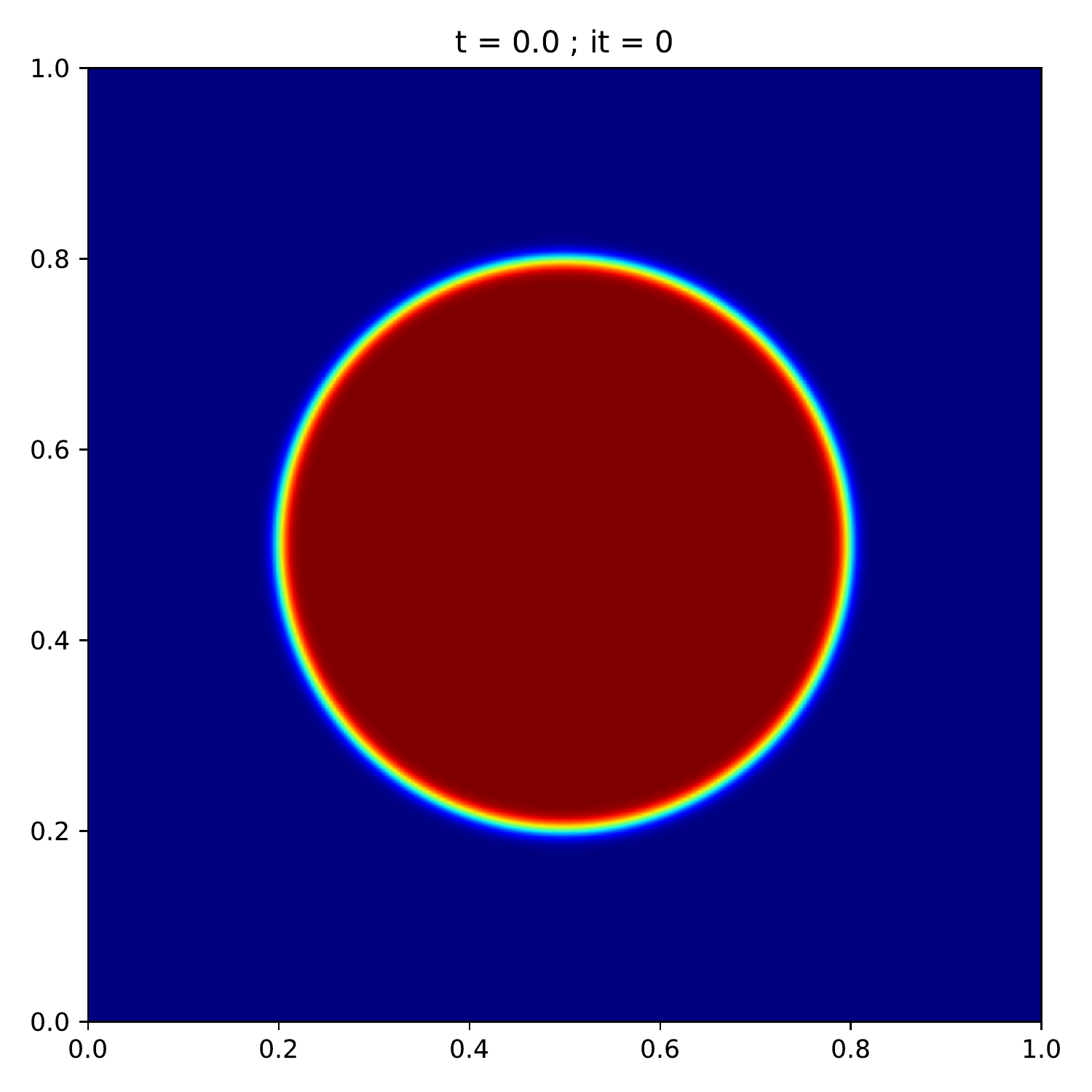}
           \includegraphics[width=0.2\textwidth]{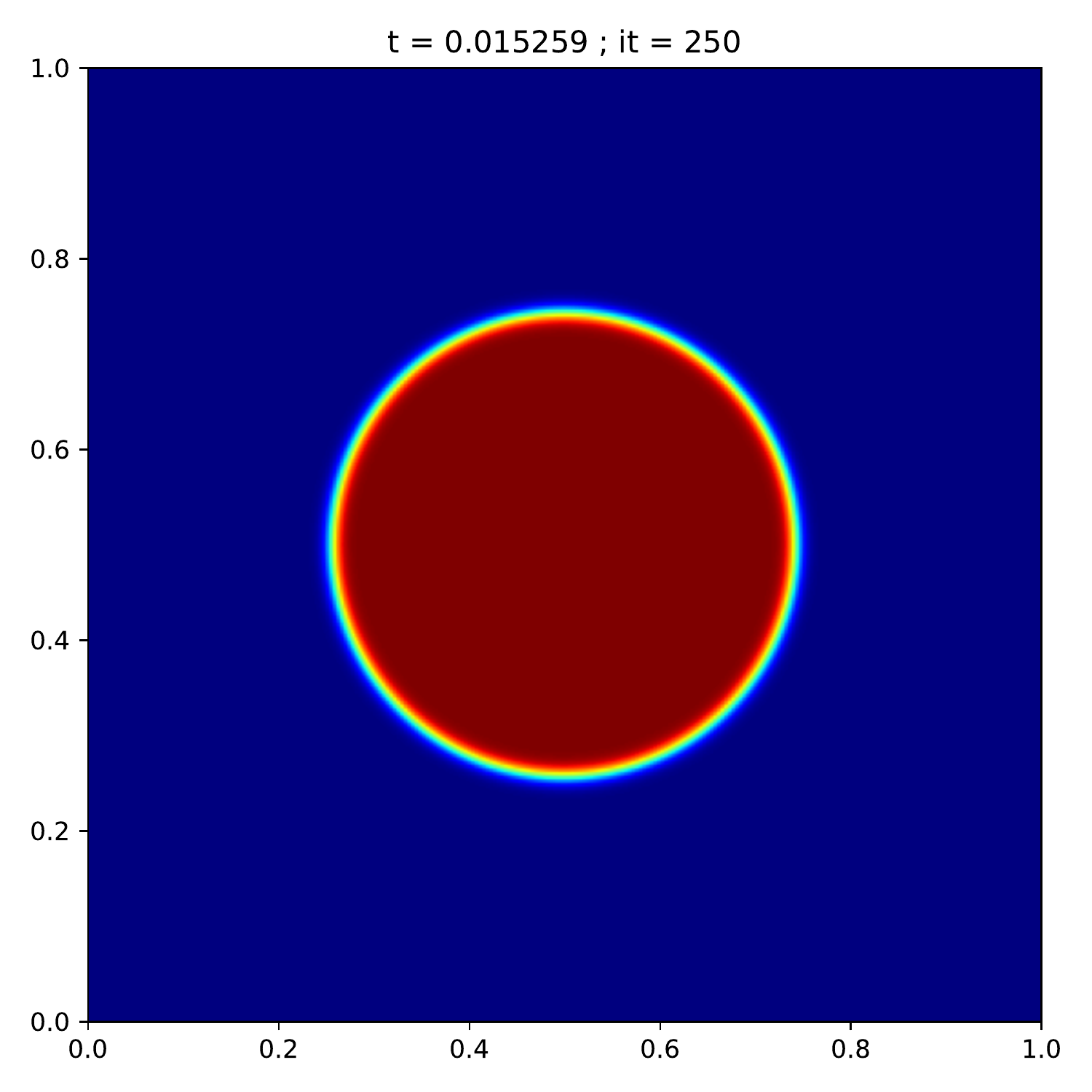}
            \includegraphics[width=0.2\textwidth]{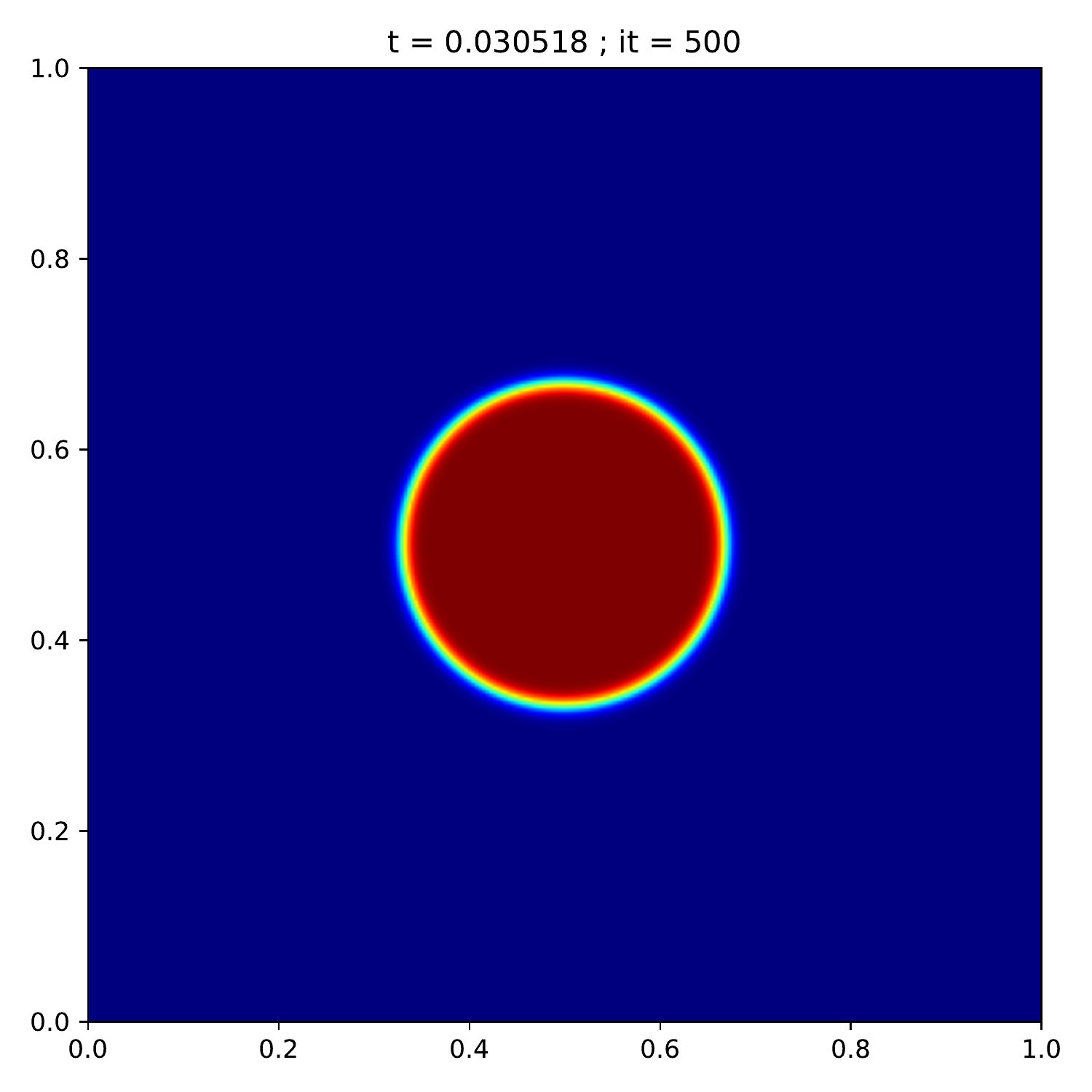}
             \includegraphics[width=0.2\textwidth]{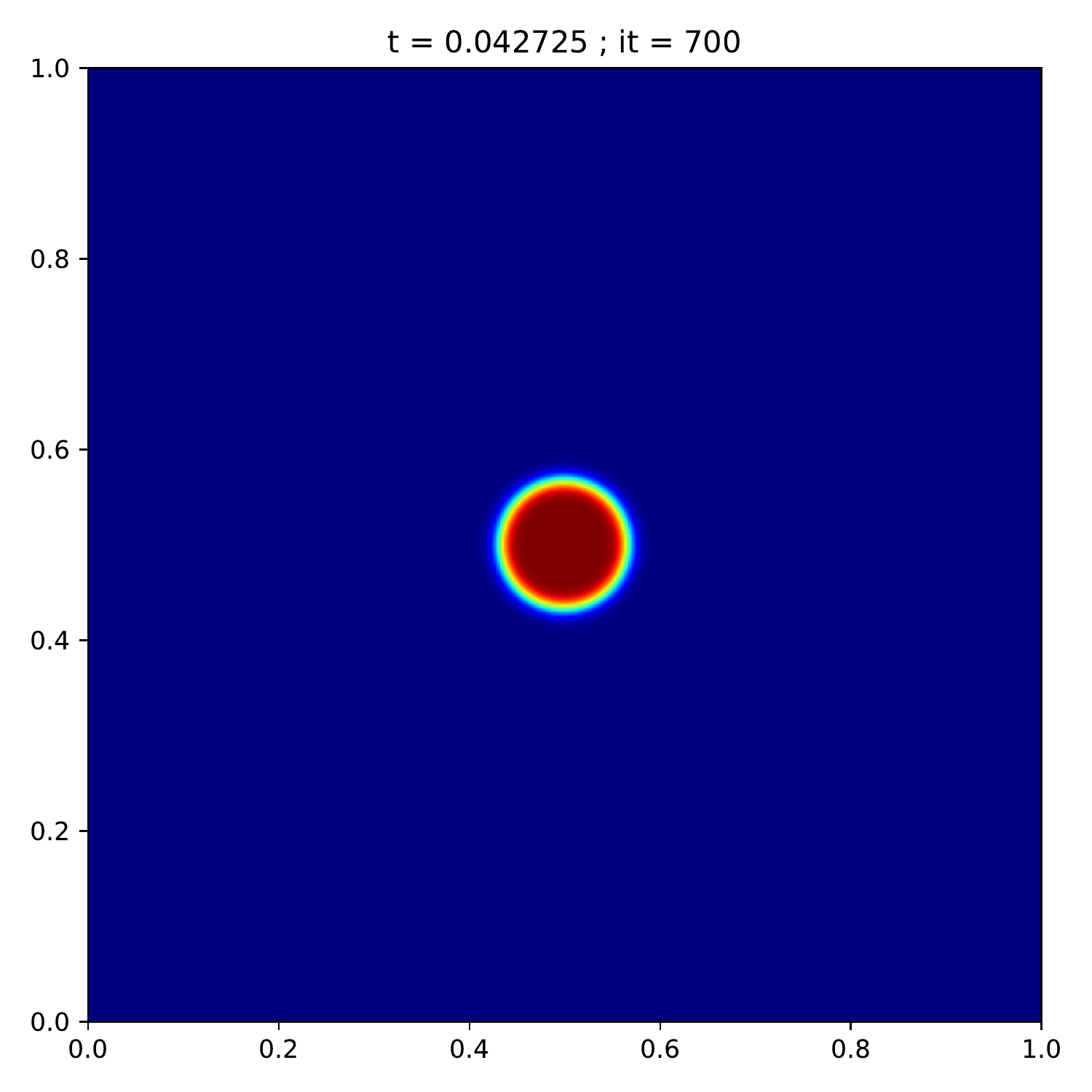}
    \caption{Comparison of the numerical semigroups $\S^{\text{AC}}_{\delta_t,\varepsilon,1}$, $\S^{\text{NN}}_{\theta,1}$ and $\S^{\text{NN}}_{\theta,2}$ for the approximation of the mean curvature 
    flow of a circle; each row corresponds to the evolution obtained at different times using,
respectively, the numerical semigroup  $\S^{\text{AC}}_{\delta_t,\varepsilon,1}$ and 
the networks $\S^{\text{NN}}_{\theta,1}$ and $\S^{\text{NN}}_{\theta,2}$. }
    \label{fig:valide_q_1}
\end{figure}

More precisely, we plot in each line of figure~\ref{fig:valide_q_1} the numerical approximation 
$u^{n}$  obtained at different times $t_n = n \delta_t$ using  respectively  the classical Lie splitting scheme $\S^{\text{AC}}_{\delta_t = \varepsilon^2,\varepsilon,1}$,
the first order network  $\S^{\text{NN}}_{\theta,1}$ and the second order one $\S^{\text{NN}}_{\theta,2}$. All three approximations of the mean curvature flow are so similar that it is difficult to observe any difference between them. At least on this particular example, the networks provide qualitatively the expected evolution.

\begin{figure}[htbp]
    \centering
    \includegraphics[width=0.45\textwidth]{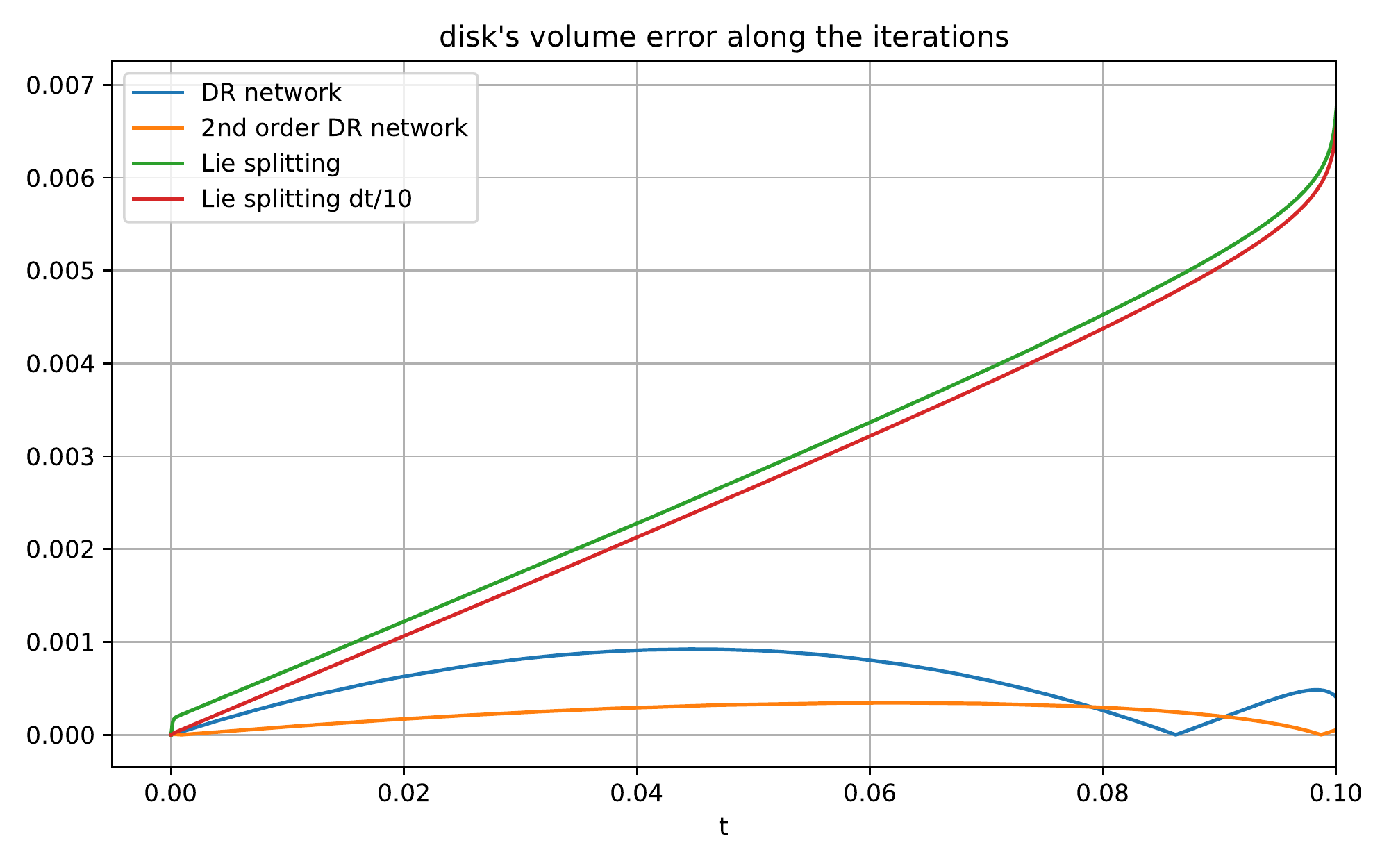}    
    \includegraphics[width=0.45\textwidth]{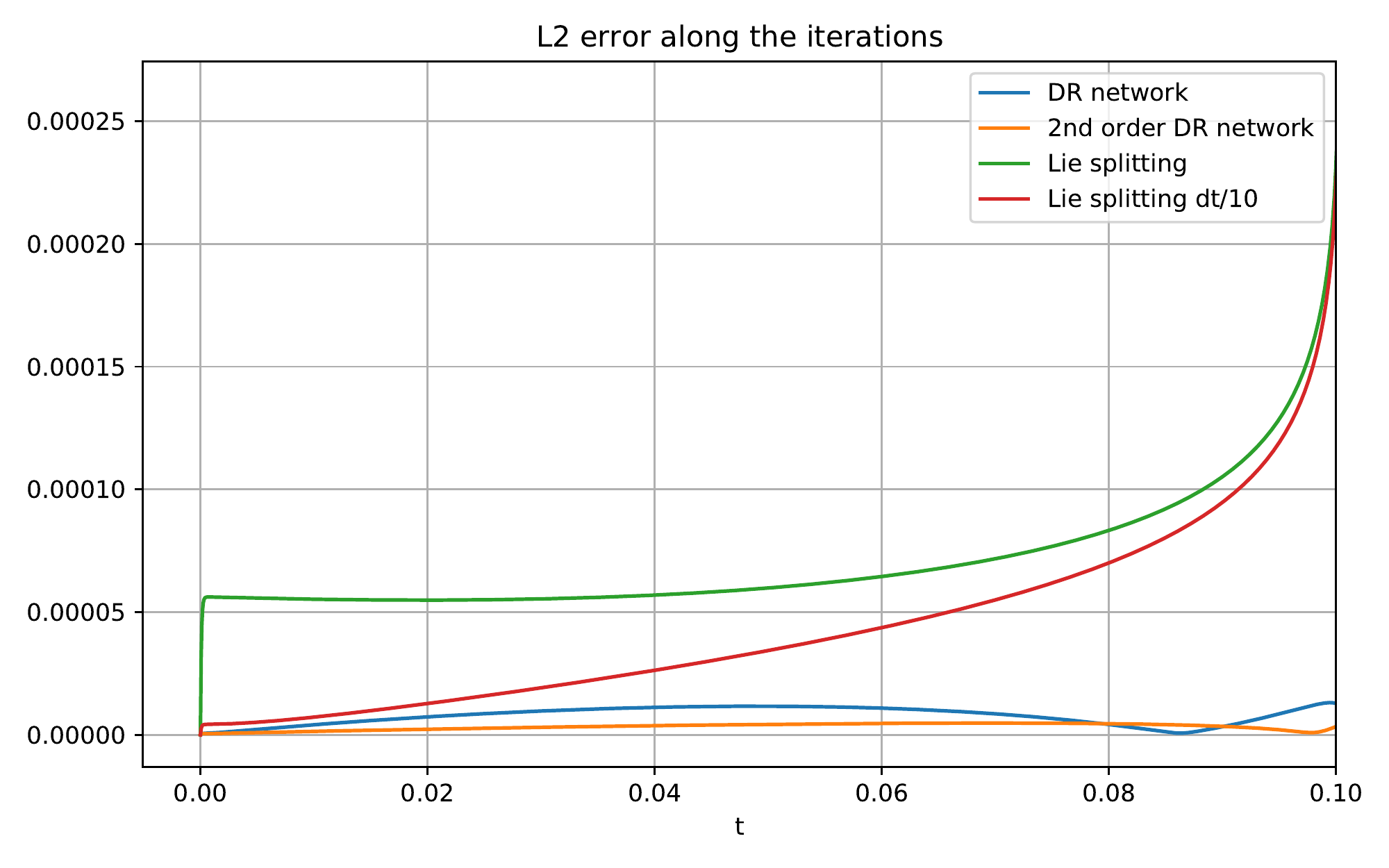}    
    \caption{
      \textcolor{black}{
    	Comparison of different schemes on the evolution of a disk by approximate mean curvature flow with initial radius $R_0 = 0.45$.
    	Plots in blue, orange, green, and red correspond to the error along iterations of, respectively, the trained networks
    	$\S^{\text{NN}}_{\theta,1}$ and $\S^{\text{NN}}_{\theta,2}$, and the schemes  $\S^{\text{AC}}_{\delta_t = \varepsilon^2,\varepsilon,1}$ 
    	and  $\S^{\text{AC}}_{\delta_t,\varepsilon,1}$ with $\delta_t = \epsilon^2/10$.
    	Left : disk's volume error, Right : $L^2$ error.}
    }
    \label{fig:radius_validation}
\end{figure}

In order to have more quantitative error estimates, we display on figure~\ref{fig:radius_validation} the error on
the integral of the phase field function $\varphi_{R(t)}$ of a circle of radius $R(t)$
during the iterations knowing that the evolution by mean curvature
of an initial disk of radius $R_0$ is a disk of radius
$R(t) = \sqrt{R_0^2 - 2 t}$, at least until the extinction time $t_0 = \frac{1}{2} R_0^2$.
More precisely, we plot:

\begin{itemize}
 \item the volume error
 $$\textcolor{black}{n \mapsto \left| \int_Q u^{n}(x) dx - \int_Q \varphi_{\sqrt{R_0^2 - 2 n \delta_t}}(x) dx\right|^2},$$ 
 \item the (squared) $L^2$ error
$$\textcolor{black}{n \mapsto  \int_Q \left( u^{n}(x)- \varphi_{\sqrt{R_0^2 - 2 n \delta_t}}(x) \right)^2 dx,}$$ 
 \end{itemize}

using 
\begin{itemize}
 \item the trained network $\S^{\text{NN}}_{\theta,1}$ in blue
 \item the trained network $\S^{\text{NN}}_{\theta,2}$ in orange
 \item the Lie splitting scheme  $\S^{\text{AC}}_{\delta_t,\varepsilon,1}$ with $\delta_t = \epsilon^2$ in green
 \item the Lie splitting scheme  $\S^{\text{AC}}_{\delta_t,\varepsilon,1}$ with $\delta_t = \epsilon^2/10$ in red
\end{itemize}

These results clearly show that both our neural networks $\S^{\text{NN}}_{\theta,1}$ and $\S^{\text{NN}}_{\theta,2}$ provide more accurate numerical schemes than the classical discretization of the Allen-Cahn equation, even using a smaller time step $\delta_t$. 

This first validation test shows also that the training of our networks manages to correct both the errors of discretization of the Allen-Cahn equation, and the modeling errors of the approximation by the Allen-Cahn equation of the mean curvature flow.  \\

\begin{figure}[!htbp]
    \centering
    \includegraphics[width=0.2\textwidth]{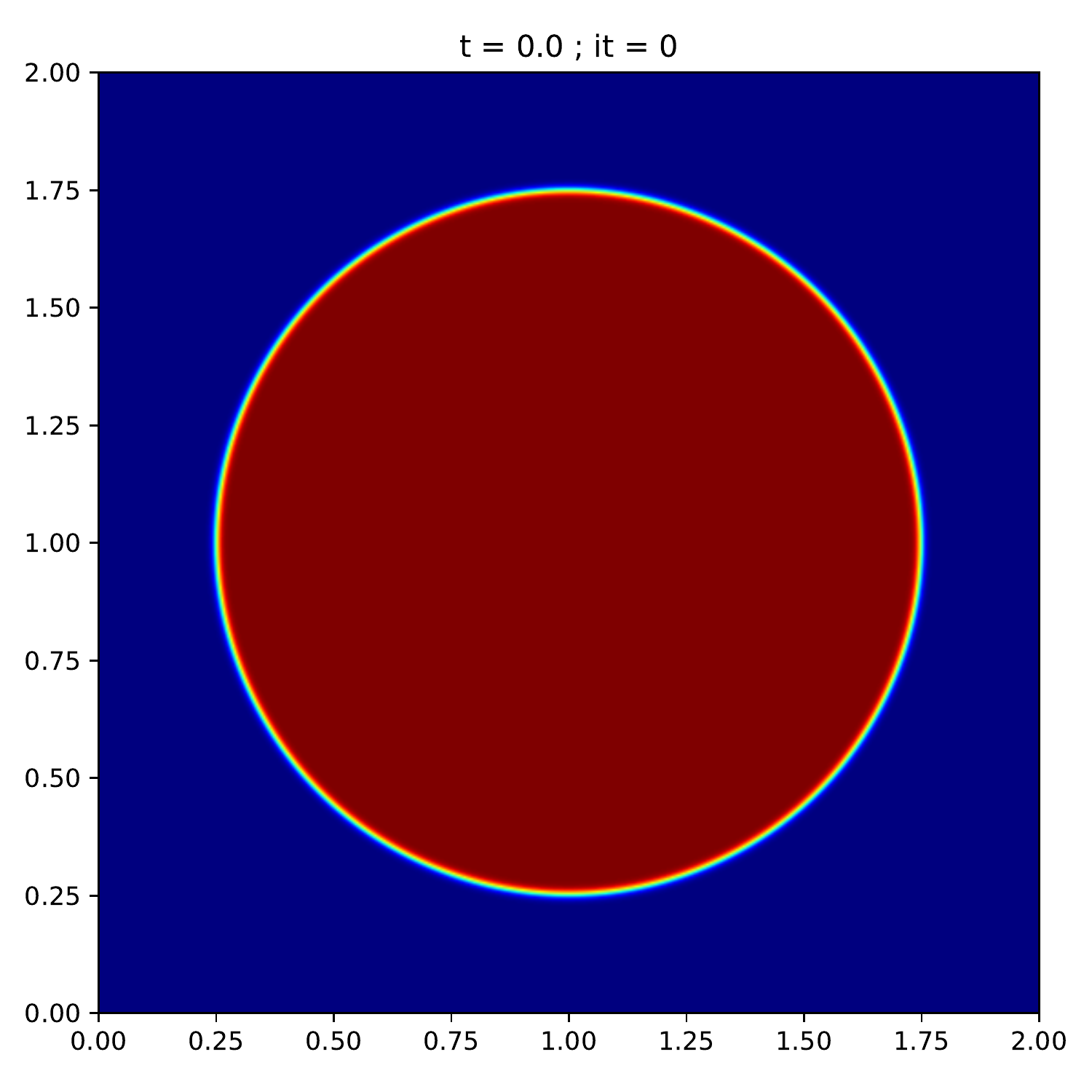}
    \includegraphics[width=0.2\textwidth]{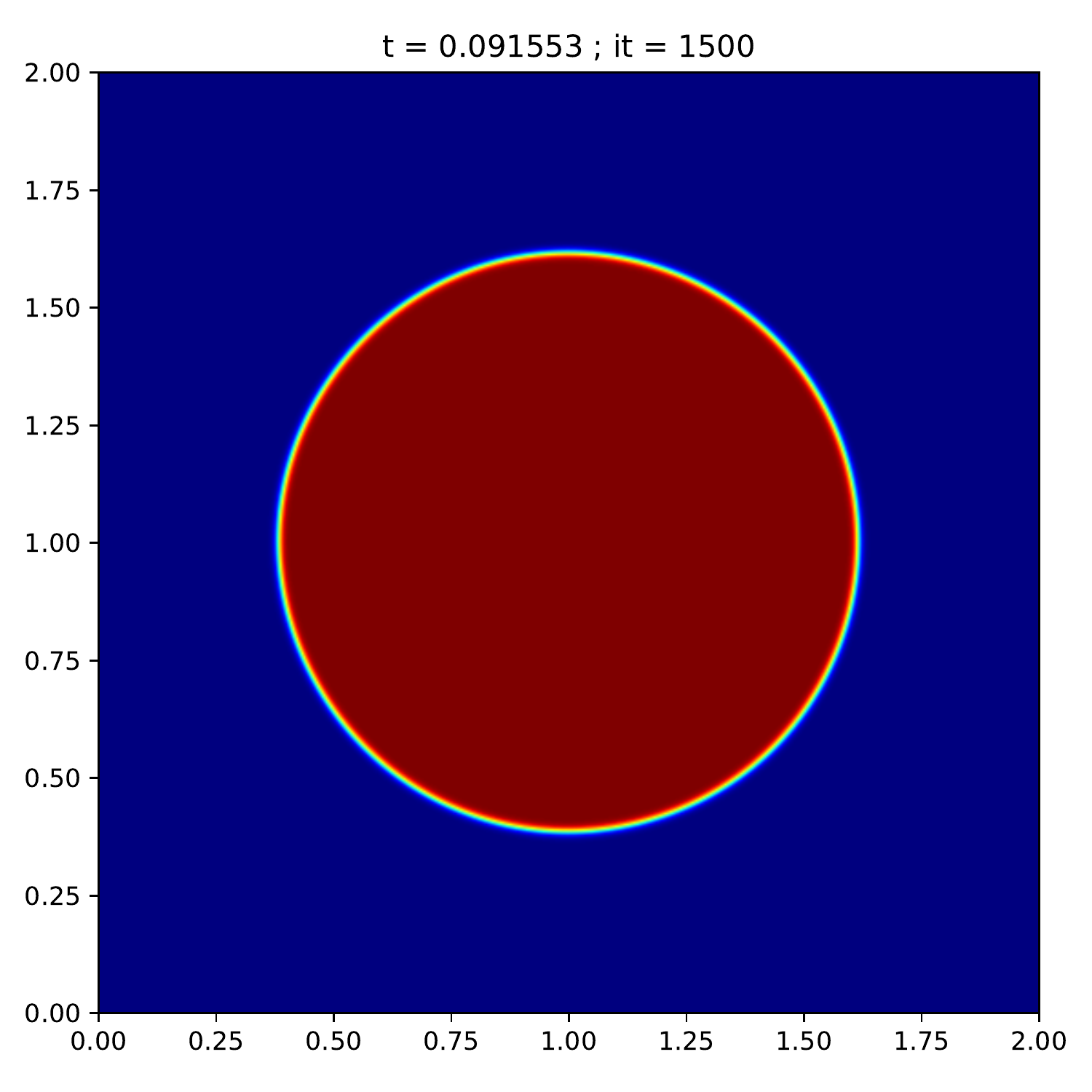}
    \includegraphics[width=0.2\textwidth]{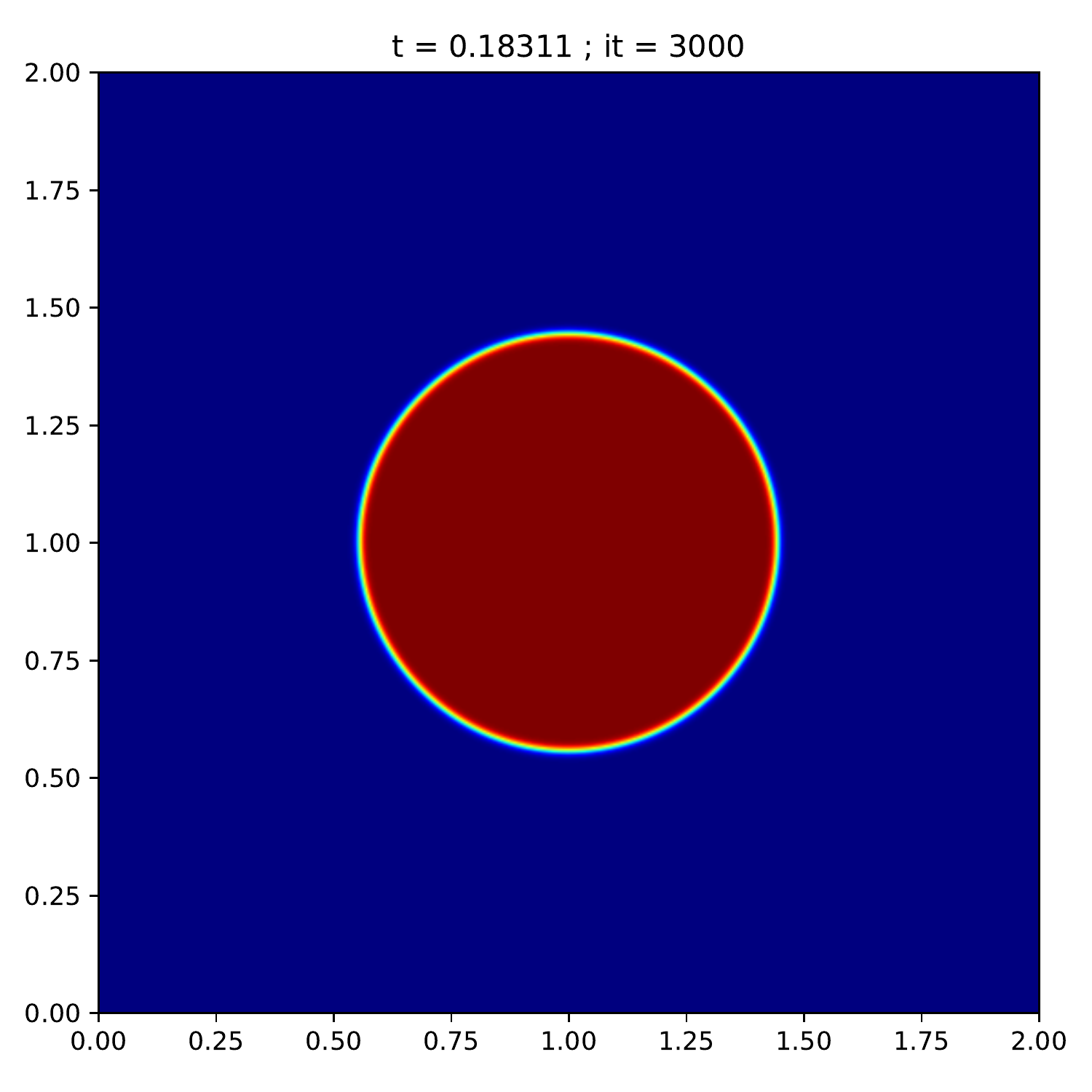}
    \includegraphics[width=0.2\textwidth]{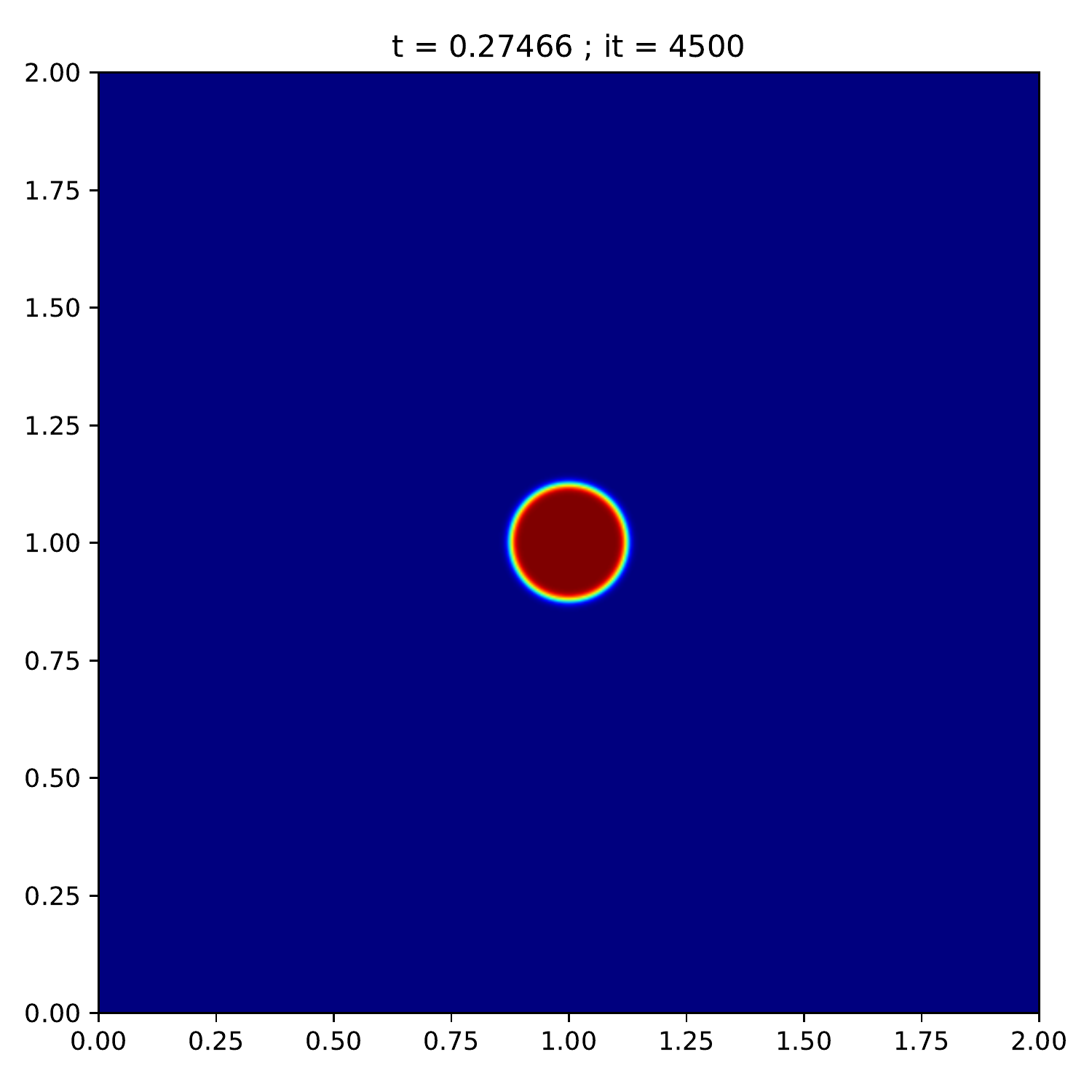}

    \caption{\textcolor{black}{Evolution along the iterations starting from an initial circle of radius $R_0 = 0.75$ on a domain $Q$ twice as large using the network $\S^{\text{NN}}_{\theta,1}$.}}
    \label{fig:valide_q1_raduis_075}
\end{figure}

\begin{figure}[!htbp]
    \centering
    \includegraphics[width=0.8\textwidth]{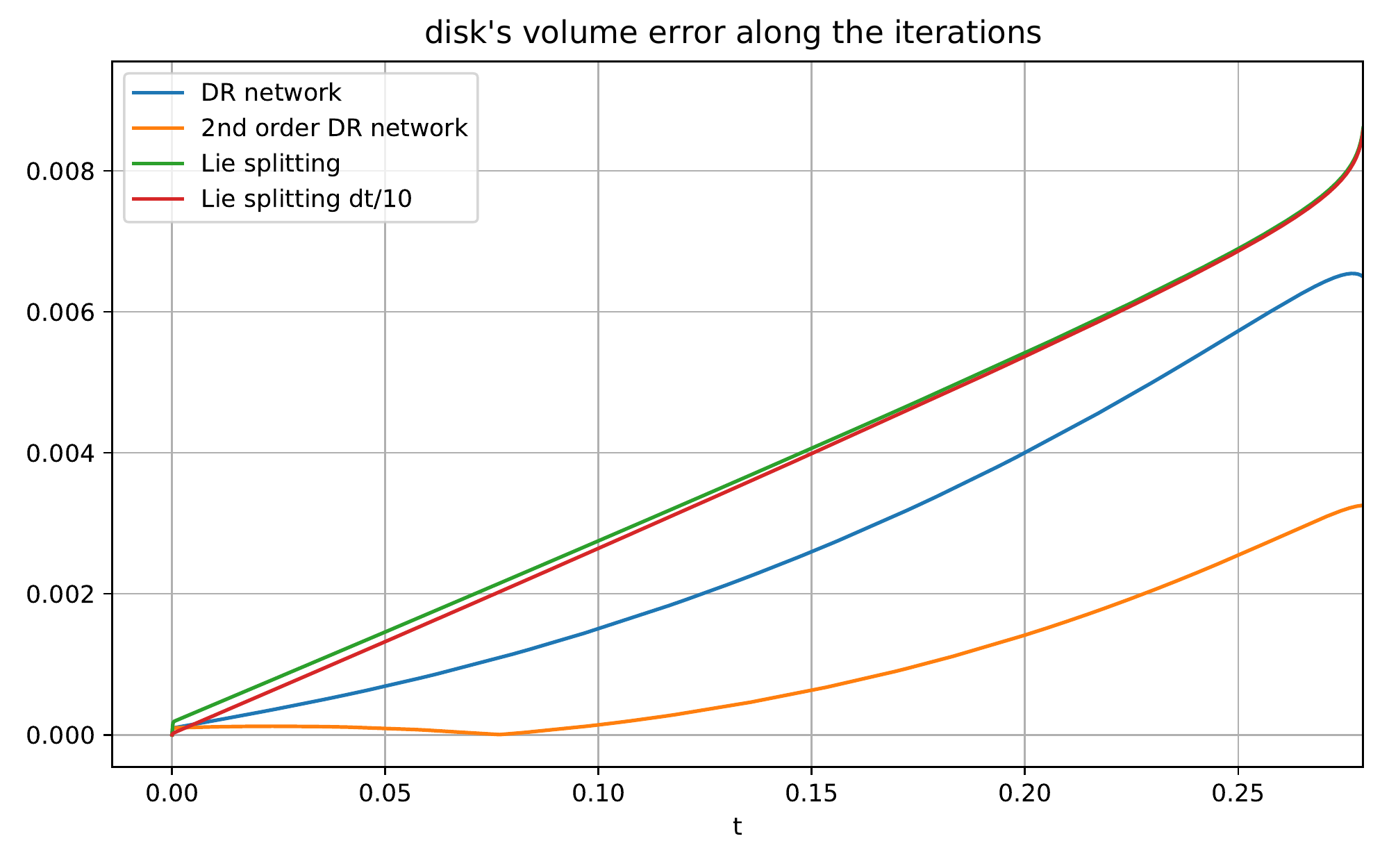}    
    \caption{
    	Comparison of different schemes on the evolution of the radius of a circle with initial radius $R_0= 0.75$ evolving by mean curvature flow.
    	Plots in blue, orange, green, and red correspond to the error
    	$n \mapsto   \left( \int_Q u^{n}(x) dx - \int_Q \varphi_{\sqrt{R_0^2 - 2 n \delta_t}}(x) dx \right)^2 $
    	along iterations of, respectively, the trained networks
    	$\S^{\text{NN}}_{\theta,1}$ and $\S^{\text{NN}}_{\theta,2}$,
    	and the schemes  $\S^{\text{AC}}_{\delta_t = \varepsilon^2,\varepsilon,1}$ 
    	and  $\S^{\text{AC}}_{\delta_t,\varepsilon,1}$ with $\delta_t = \epsilon^2/10$.
    }
    \label{fig:radius_validation_prim}
\end{figure}

\textcolor{black}{We now repeat the experiment considering instead a starting radius $R_0$ which does not belong to the training interval $[0.05,0.45]$.
This will be useful in assessing the ability of our neural networks to generalize well, i.e to handle properly new configurations which are not in the training database. We use a twice as large computational domain 
$Q = [0,2]^2$ which can contain an initial circle of radius $R_0 = 0.75$. We use the same networks as before without requiring any re-training. As the parameters $\delta_x$, $\varepsilon$ and 
$\delta_t$ are fixed, observe that due to the visualization scale, the width of the diffuse interface  is twice as small on each image
of figure  \ref{fig:valide_q1_raduis_075} in comparison with figure  \ref{fig:valide_q_1}.  In figure \ref{fig:radius_validation_prim} we plot the volume error along iterations associated with the estimated radius $R(t)$ for each method. As observed previously, the second order network $S^{NN}_{2,\theta}$
is the most accurate method, although its performance has marginally deteriorated as compared to the previous experiment.  Still, these results confirm that our networks generalize well.
}

\begin{figure}[htbp]
    \centering
    \includegraphics[width=0.2\textwidth]{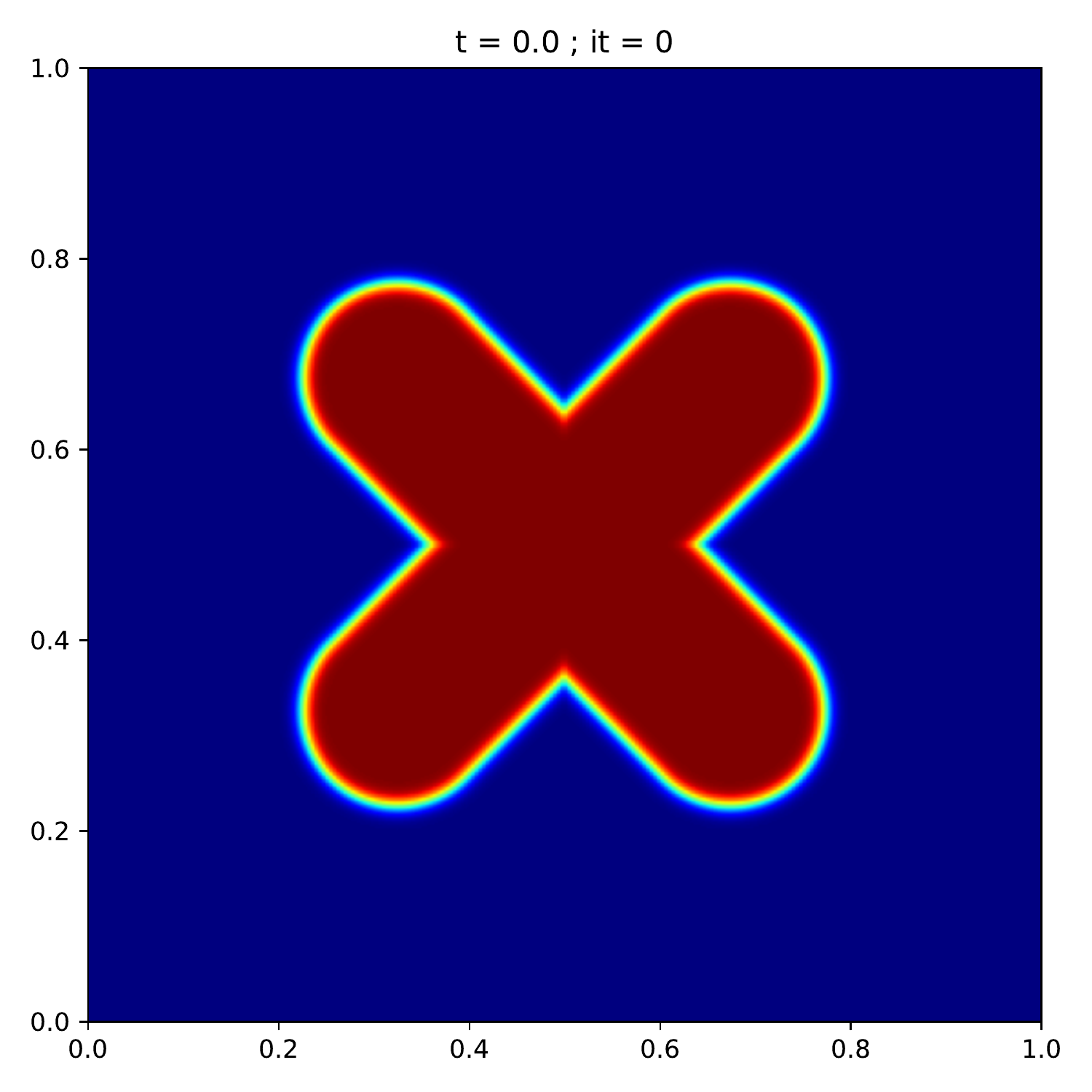}
    \includegraphics[width=0.2\textwidth]{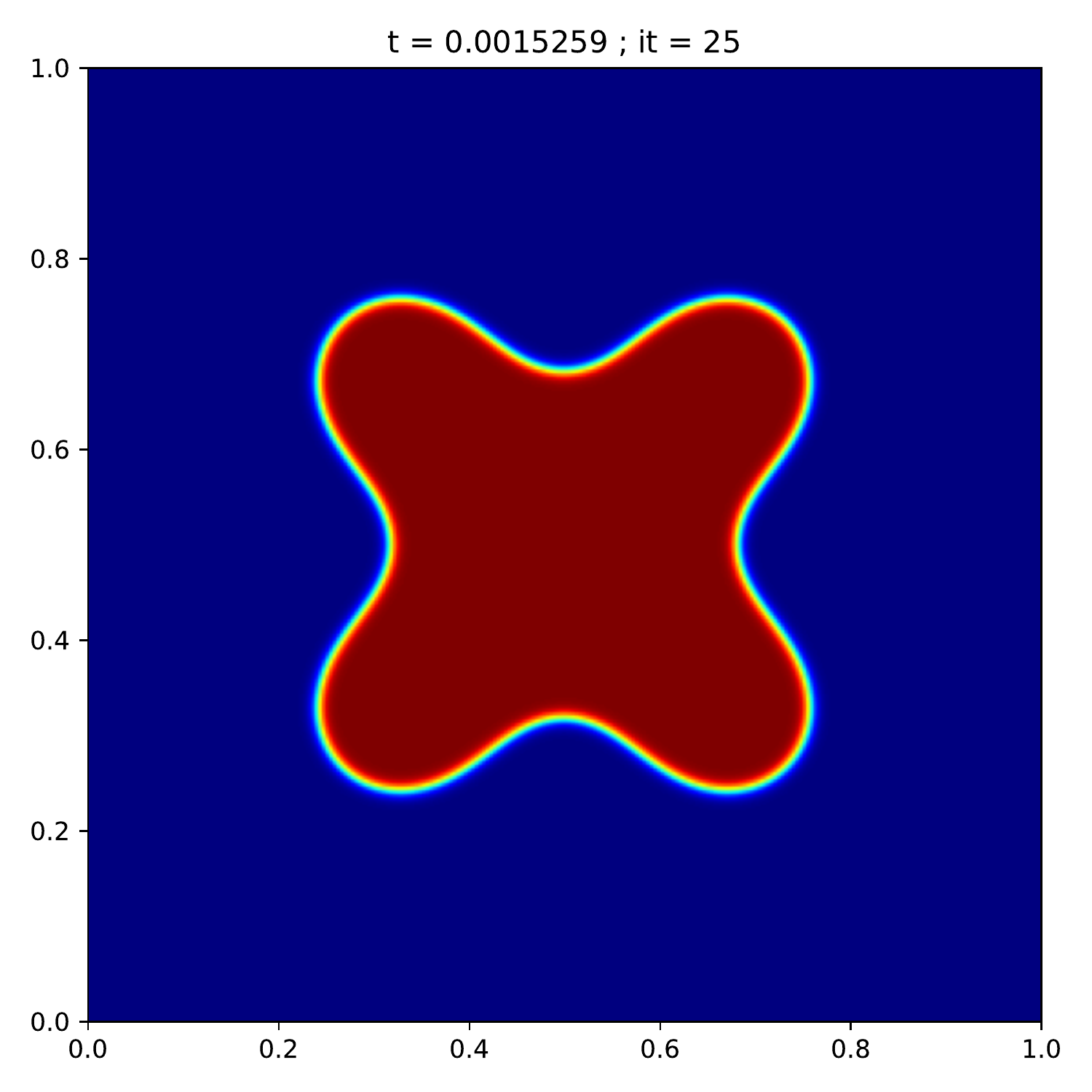}
    \includegraphics[width=0.2\textwidth]{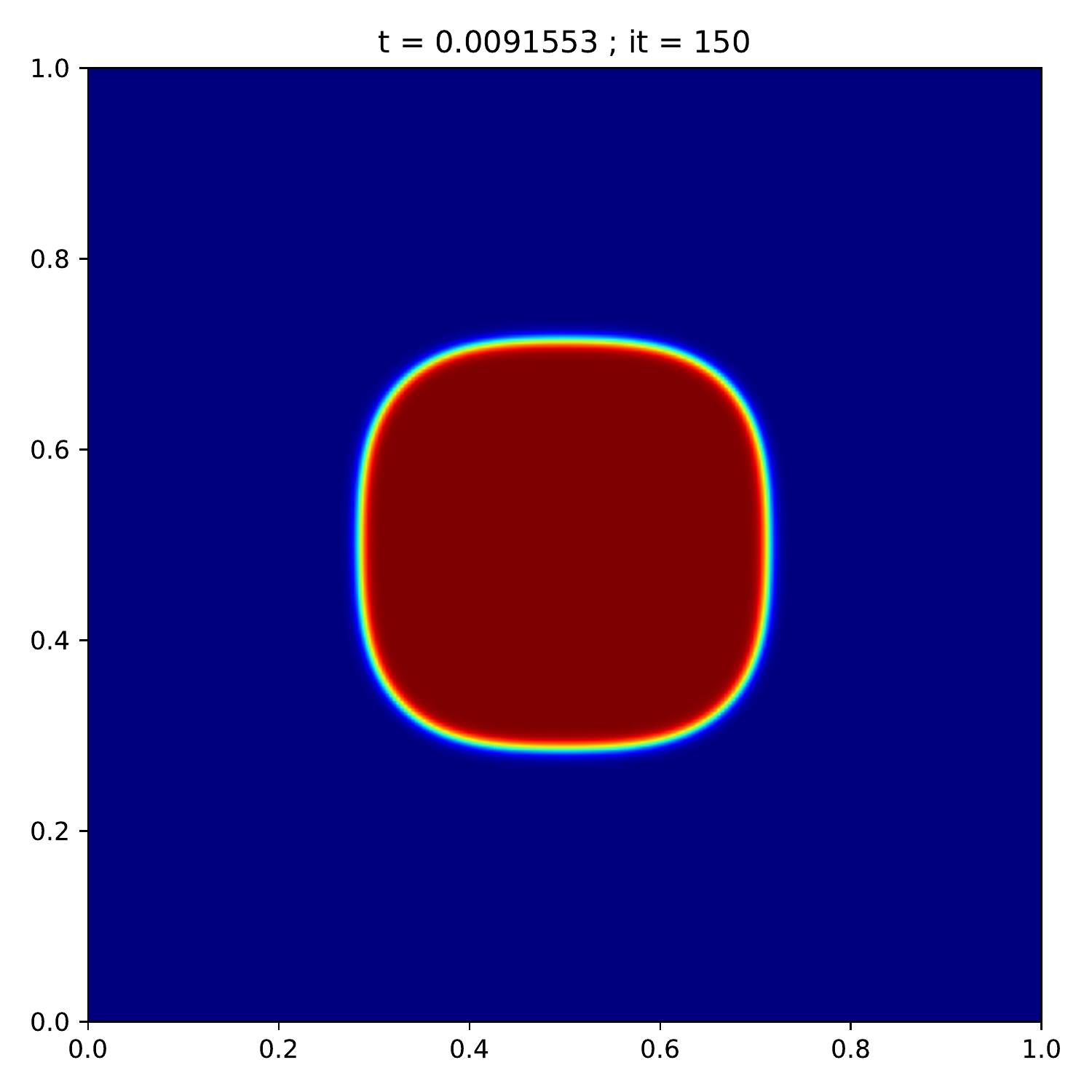}
    \includegraphics[width=0.2\textwidth]{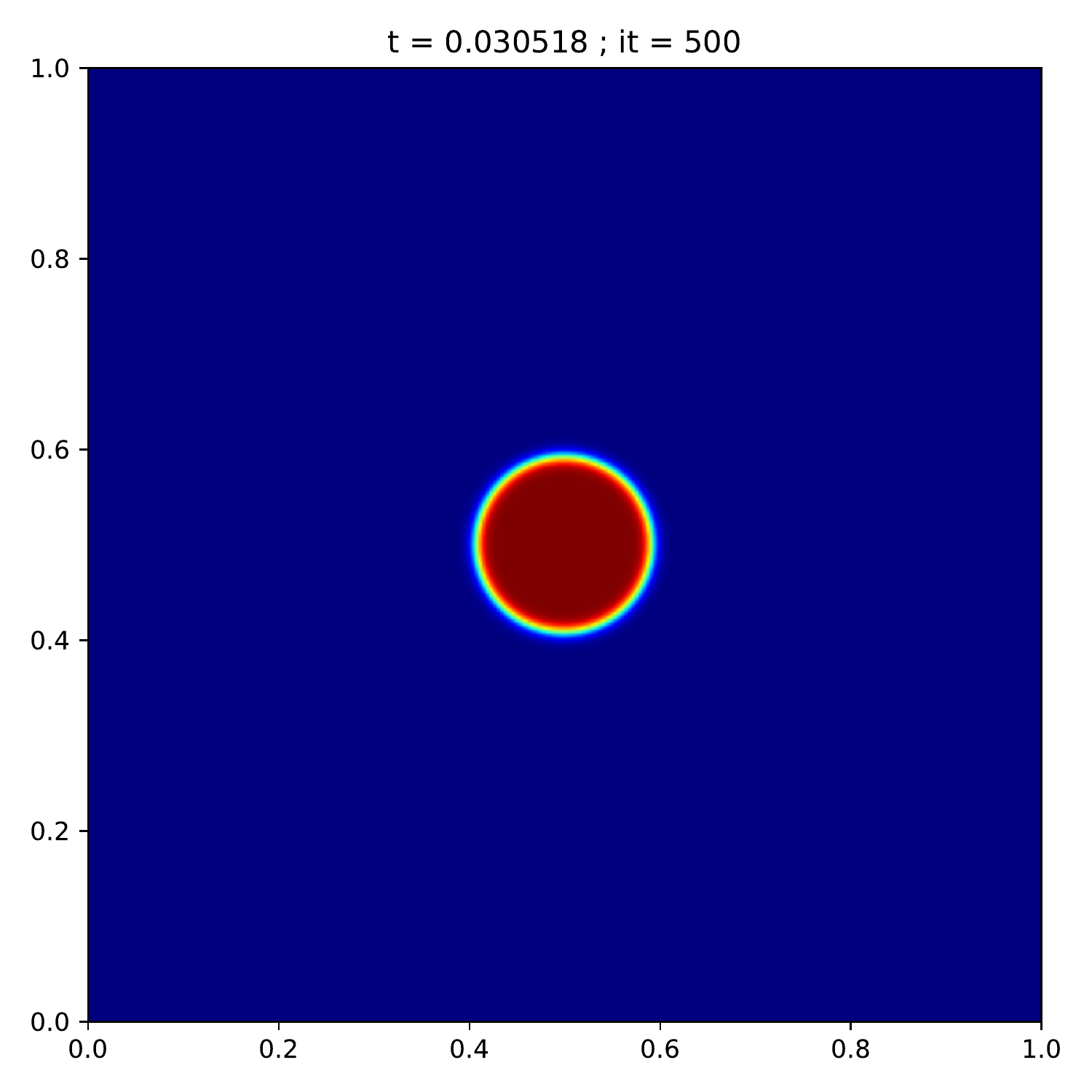} \\
    \includegraphics[width=0.2\textwidth]{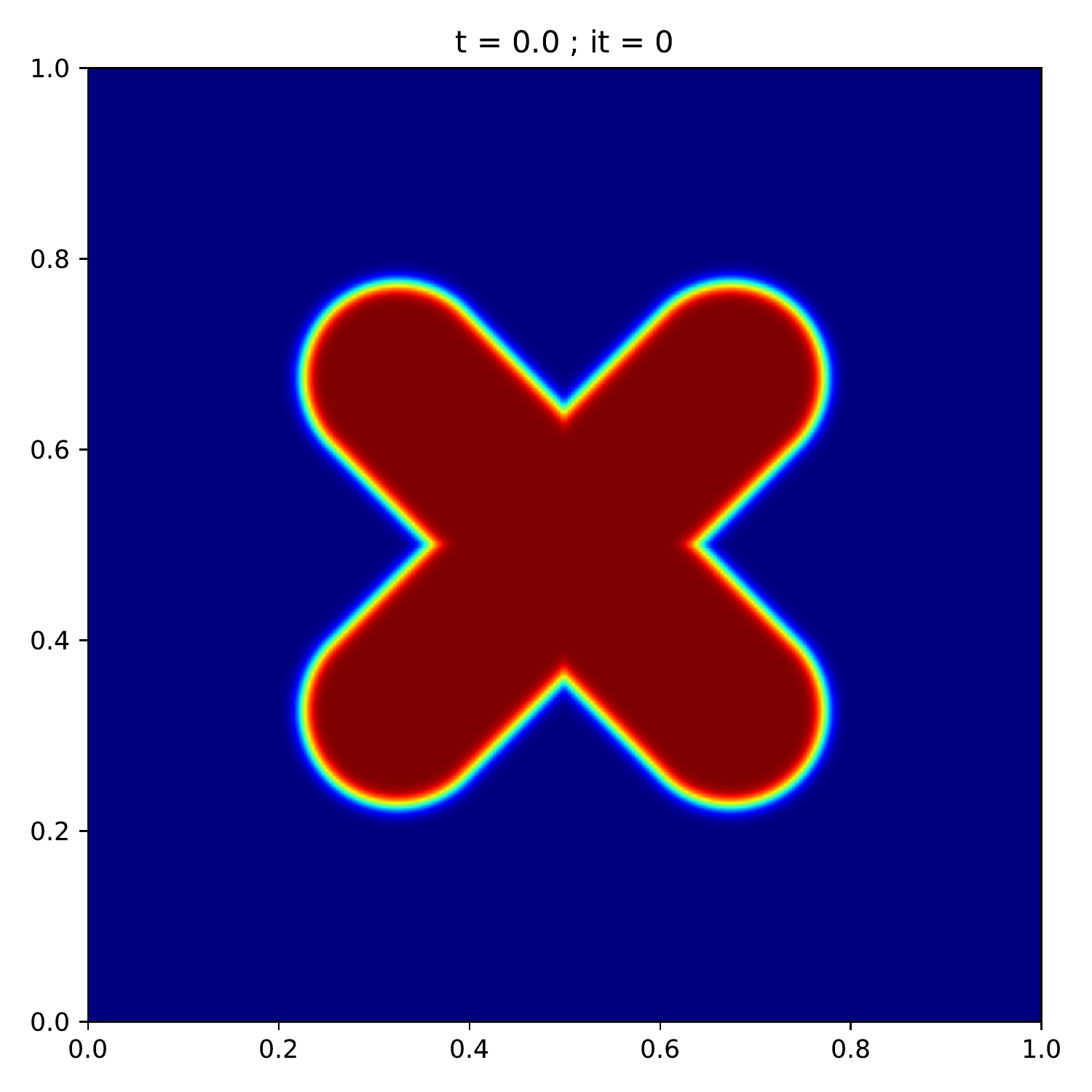}
     \includegraphics[width=0.2\textwidth]{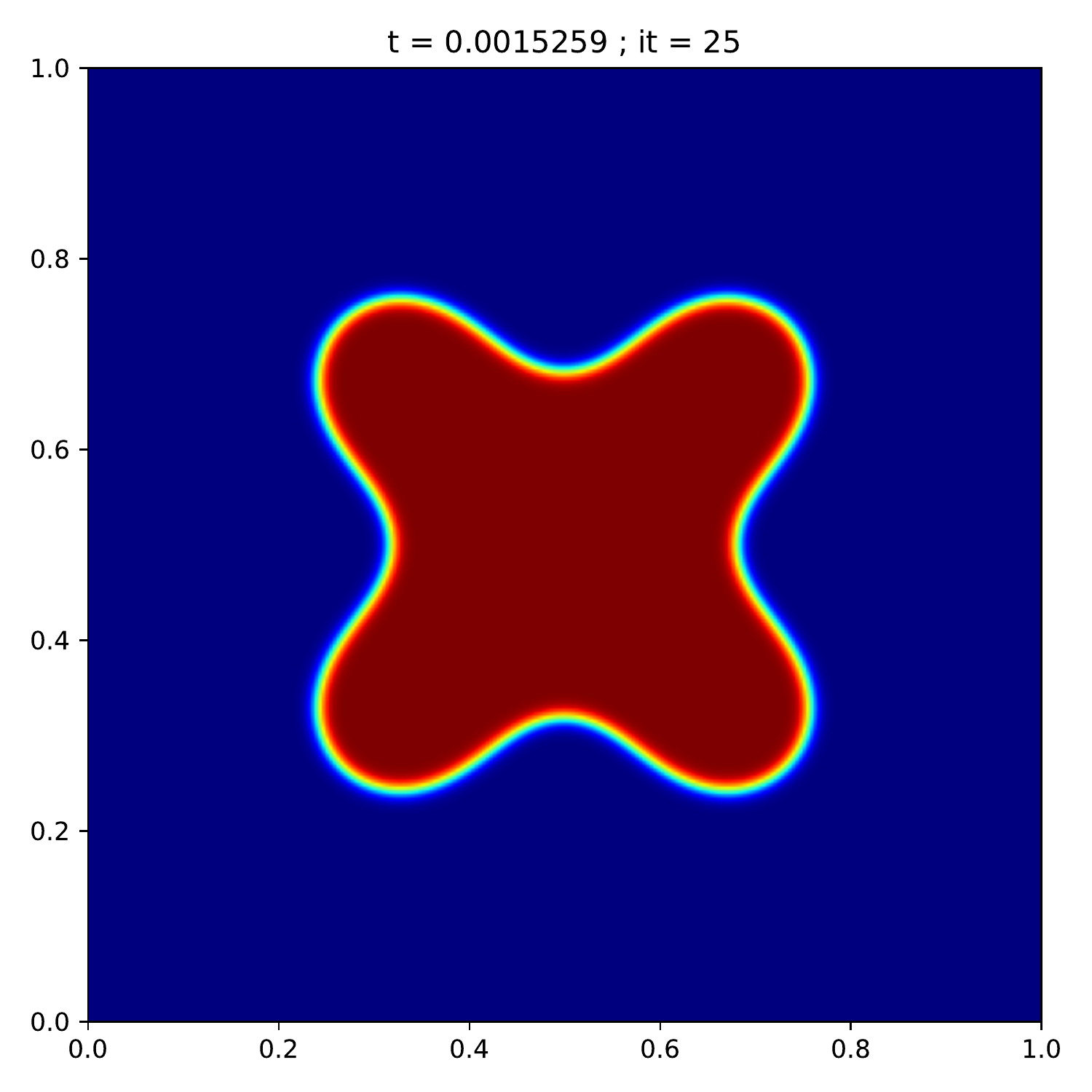}
      \includegraphics[width=0.2\textwidth]{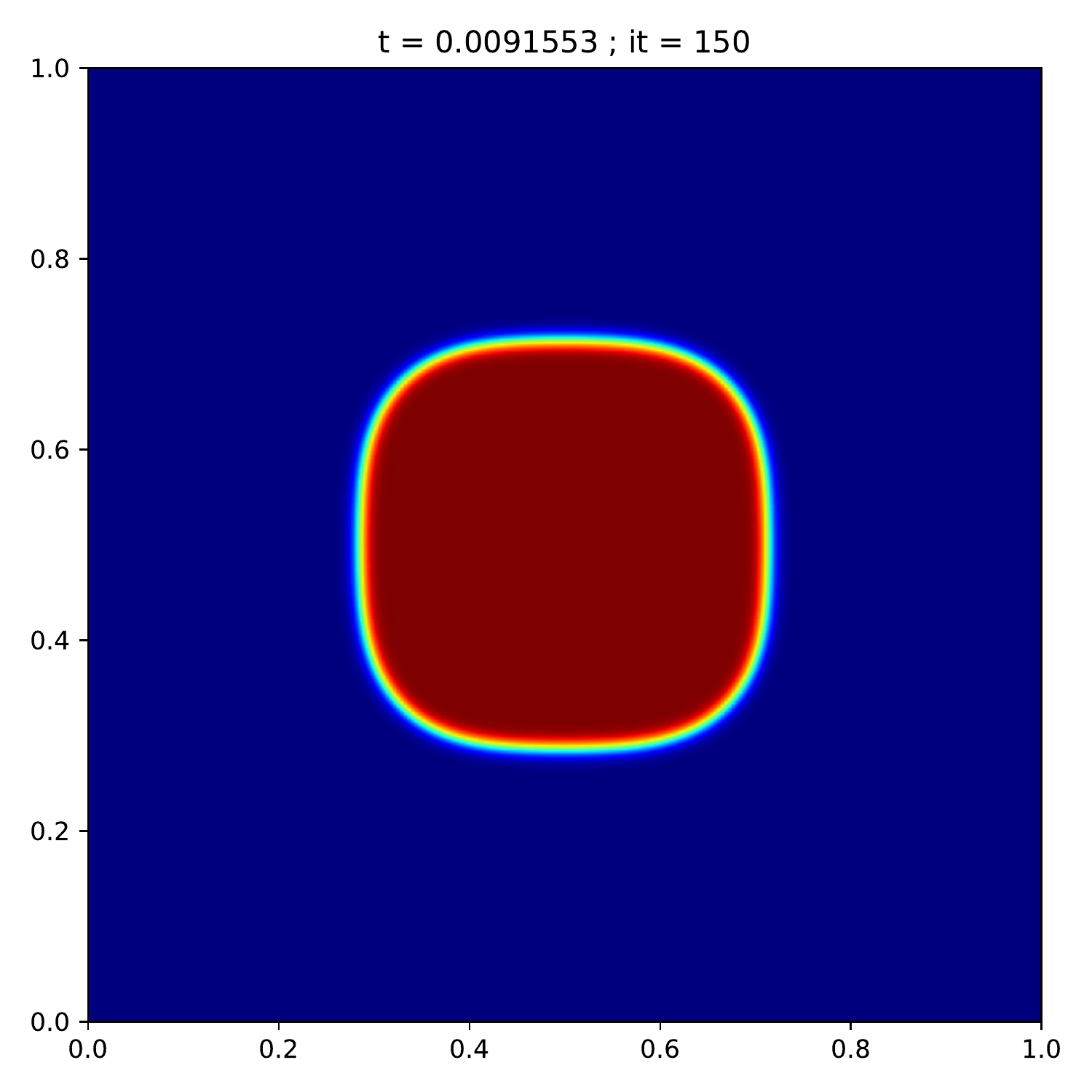}
       \includegraphics[width=0.2\textwidth]{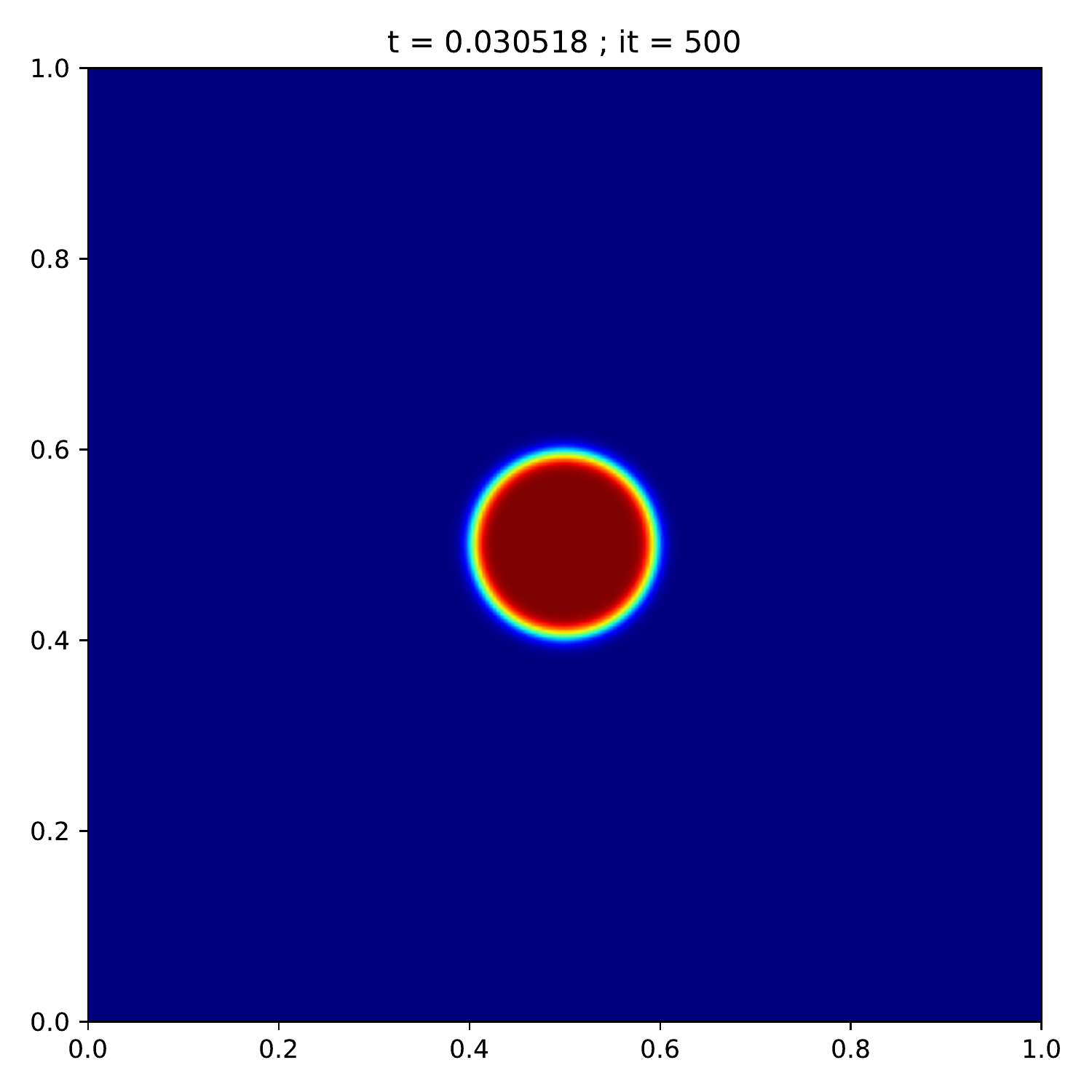} \\
         \includegraphics[width=0.2\textwidth]{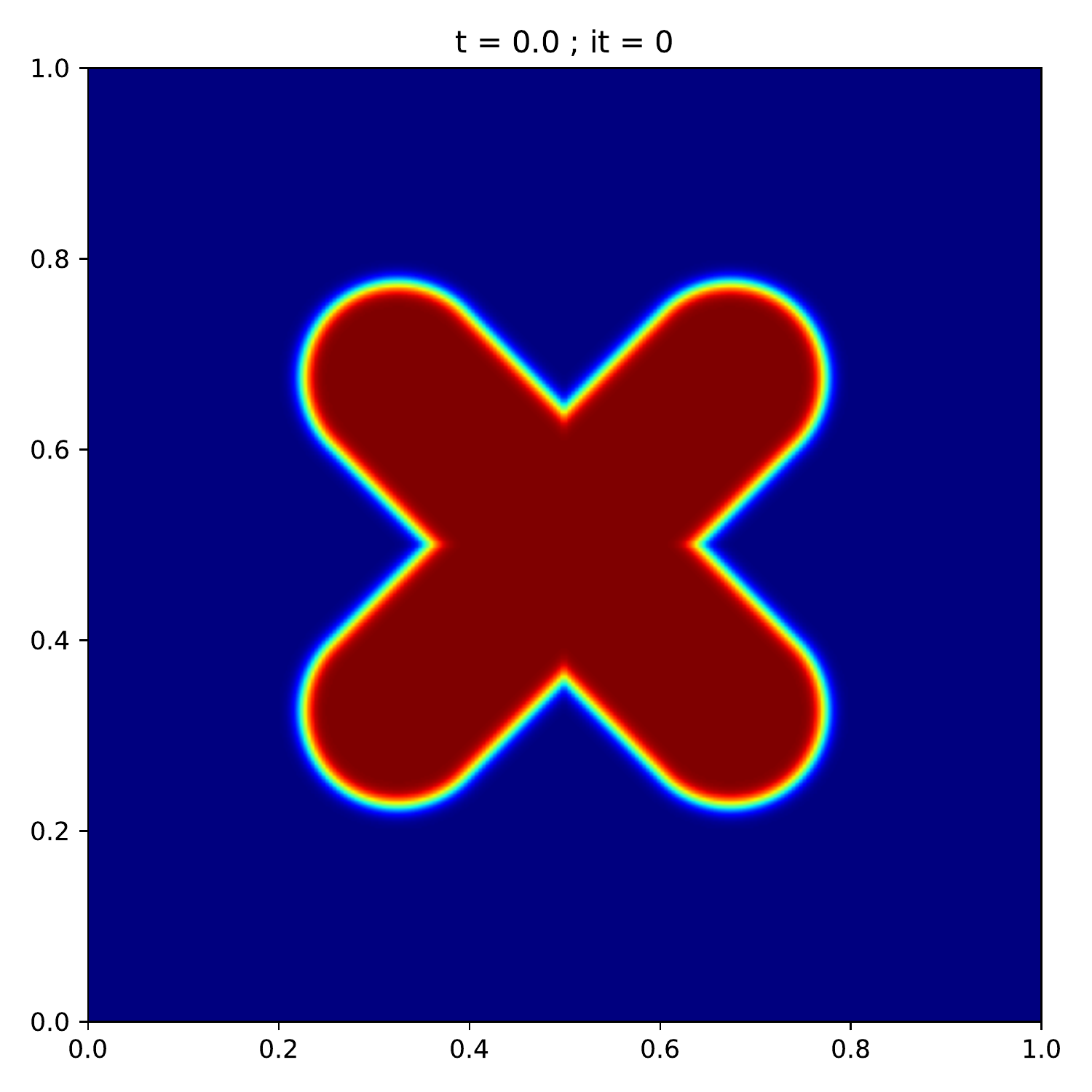}
         \includegraphics[width=0.2\textwidth]{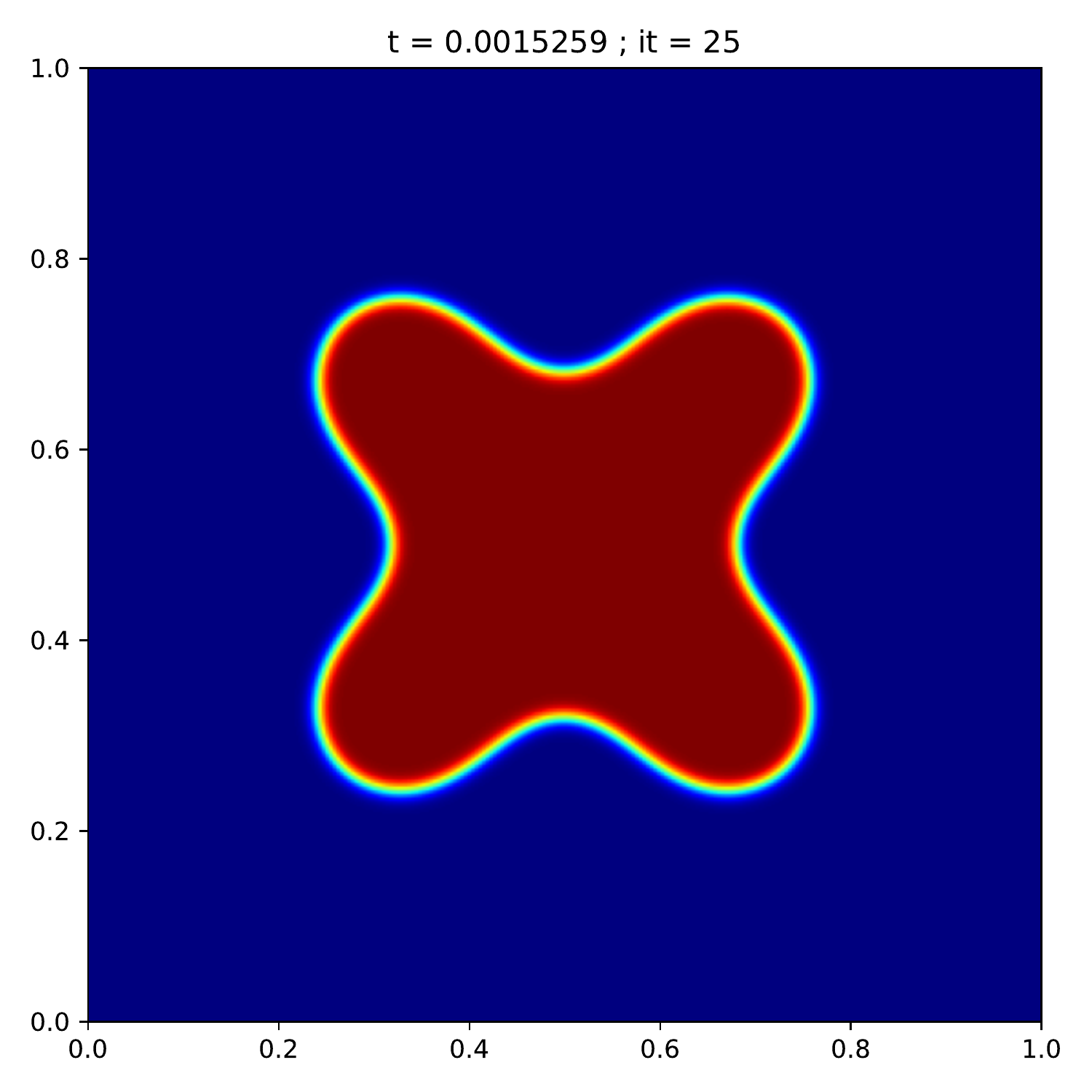}
         \includegraphics[width=0.2\textwidth]{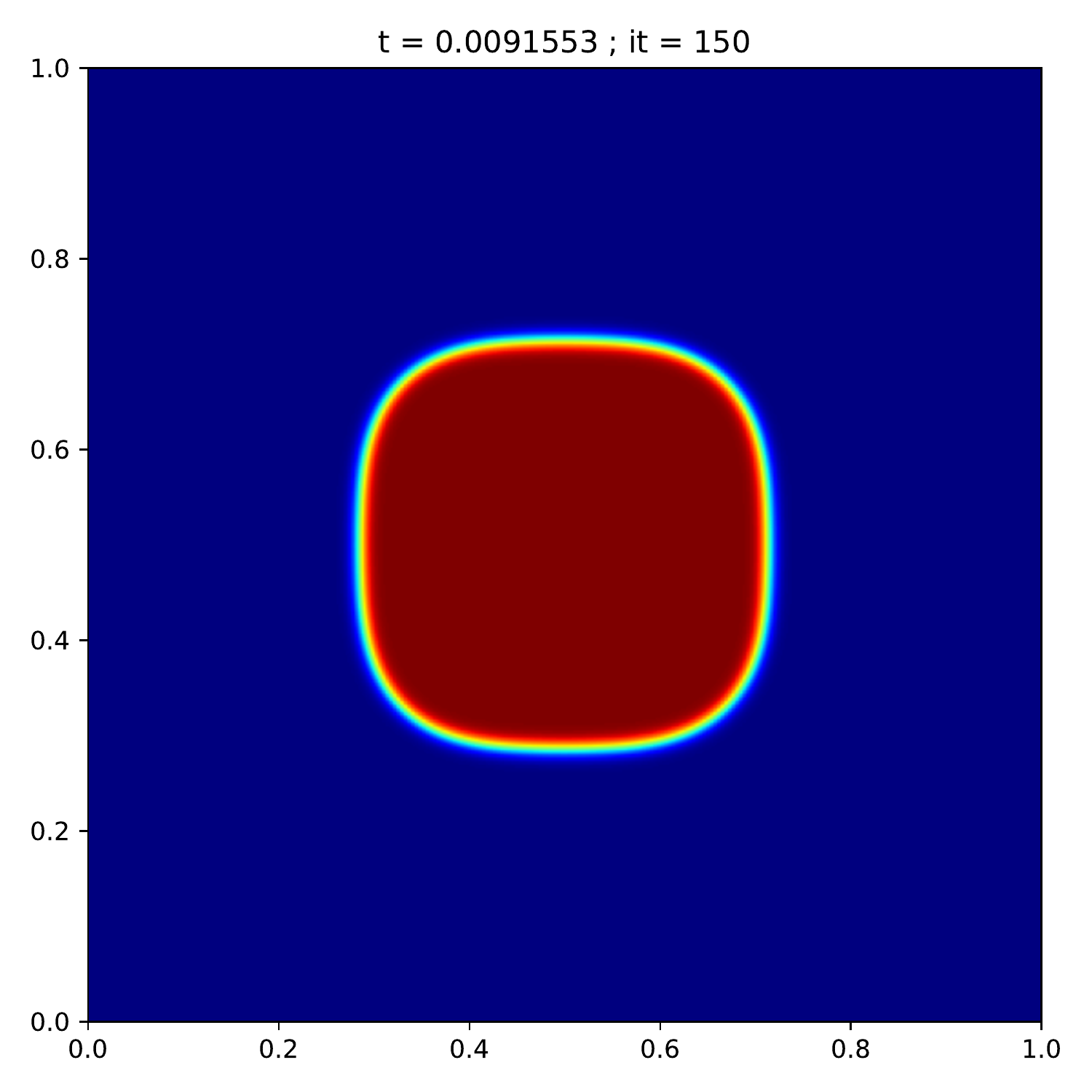}
         \includegraphics[width=0.2\textwidth]{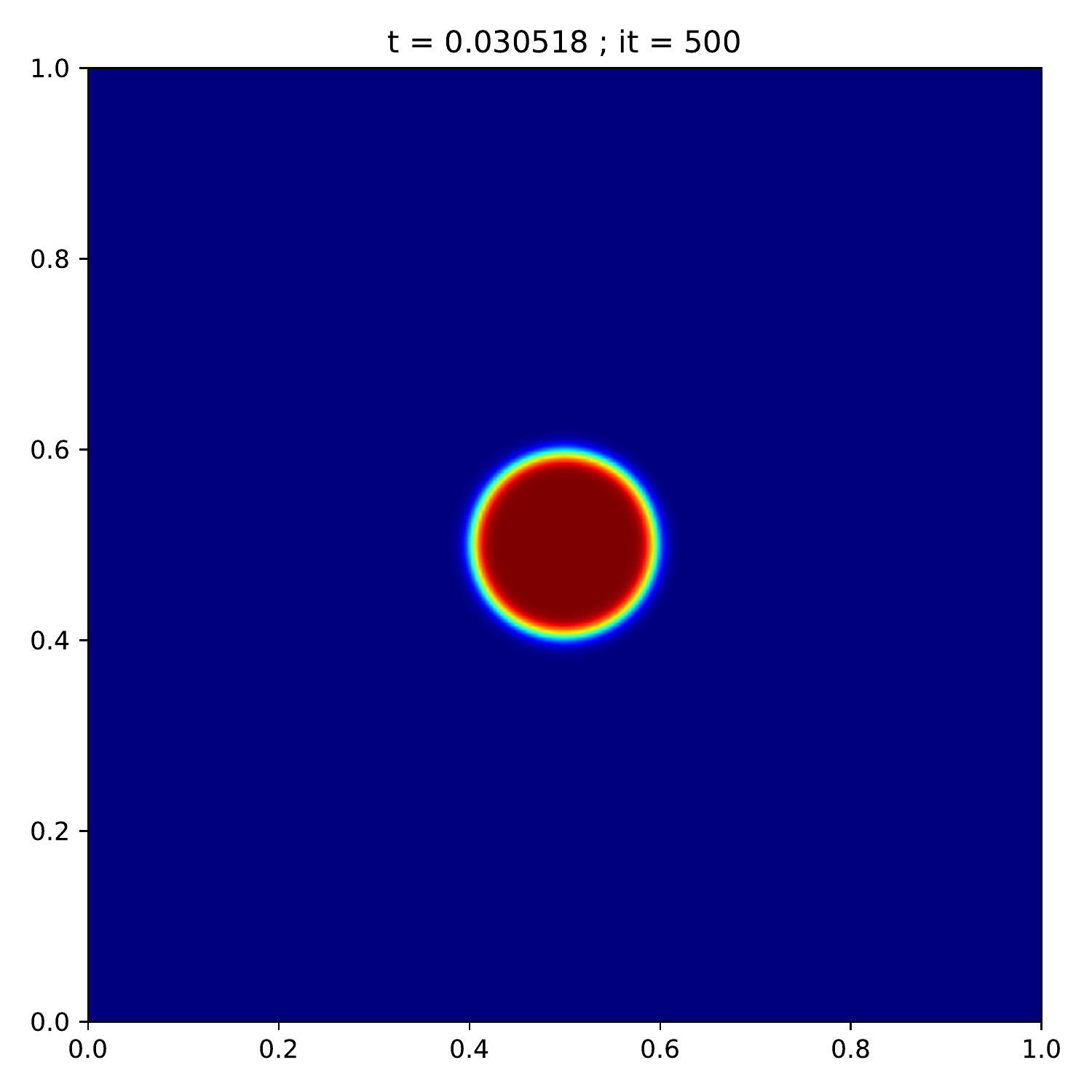}
    \caption{Comparison of numerical semigroups $\S^{\text{AC}}_{\delta_t,\varepsilon,1}$, $\S^{\text{NN}}_{\theta,1}$ and $\S^{\text{NN}}_{\theta,2}$ to approximate the mean curvature 
    flow of a non convex initial set; each line corresponds respectively to the evolution obtained at different times $t$ using, respectively,
    the numerical semigroups  $\S^{\text{AC}}_{\delta_t,\varepsilon,1}$, $\S^{\text{NN}}_{\theta,1}$ and $\S^{\text{NN}}_{\theta,2}$.}
    \label{fig:valide_q_2}
\end{figure}

We plot in figure~\ref{fig:valide_q_2} a last numerical comparison in the case of a non-convex 
initial set $\Omega(0)$ and, as before, we can clearly observe that the flows are qualitatively very similar.
 
 In conclusion, the use of neural networks not only allows us to obtain good approximations of the mean 
 curvature motion but, as shown by the quantitative comparisons above, these approximations are also of better quality than what can be obtained using  a phase field 
 approach and a solution to the Allen-Cahn equation. 

 These results are therefore very encouraging and show even more the interest of neural networks for phase field approximations in 
 the even worser case where the models would not be as accurate as the Allen-Cahn equation or the numerical schemes would 
 not be as efficient or accurate.

\subsection{Non orientable mean curvature flow $t\mapsto\Gamma(t)$ and approximation of $S_{\delta_t,\varepsilon}^{q'}$}~\\
Let us describe some numerical results on the approximation of the $S_{\delta_t,\varepsilon}^{q'}$ semigroup, still using the same architecture for the two neural networks $\S^{\text{NN}}_{\theta,1}$ and $\S^{\text{NN}}_{\theta,2}$ but training on data built from evolutions of circles in the exact non-oriented phase field representation. Let us recall here that to the limit of our knowledge,
there is no phase field model allowing to approximate such a flow by solving an Allen-Cahn-type PDE.
As before, the idea is to train both our networks on a database still made of circles evolution on a $\delta_t$ time step, 
but using the $q'$ profile instead of $q$.

\begin{figure}[htbp]
    \centering
    \includegraphics[width=1\textwidth]{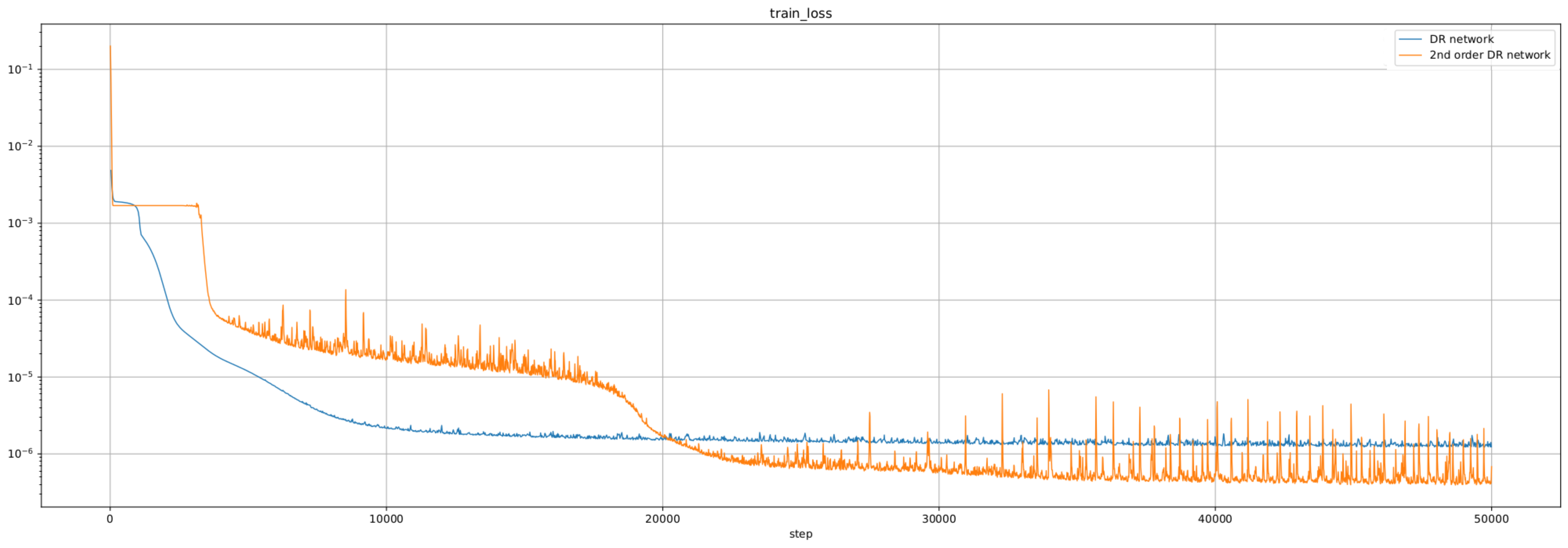}
    \caption{Learning process by optimizing the training losses $J_k$ with $k=5$ for the networks $\S^{\text{NN}}_{\theta,1}$ and $\S^{\text{NN}}_{\theta,2}$, respectively plotted in blue and orange.}
    \label{fig:Learning_process_qprim}
\end{figure}

We first plot in figure~\ref{fig:Learning_process_qprim} the evolution of the training loss energy during the training process, respectively in blue and orange for the networks $\S^{\text{NN}}_{\theta,1}$ and $\S^{\text{NN}}_{\theta,2}$. Here, both networks seem to succeed in learning the flow, although the train loss value seems much better for the second network $\S^{\text{NN}}_{\theta,2}$. 

\begin{figure}[htbp] 
    \centering
    \includegraphics[width=0.2\textwidth]{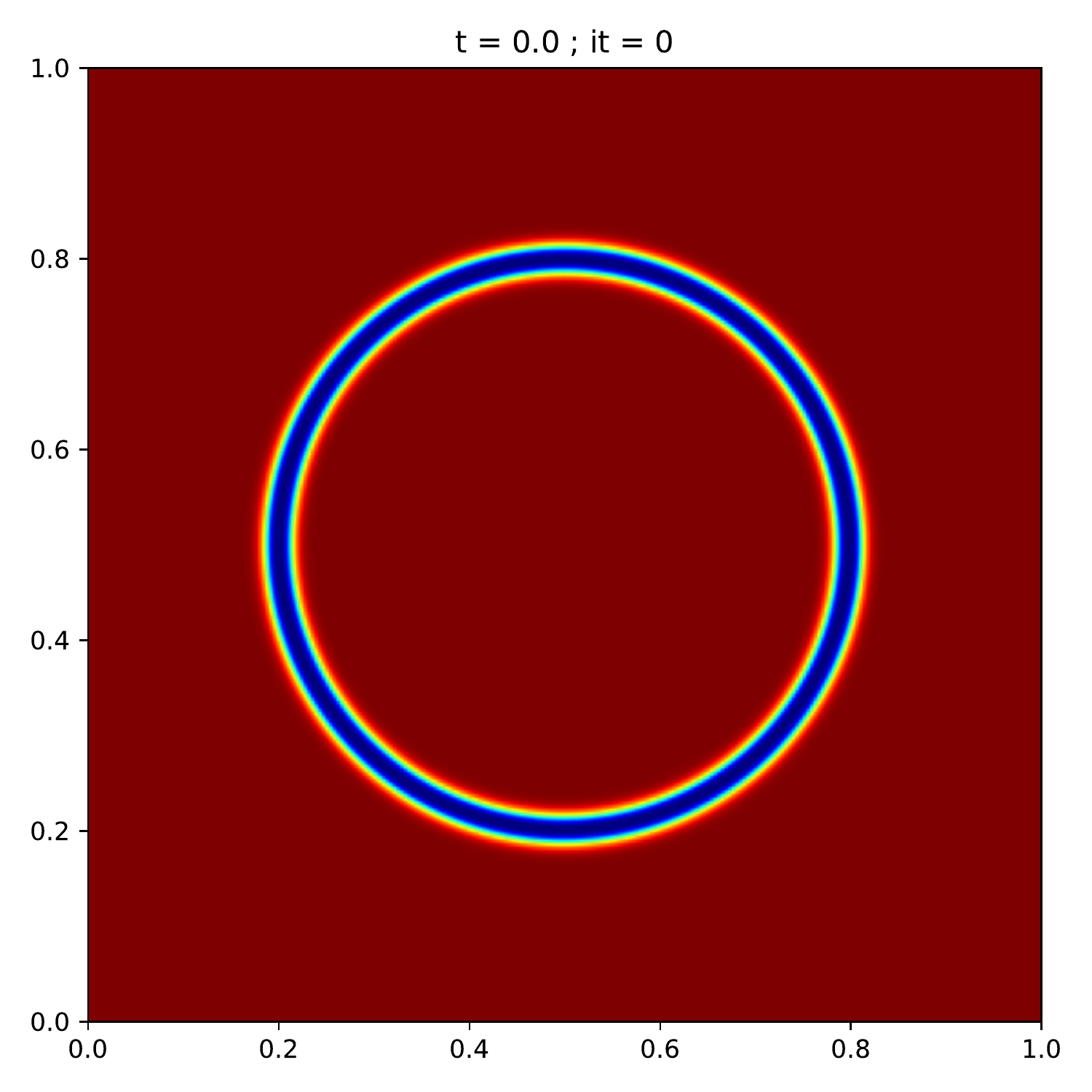}
    \includegraphics[width=0.2\textwidth]{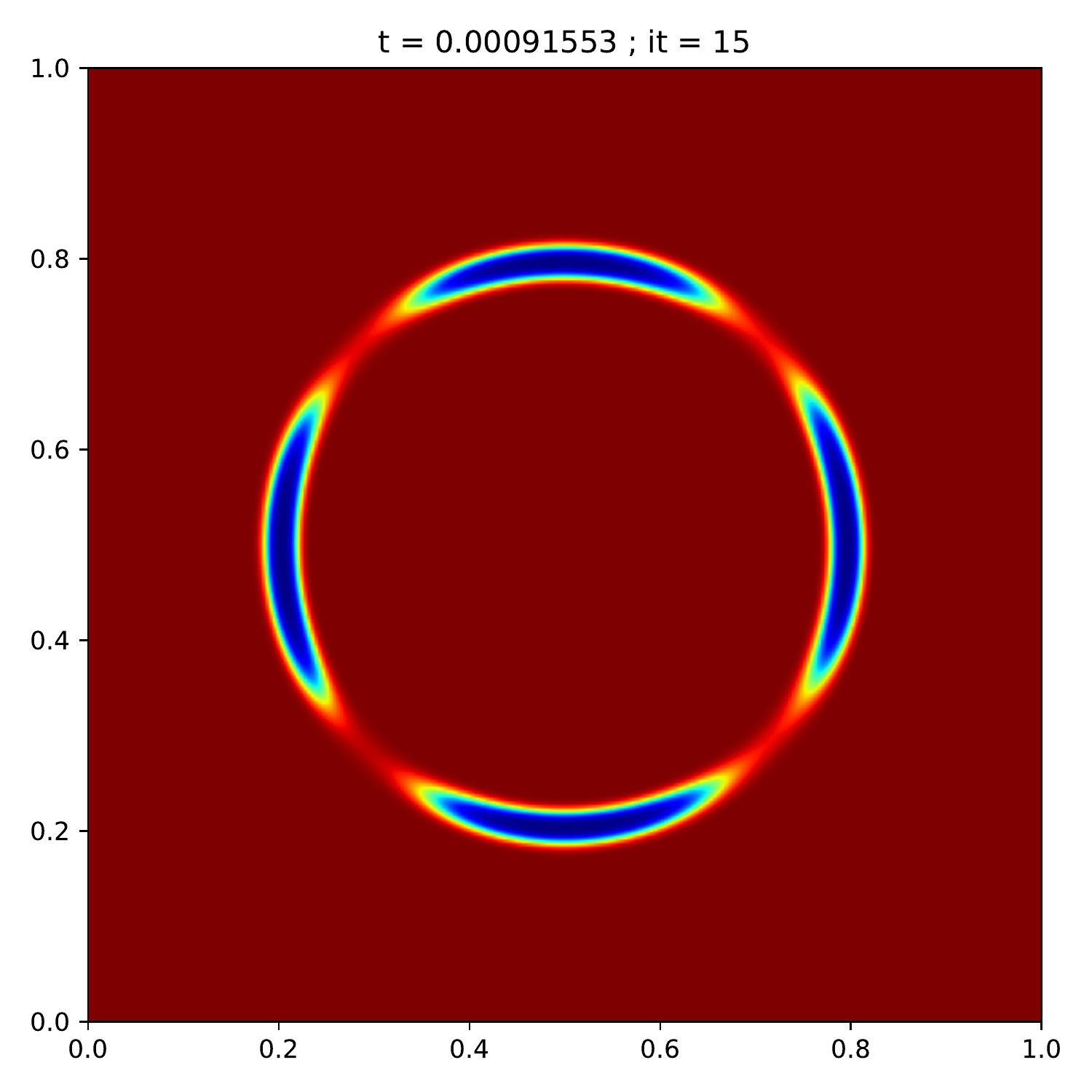}
    \includegraphics[width=0.2\textwidth]{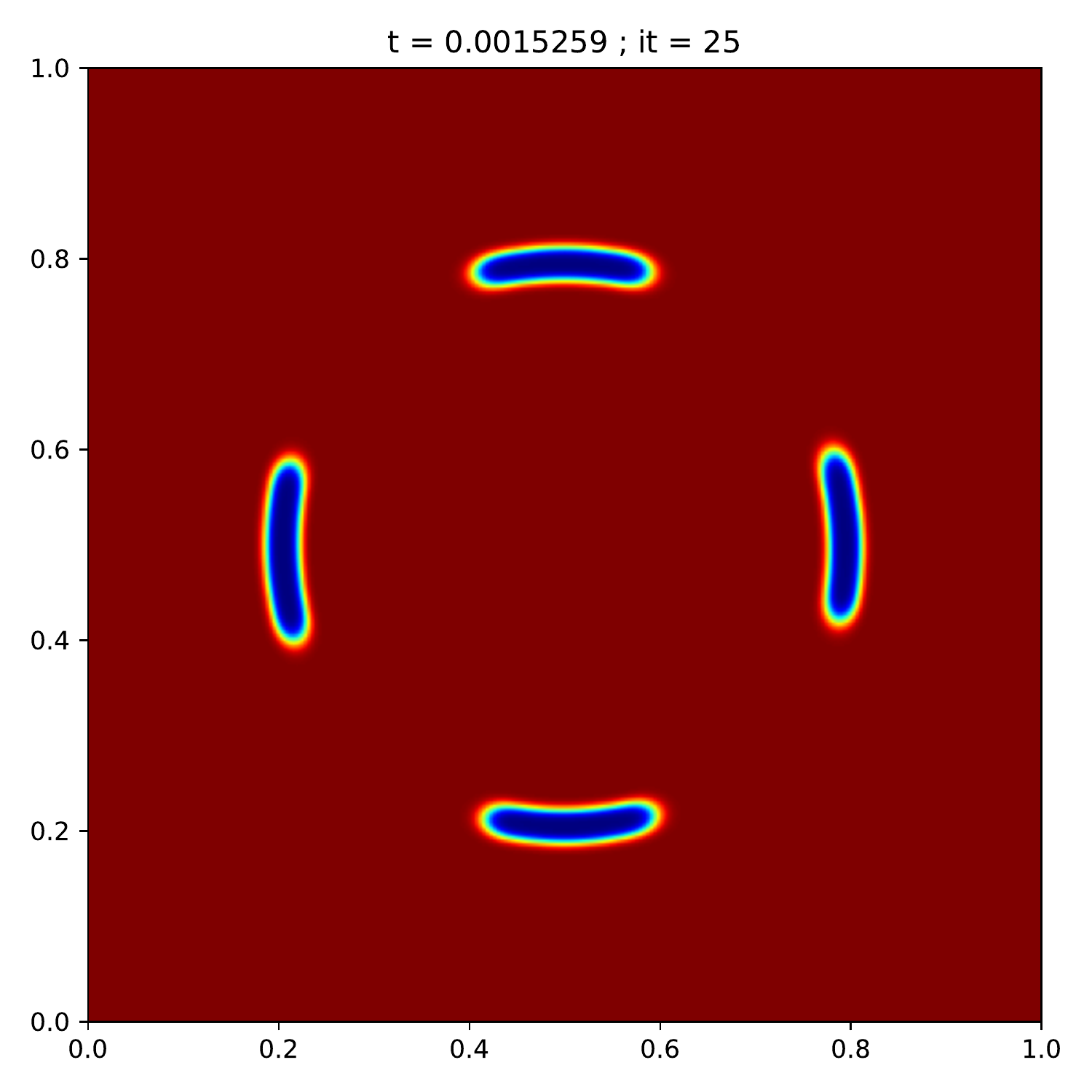}
    \includegraphics[width=0.2\textwidth]{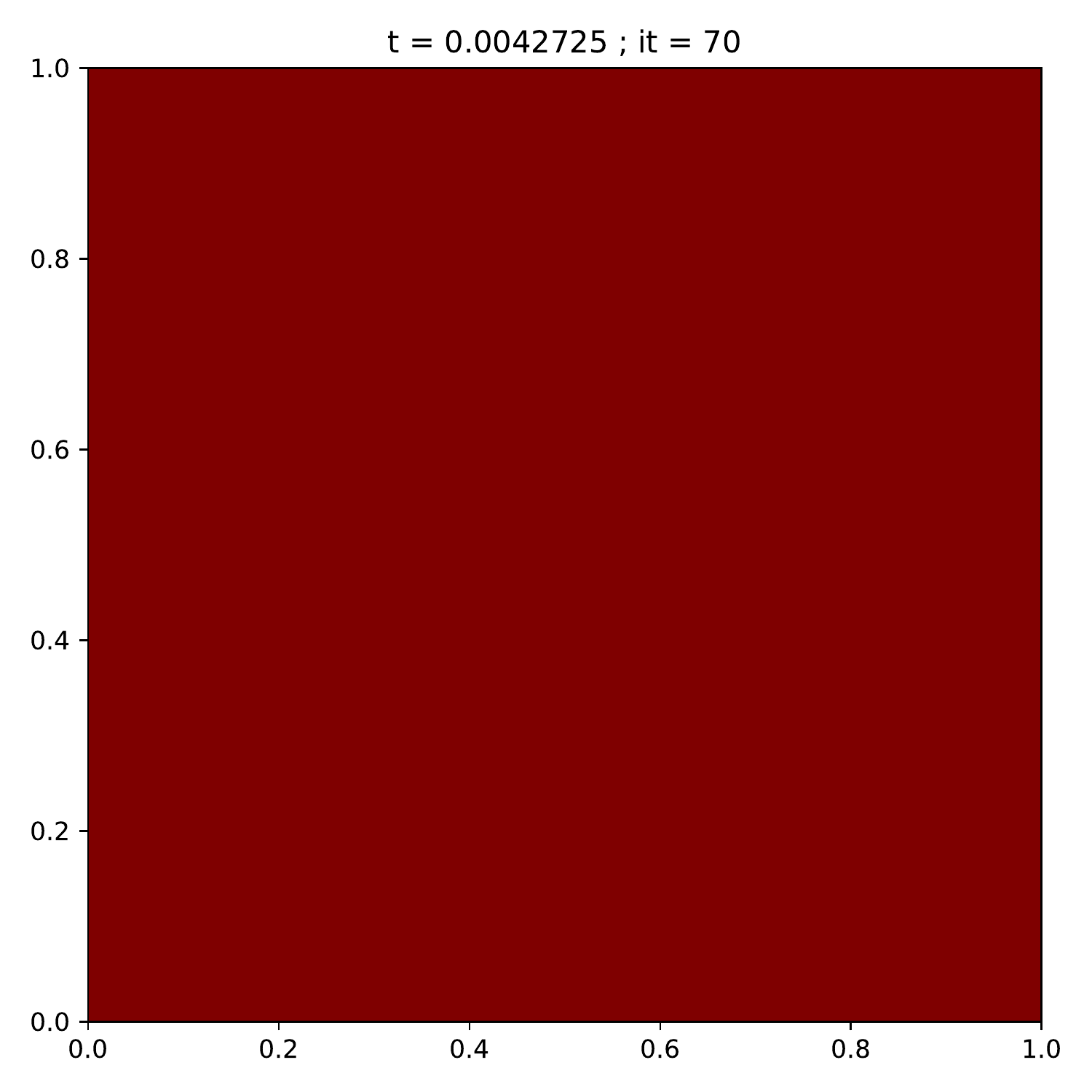} \\
    \includegraphics[width=0.2\textwidth]{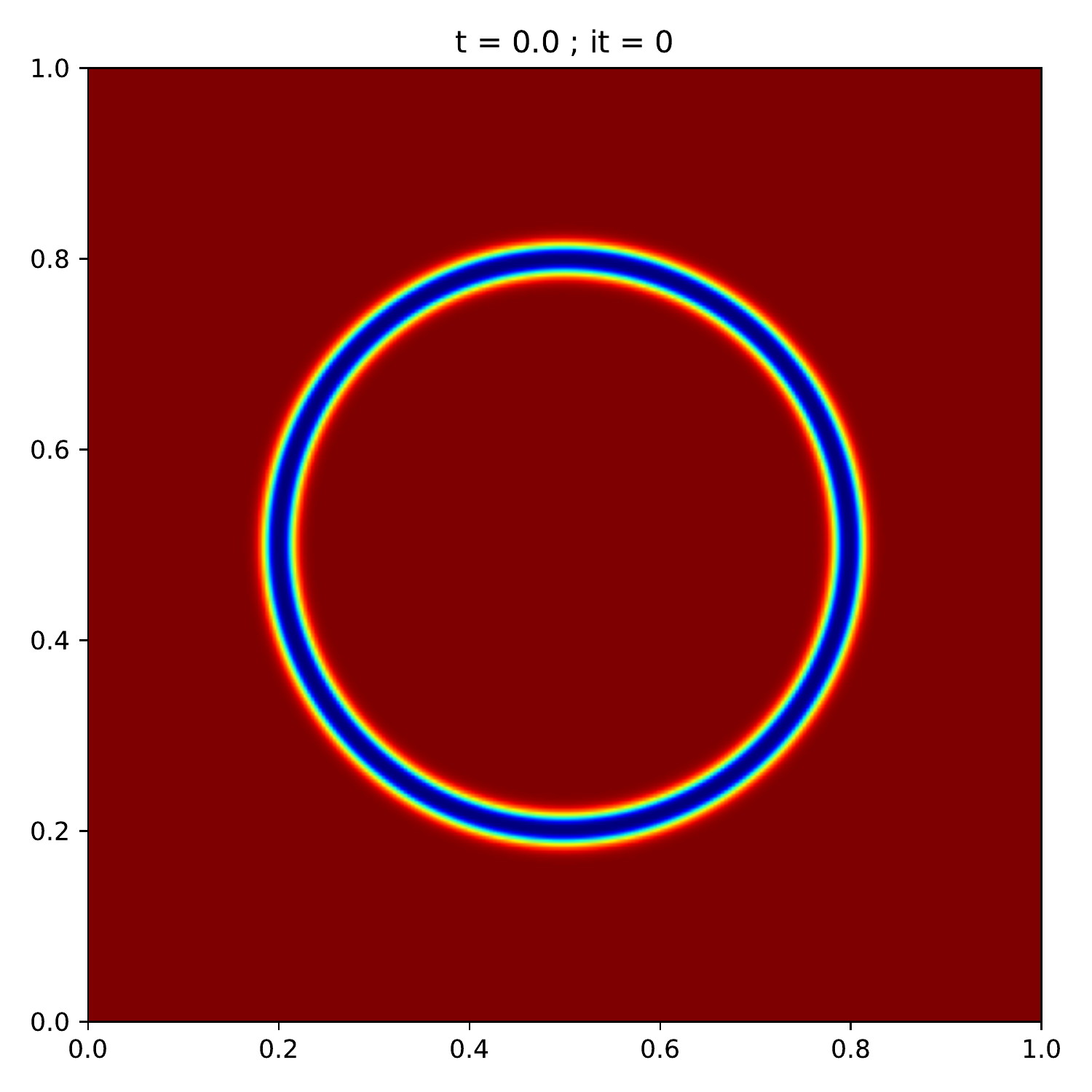}
    \includegraphics[width=0.2\textwidth]{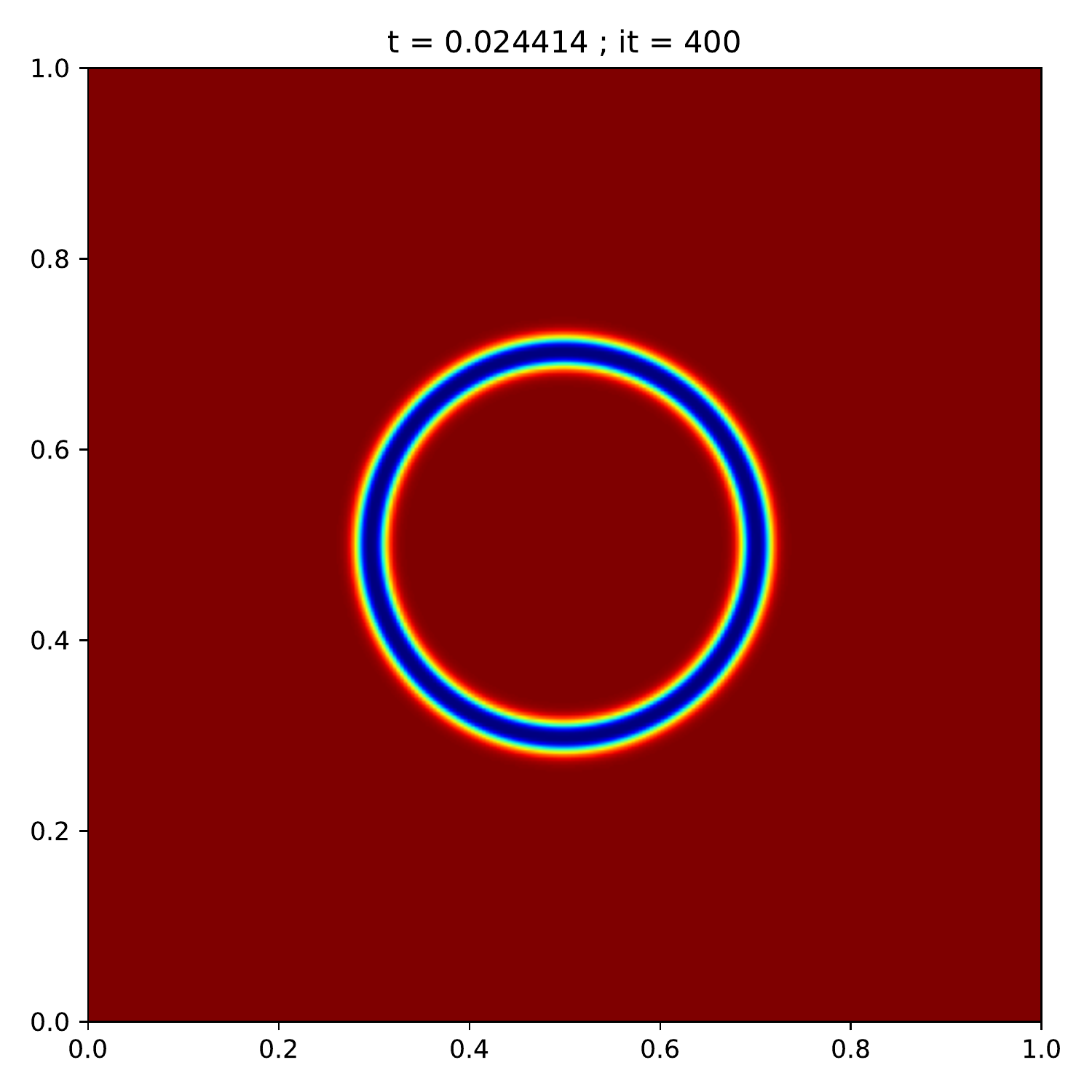}
    \includegraphics[width=0.2\textwidth]{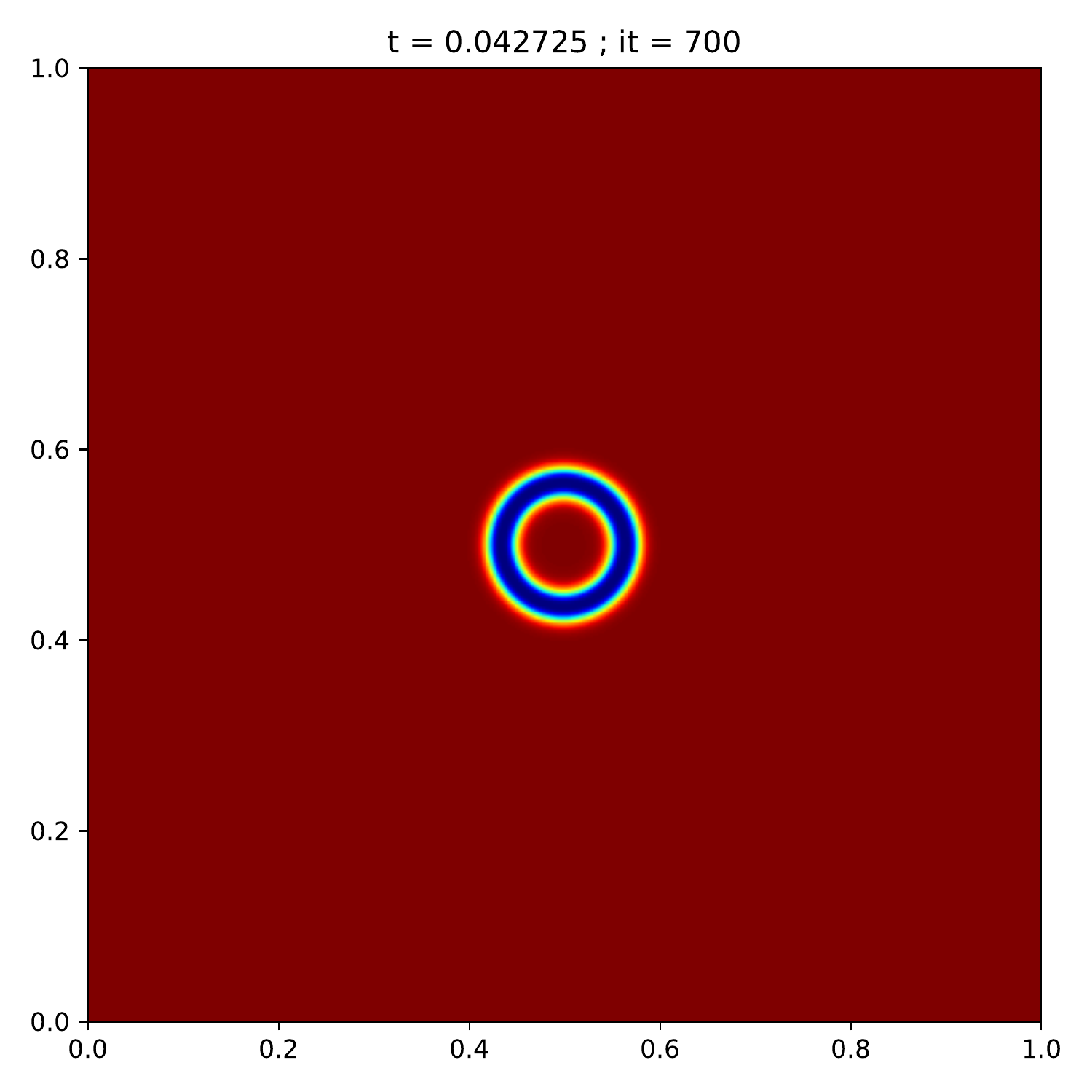}
    \includegraphics[width=0.2\textwidth]{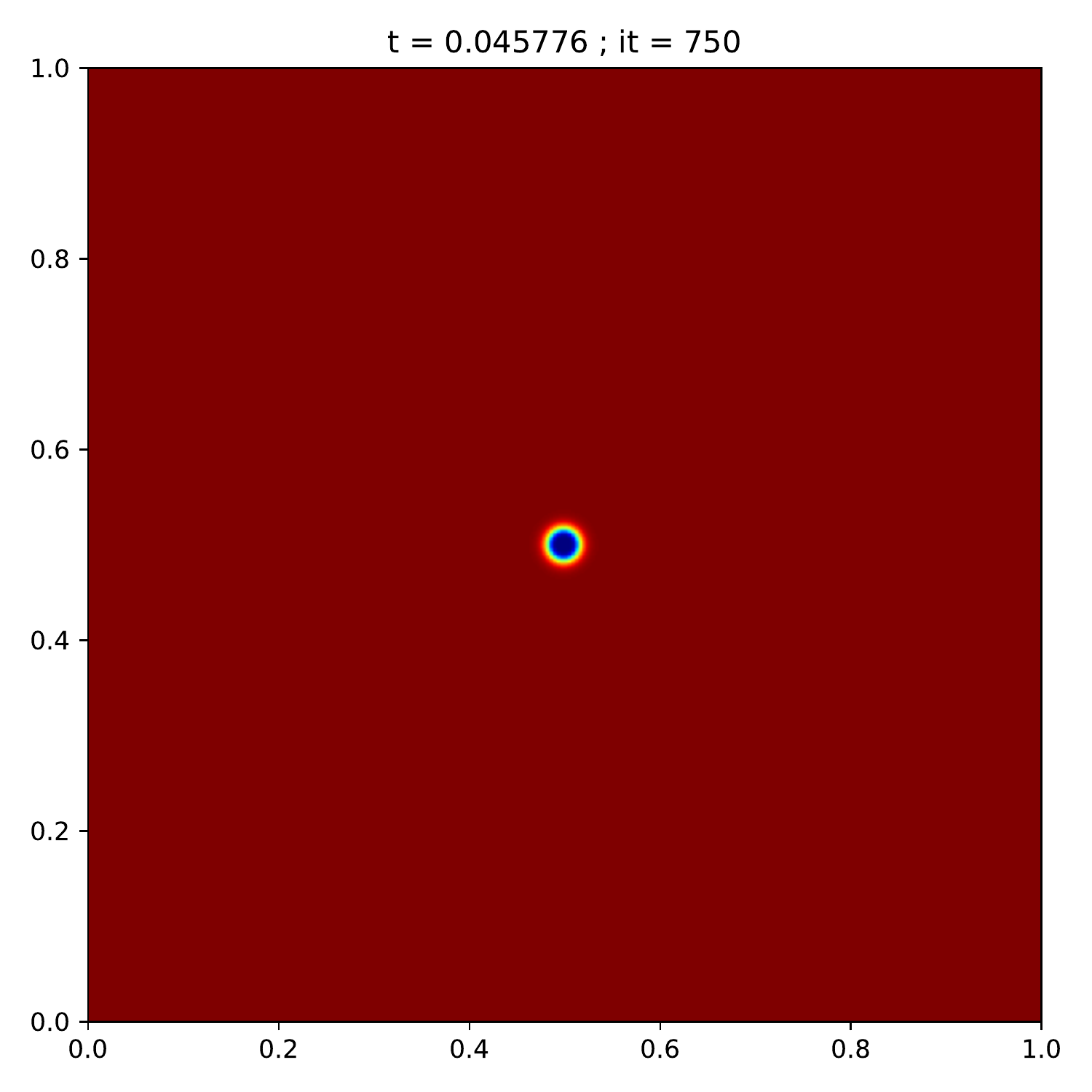}
    \caption{Numerical comparison of networks $\S^{\text{NN}}_{\theta,1}$ and $\S^{\text{NN}}_{\theta,2}$ to approximate the mean curvature 
    flow of a circle in the non oriented setting; each row corresponds to the evolution 
    obtained at different times using 
    $\S^{\text{NN}}_{\theta,1}$ and $\S^{\text{NN}}_{\theta,2}$, respectively.}
    \label{fig:One_non_oriented}
\end{figure}

\begin{figure}[htbp]
    \centering
    \includegraphics[width=0.6\textwidth]{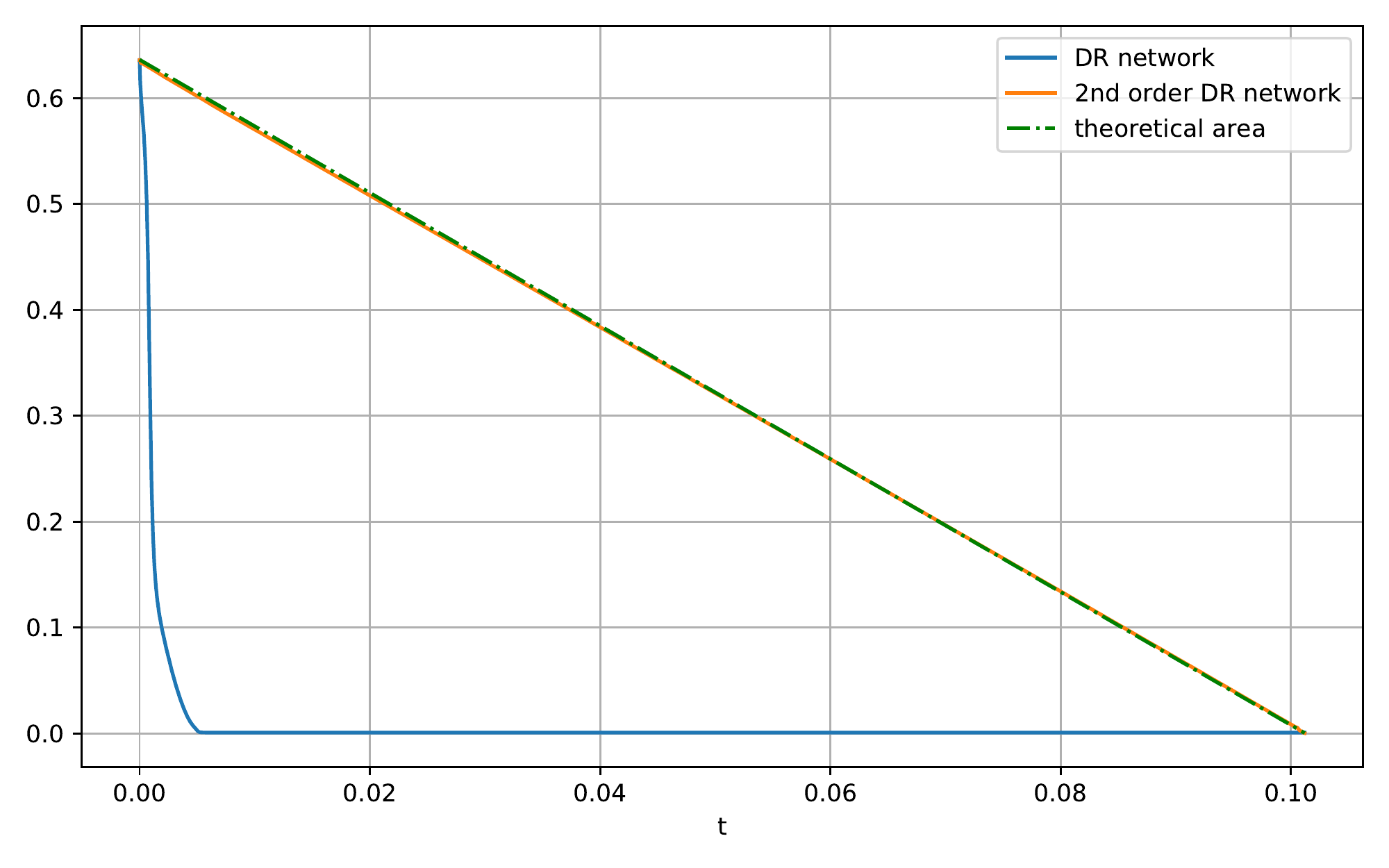}    
    \caption{Comparison of different schemes on the evolution of the radius of a circle evolving by mean curvature flow. Plots in blue, orange  correspond to $ n \mapsto \pi (\frac{1}{ 2 \pi \varepsilon} \int u^{n}(x) dx)^2$ along iterations using respectively $\S^{\text{NN}}_{\theta,1}$ and $\S^{\text{NN}}_{\theta,2}$.
    In green, the theoretical evolution of the circle area $ n \mapsto  \pi (R_0^2 - 2 n \delta_t) $. }
    \label{fig:non_oriented_radius_validation}
\end{figure}

In contrast, when we test both networks to approximate the motion by the mean curvature of a circle,
we clearly observe on figure~\ref{fig:One_non_oriented} that the first order network  
$\S^{\text{NN}}_{\theta,1}$ does not manage  to keep the profile $q'$ in a stable way
and to decrease the radius of the circle. The good news is that the second network $\S^{\text{NN}}_{\theta,2}$ 
leads to a numerical scheme that is stable enough
to reproduce well the evolution of the circle while keeping the $q'$ profile along the evolution. 
In order to have a more quantitative criterion on the evolution of the circle, we draw on figure~\ref{fig:non_oriented_radius_validation} the evolution $ n \mapsto \pi (\frac{1}{ 2 \pi \varepsilon} \int u^{n}(x) dx)^2$ which should correspond to the evolution of the area of a circle of radius $R(t)$ given by $t \mapsto \pi (R_0^2 - 2t)$. 
This is indeed clearly the case  on figure~\ref{fig:non_oriented_radius_validation} which
shows that the obtained evolution law corresponds well to the motion by mean curvature.

In order to test the reliability of the network $\S^{\text{NN}}_{\theta,2}$ to approximate 
 motions by mean curvature in a general way,
we also compare the evolution obtained with the scheme deriving from our network 
with the one deriving from the discretization  of the Allen-Cahn equation.

\begin{figure}[htbp]
    \centering
             \includegraphics[width=0.2\textwidth]{{figures/ResidualReaction_it0_allenCahn3}}
         \includegraphics[width=0.2\textwidth]{{figures/ResidualReaction_it25_allenCahn3}}
         \includegraphics[width=0.2\textwidth]{{figures/ResidualReaction_it150_allenCahn3}}
         \includegraphics[width=0.2\textwidth]{{figures/ResidualReaction_it500_allenCahn3}} \\
    \includegraphics[width=0.2\textwidth]{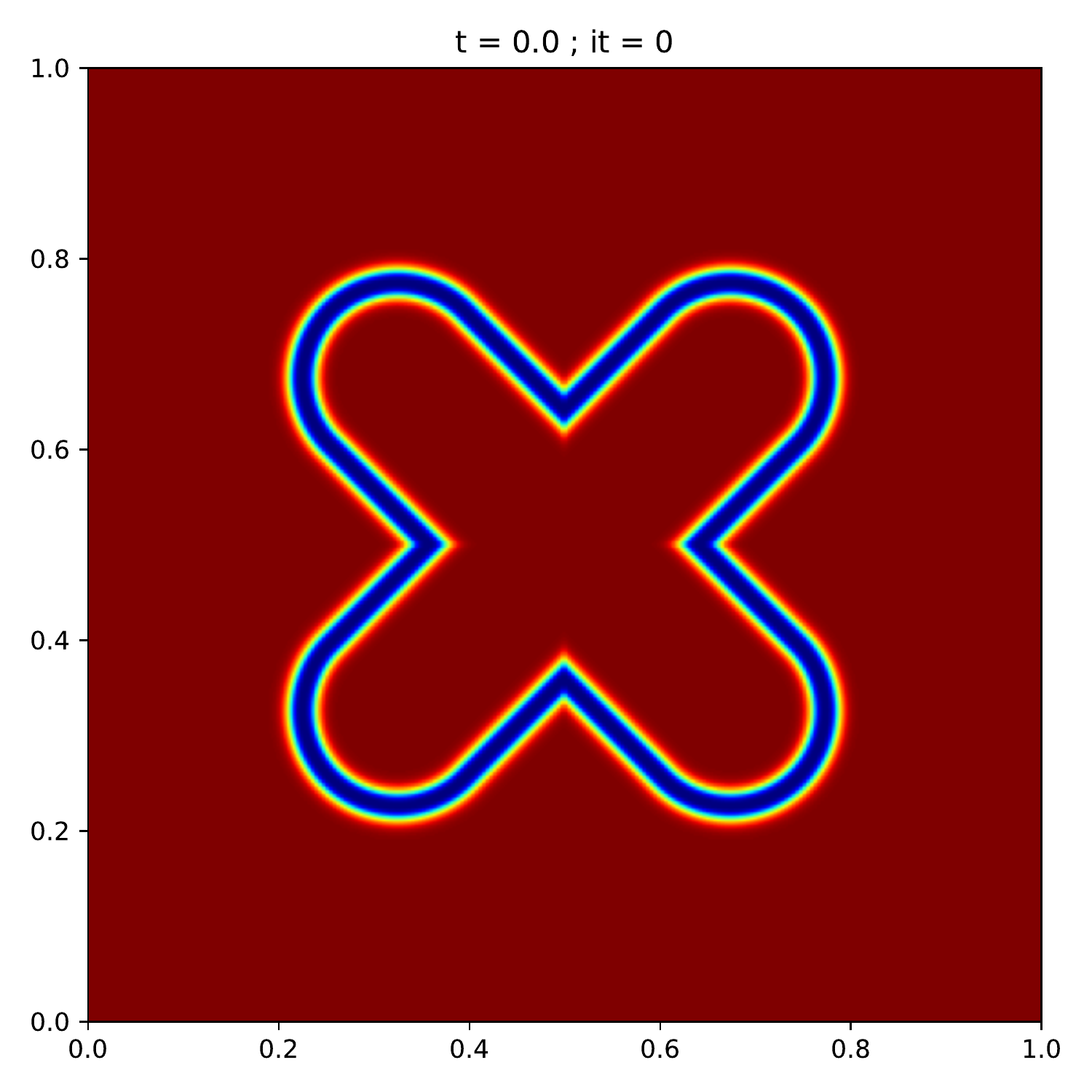}
    \includegraphics[width=0.2\textwidth]{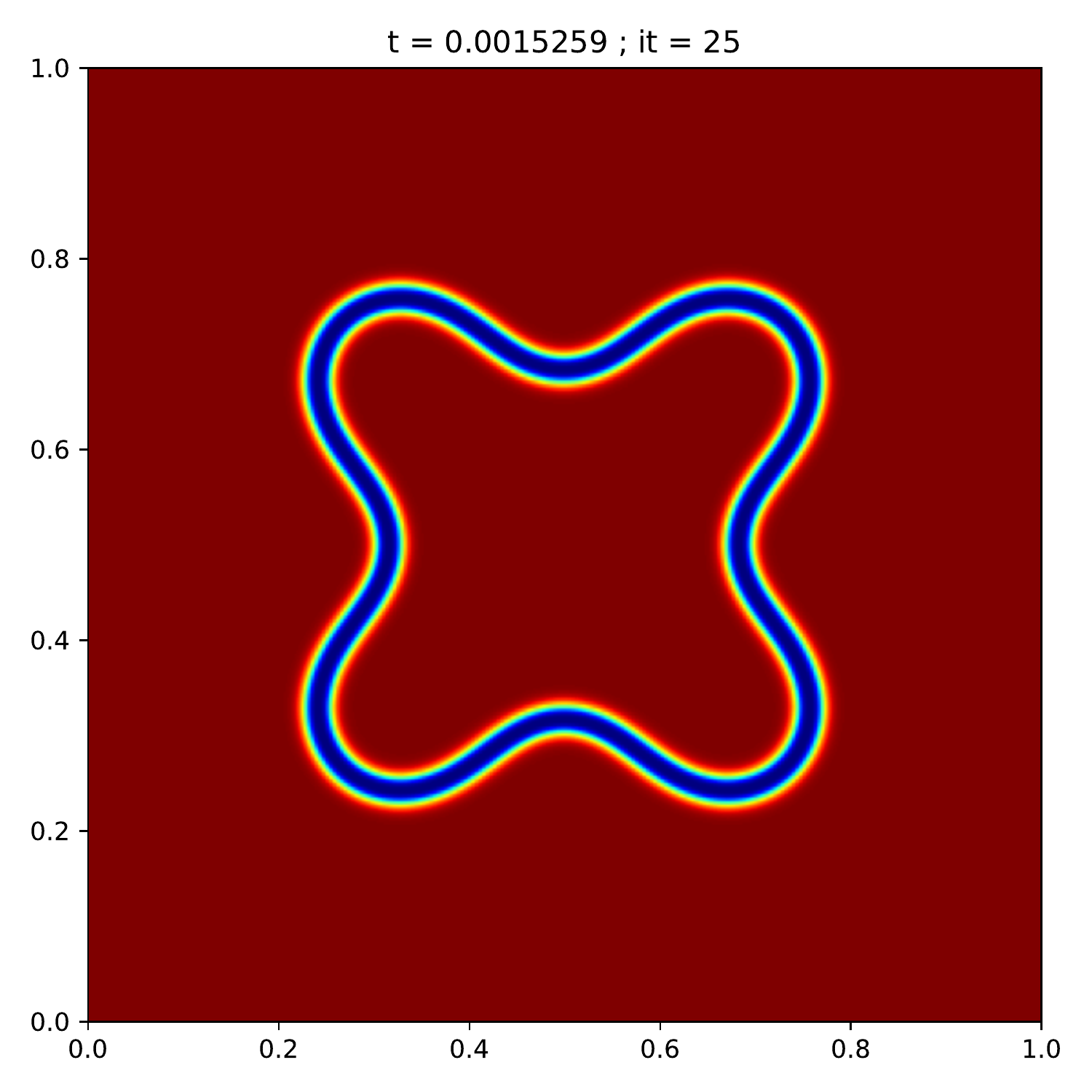}
    \includegraphics[width=0.2\textwidth]{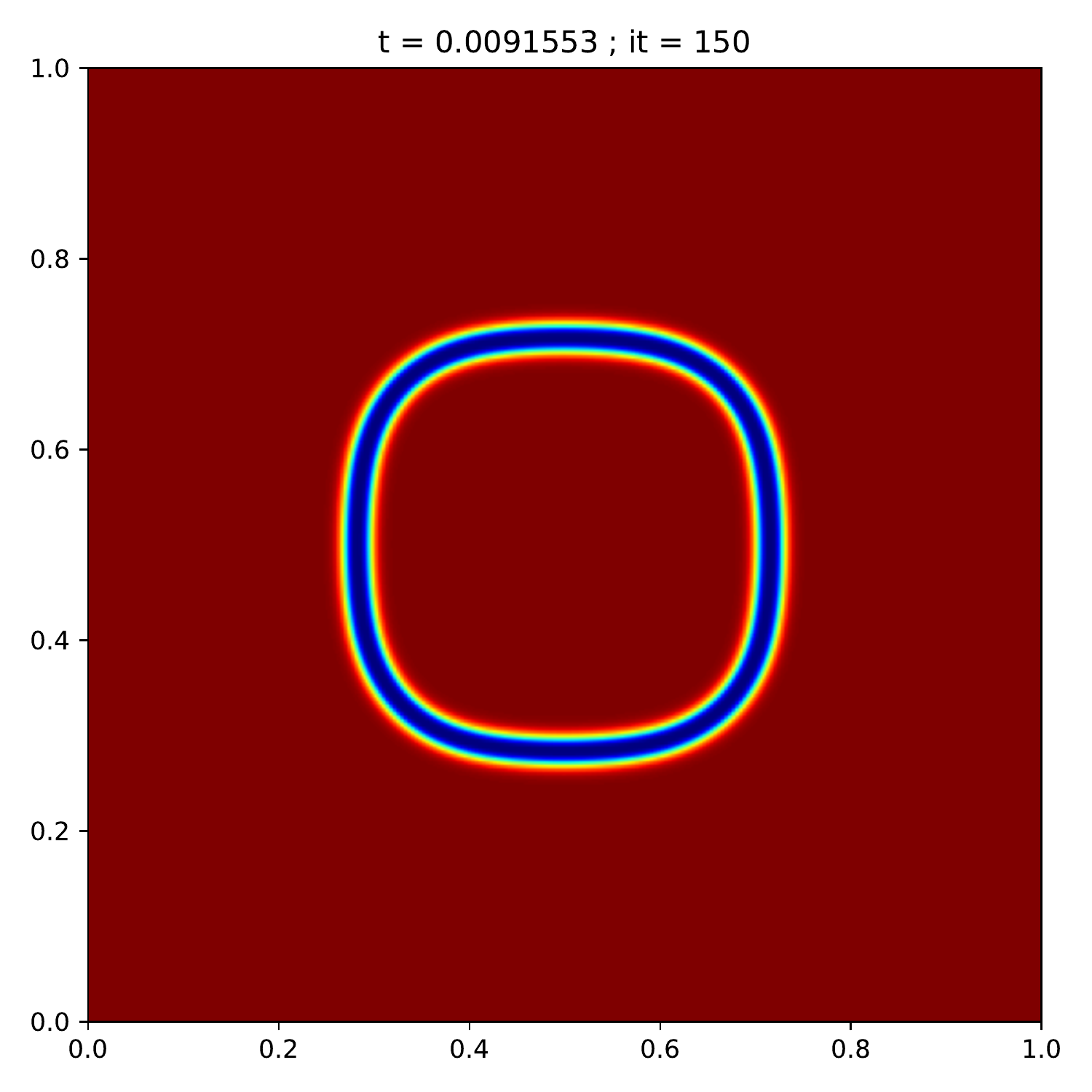}
    \includegraphics[width=0.2\textwidth]{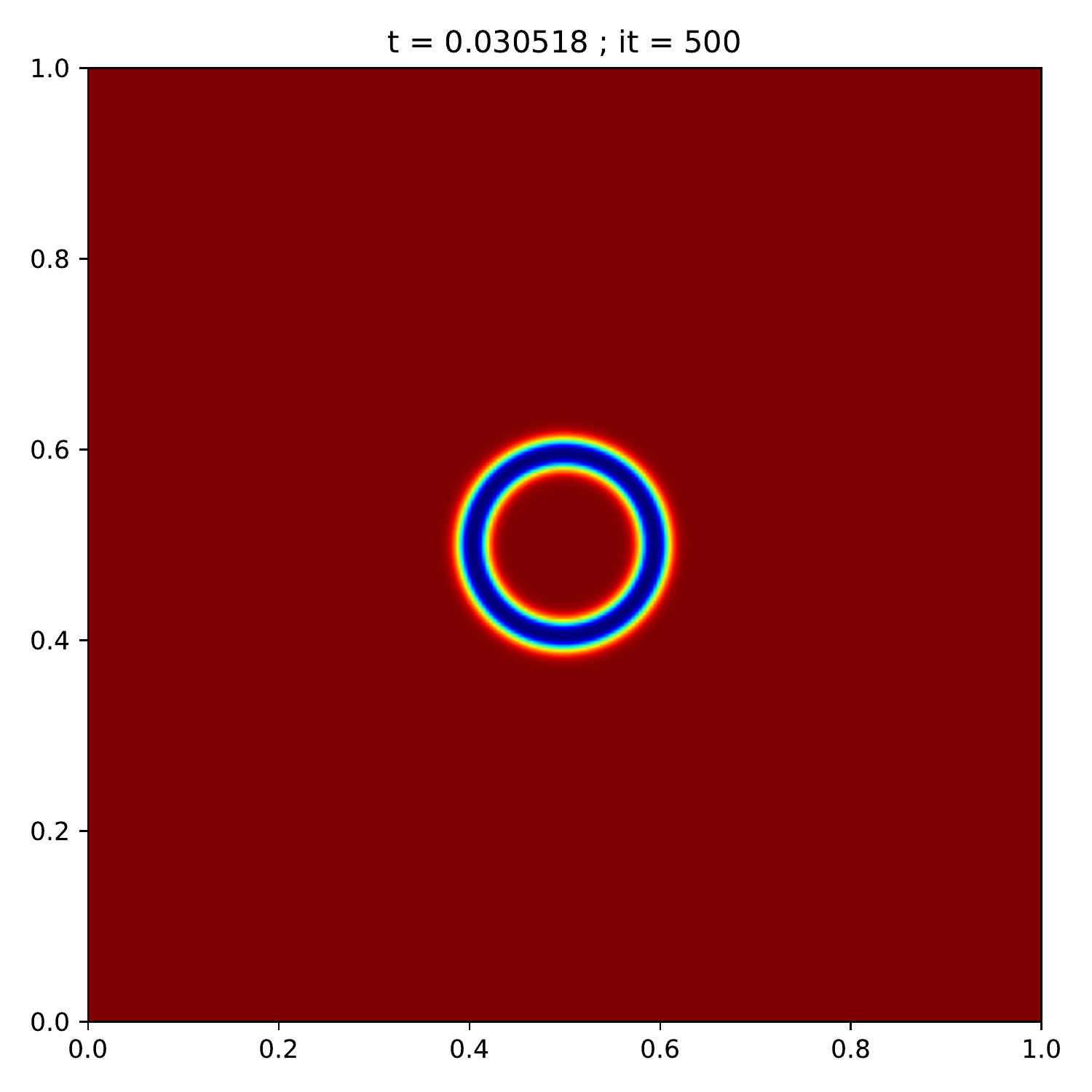}
    \caption{Approximation of the mean curvature flow of a non convex initial set: numerical comparison of the classical discretisation
     of the Allen-Cahn equation $\S^{\text{AC}}_{\delta_t = \varepsilon^2, \varepsilon,1}$ (first line)  
     with  the network $\S^{\text{NN}}_{\theta,2}$ trained on the non oriented  phase field $q'(\operatorname{dist}(x,\Gamma)/\varepsilon)$ (second line). 
     In each line, we plot the solution $u^{n}$ at different times~$t_n$.} 
    \label{fig:non_oriented_non_convex_evolution}
\end{figure}

Here, although the starting set $\Omega(0)$ is the same in both numerical experiments, the initial phase field solution $u^{0}$ 
is different and depends on the profile used. 
We plot on figure \ref{fig:non_oriented_non_convex_evolution}   the solution $u^{n}$ computed at different times $t_n$ using either the Allen-Cahn discrete semigroup $\S^{\text{AC}}_{\delta_t = \varepsilon^2, \varepsilon,1}$  or the network $\S^{\text{NN}}_{\theta,2}$. We clearly observe a similar flow that validates the methodology 
and the use of our neural networks to  approximate the motion by mean curvature in the case 
of a non-orientable set.
 
\begin{figure}[!htbp]
    \centering
    \includegraphics[width=0.2\textwidth]{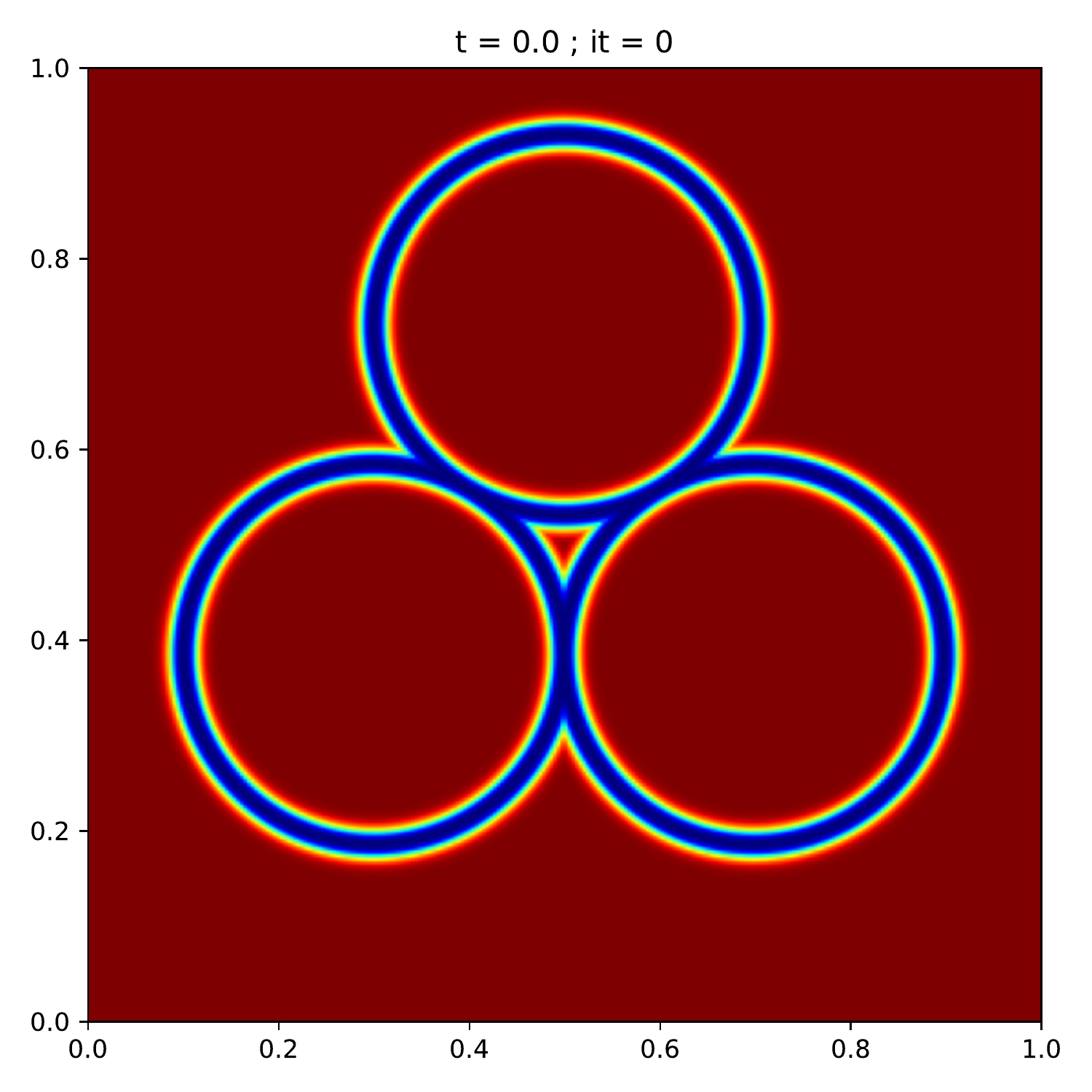}
    \includegraphics[width=0.2\textwidth]{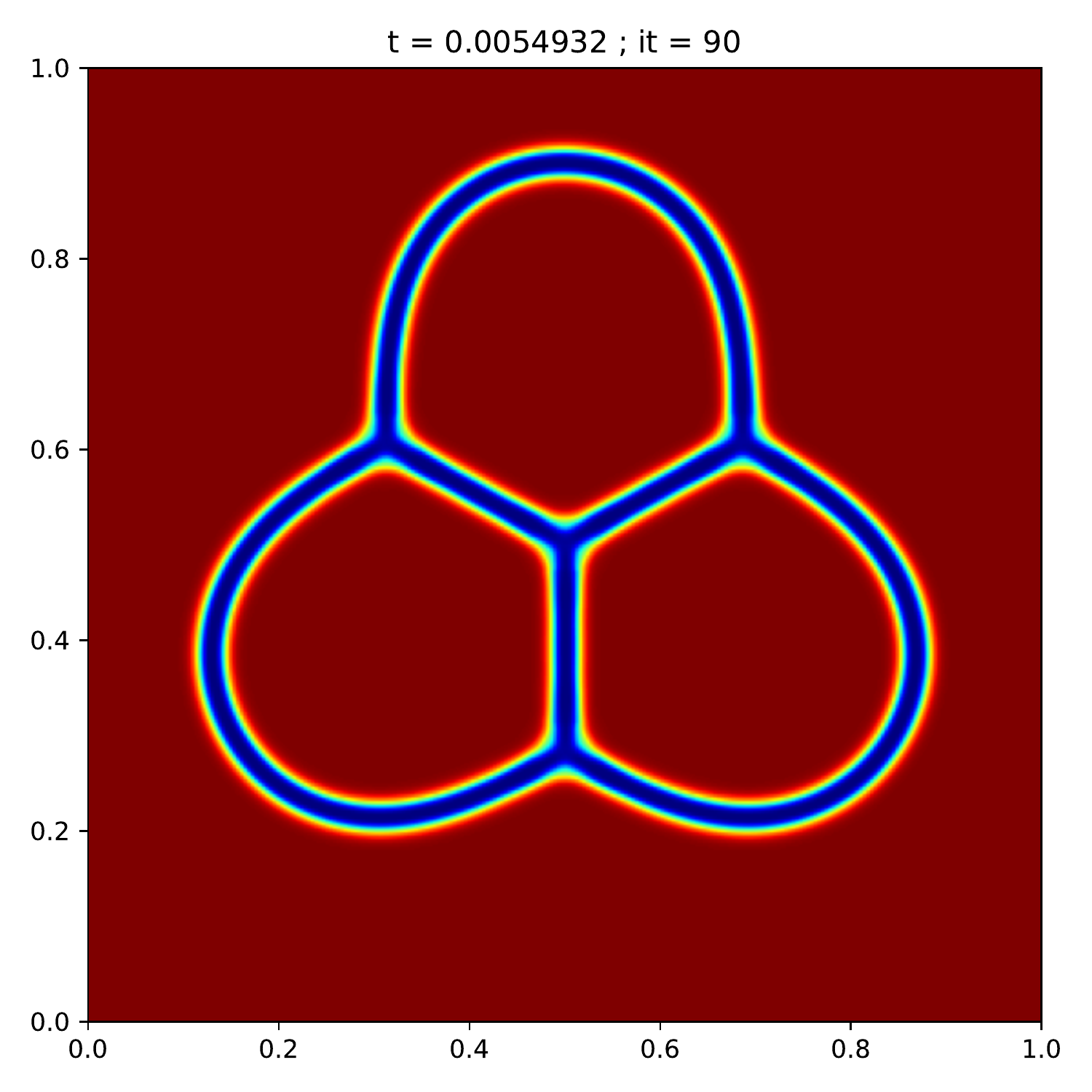}
    \includegraphics[width=0.2\textwidth]{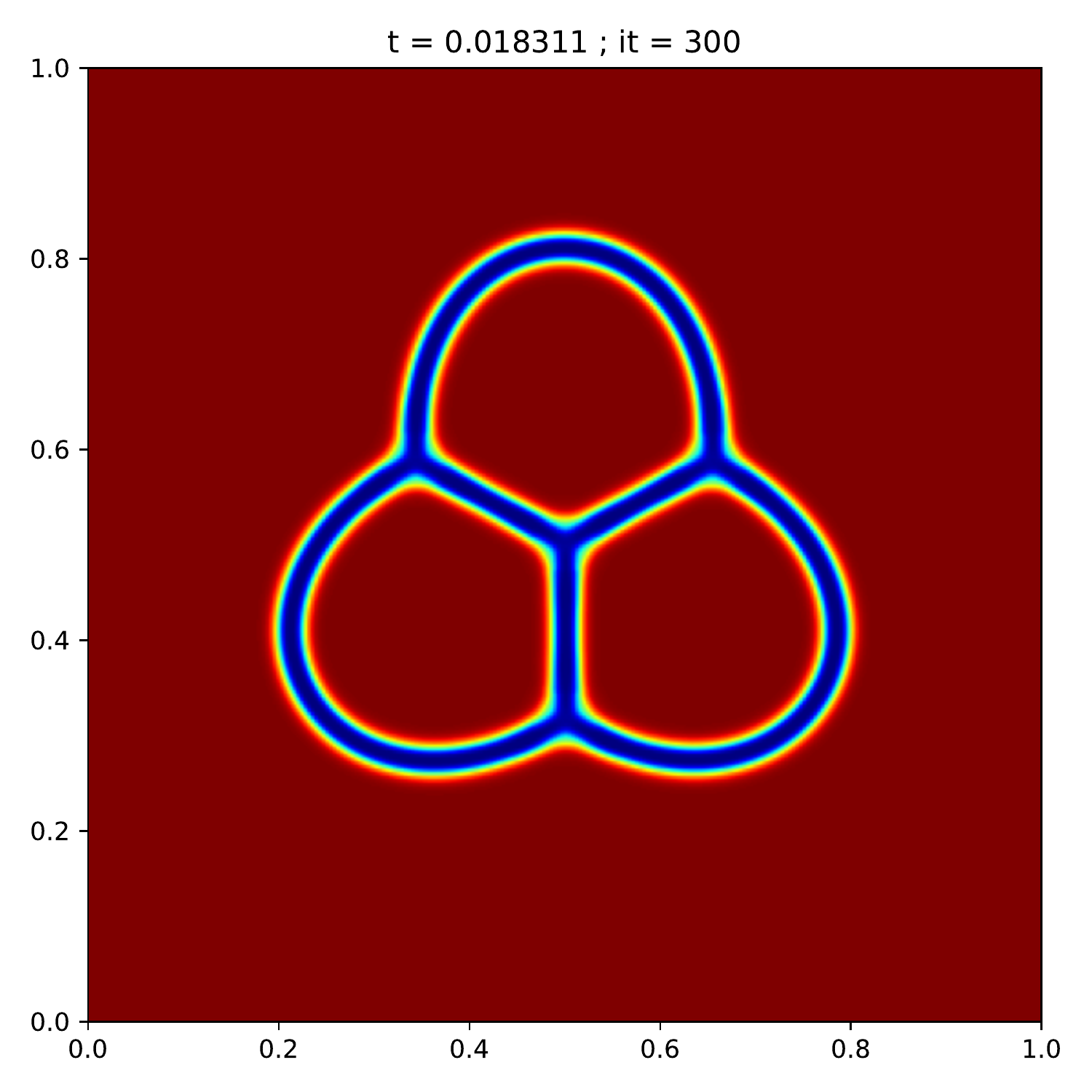}
    \includegraphics[width=0.2\textwidth]{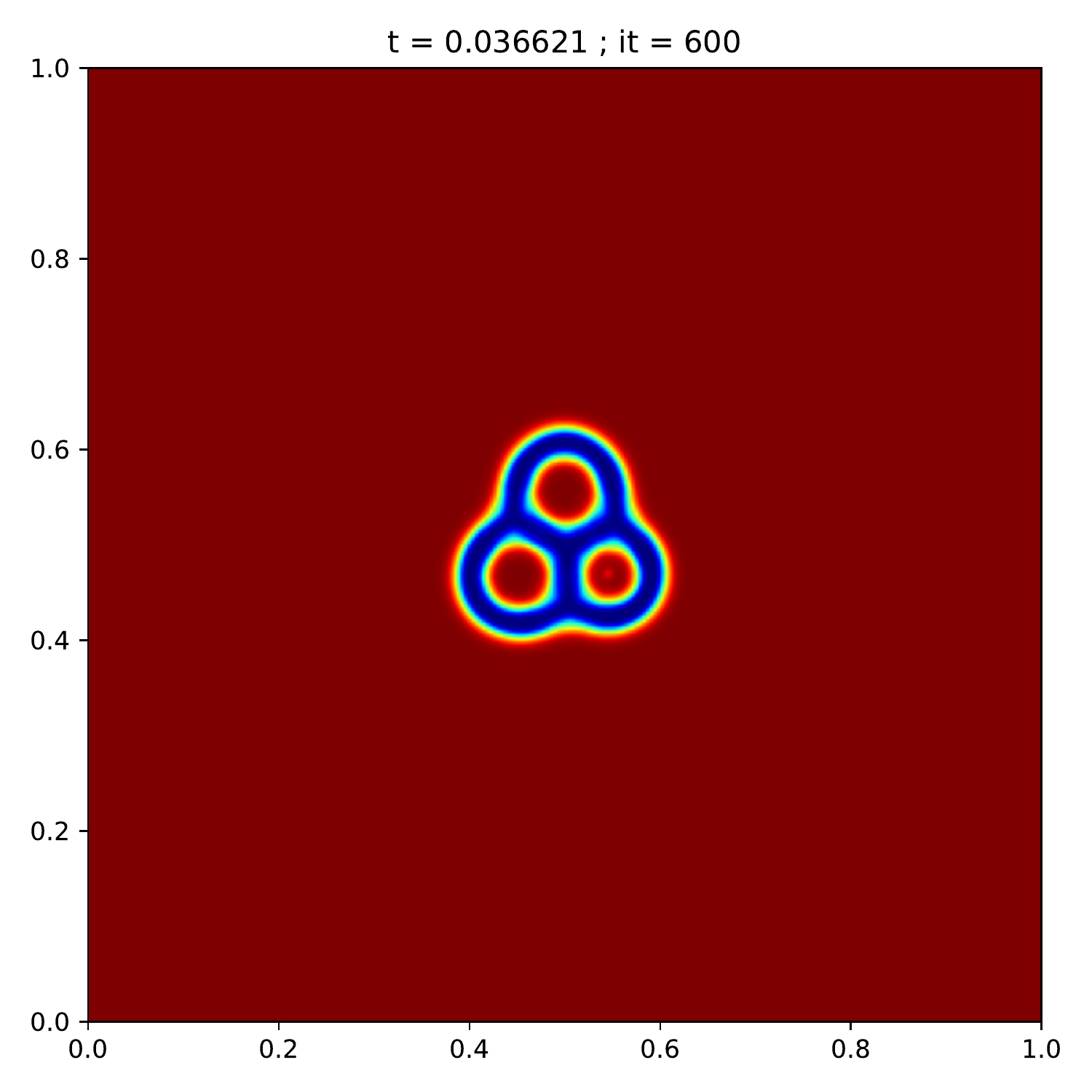}\\
   \includegraphics[width=0.2\textwidth]{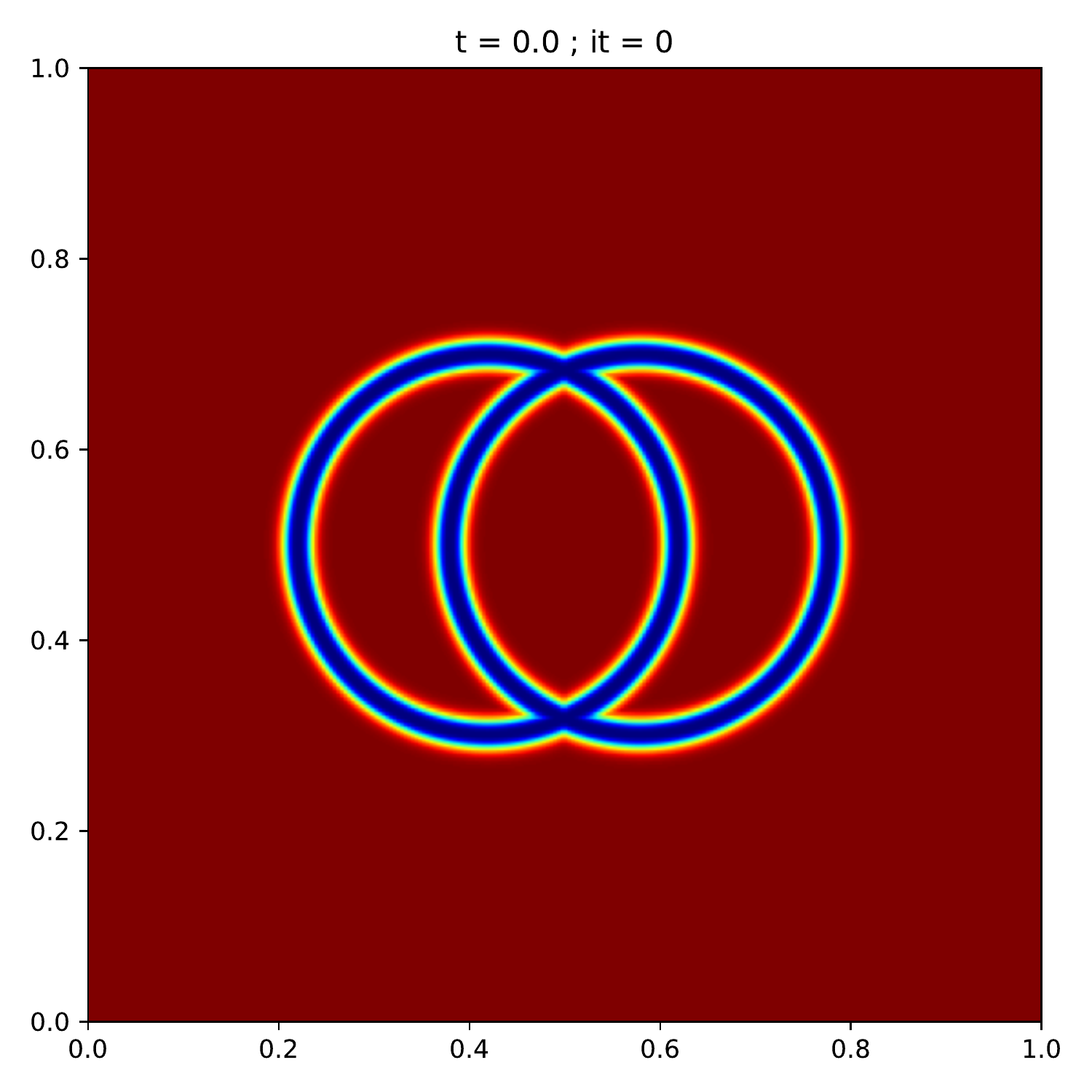}
   \includegraphics[width=0.2\textwidth]{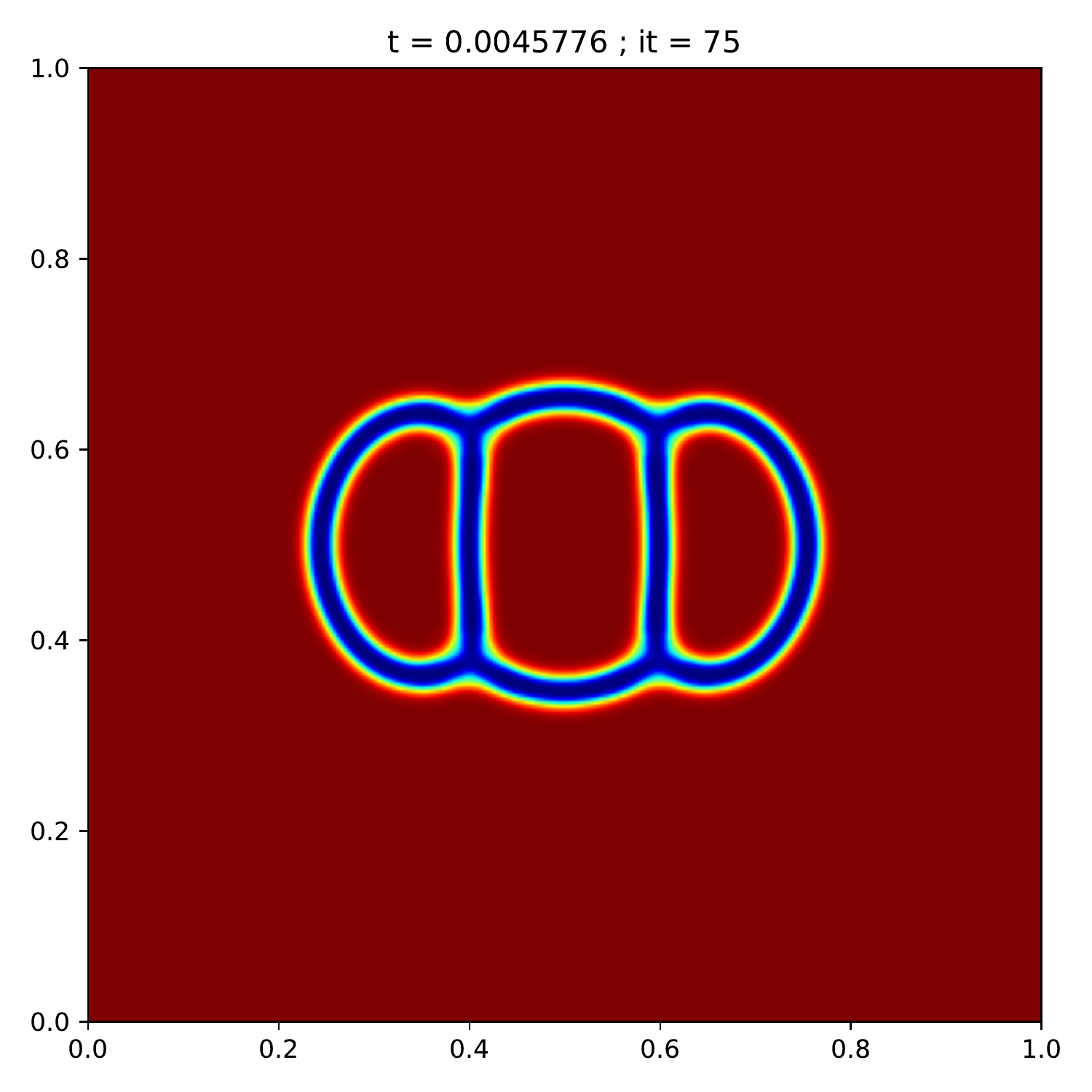}
   \includegraphics[width=0.2\textwidth]{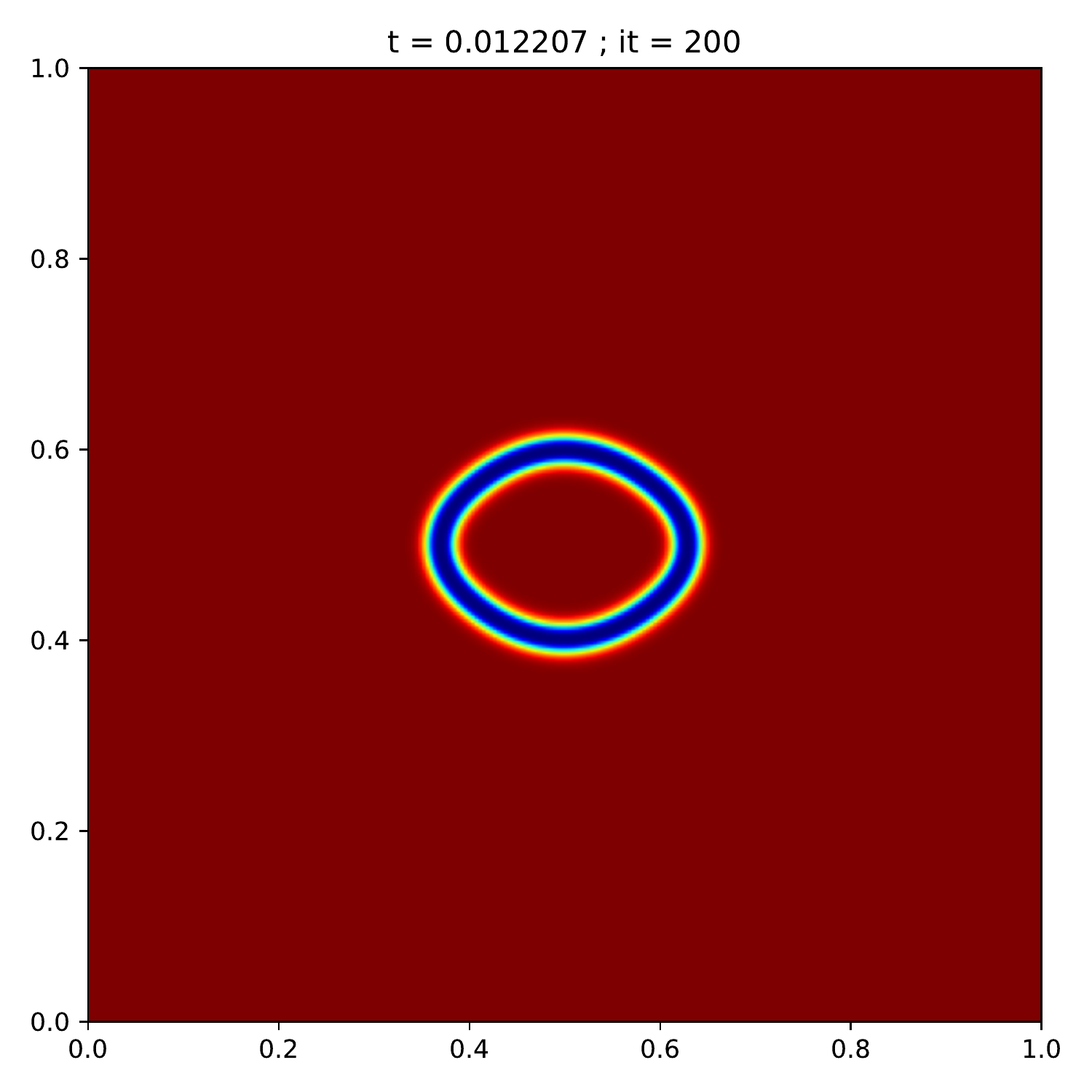}
   \includegraphics[width=0.2\textwidth]{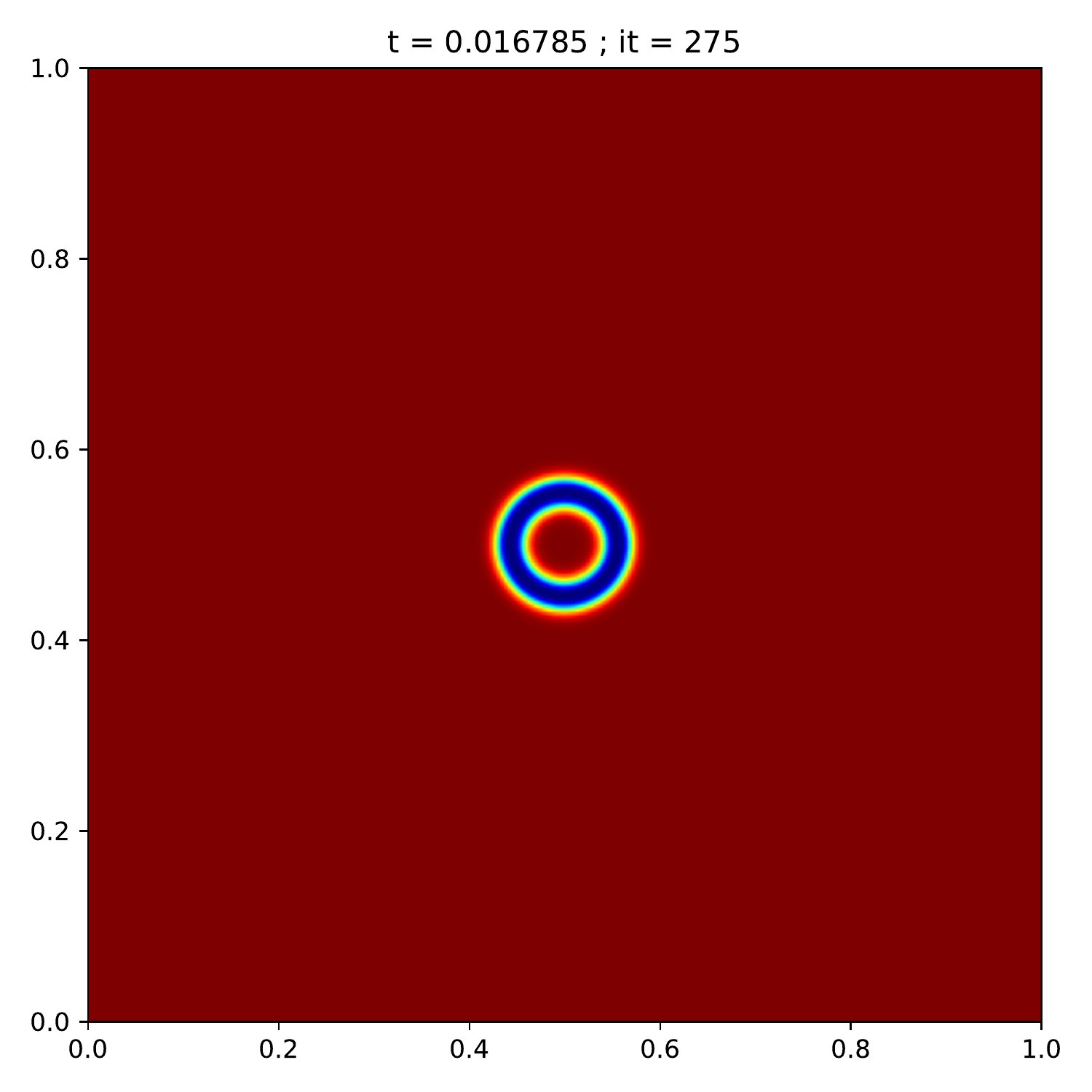}\\
     \includegraphics[width=0.2\textwidth]{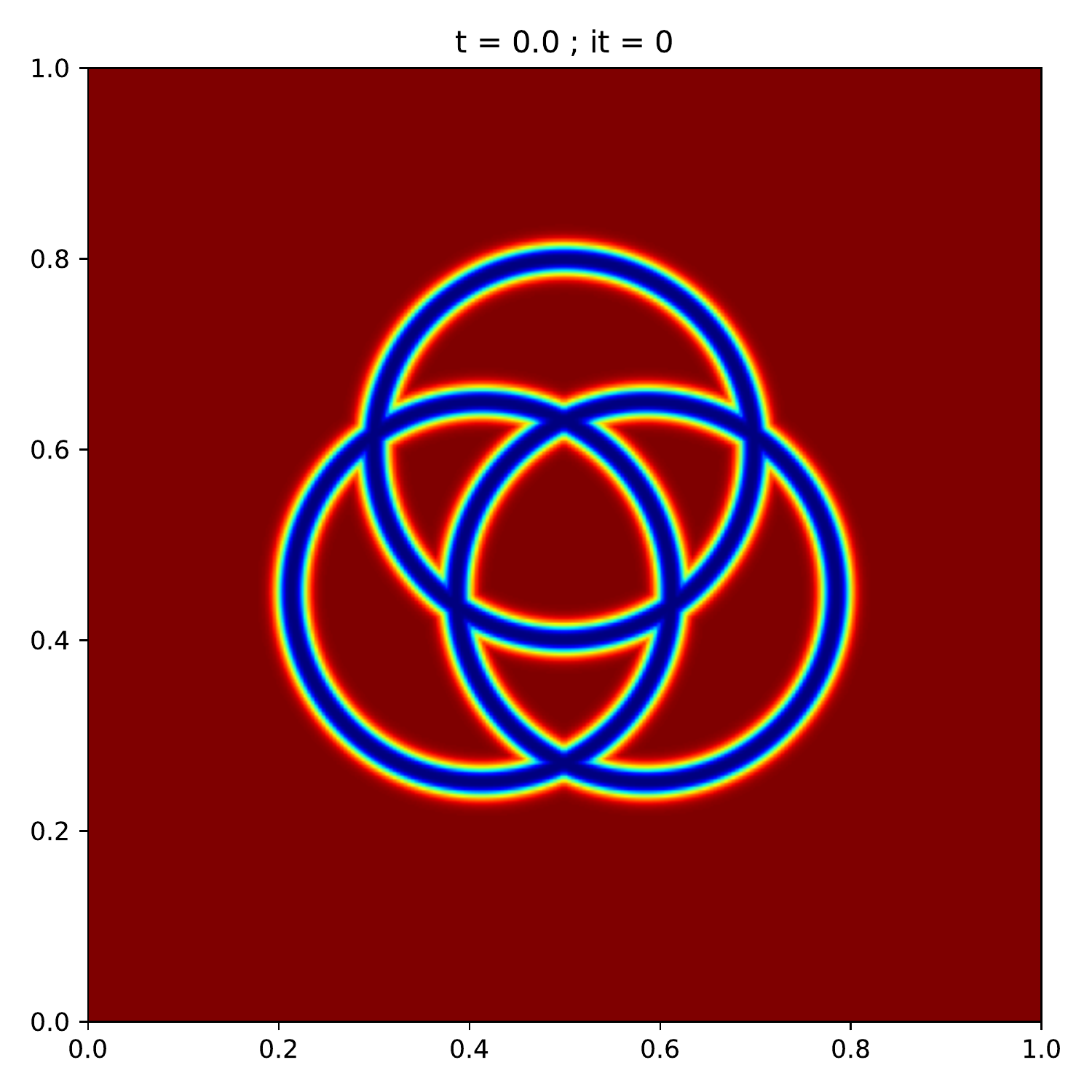}
   \includegraphics[width=0.2\textwidth]{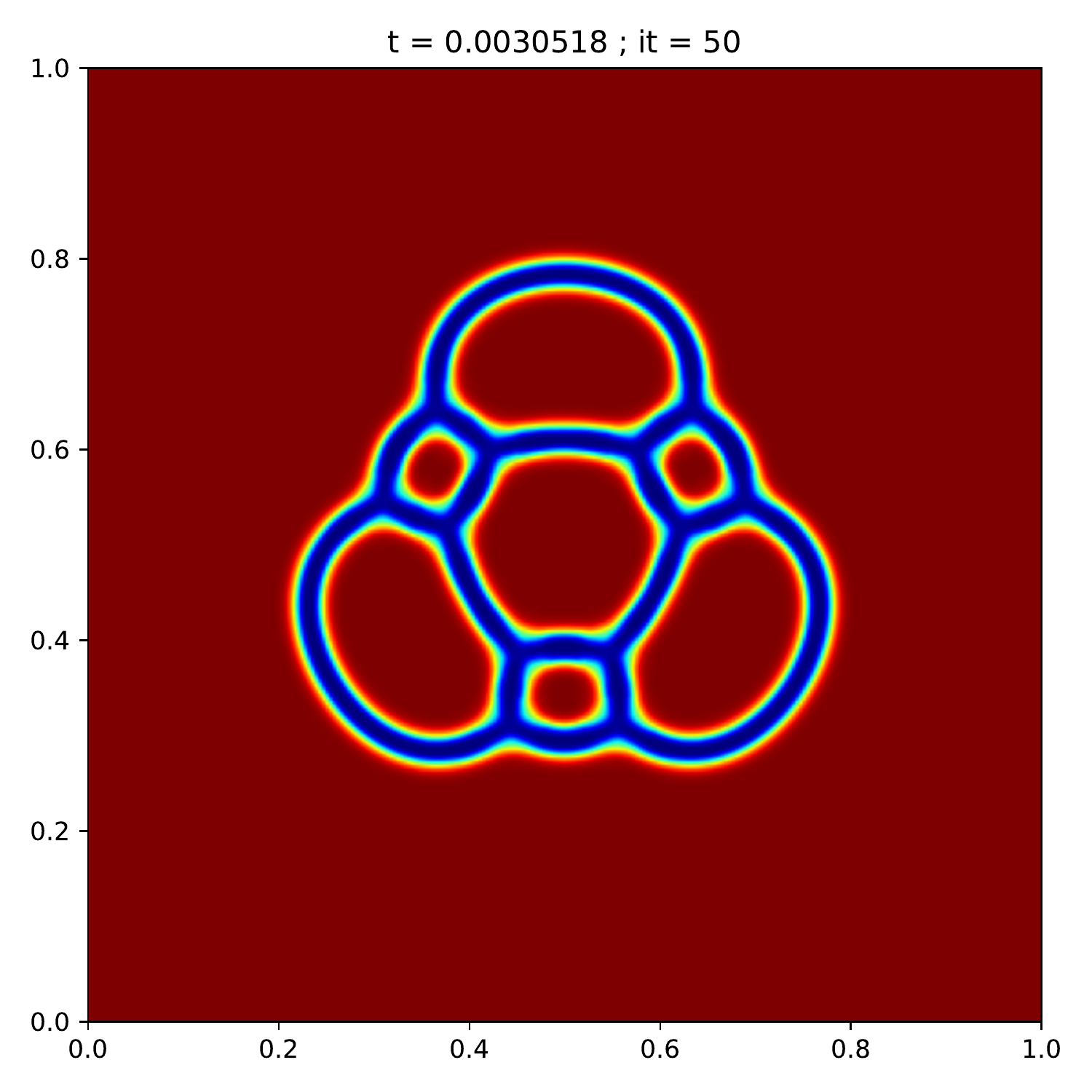}
   \includegraphics[width=0.2\textwidth]{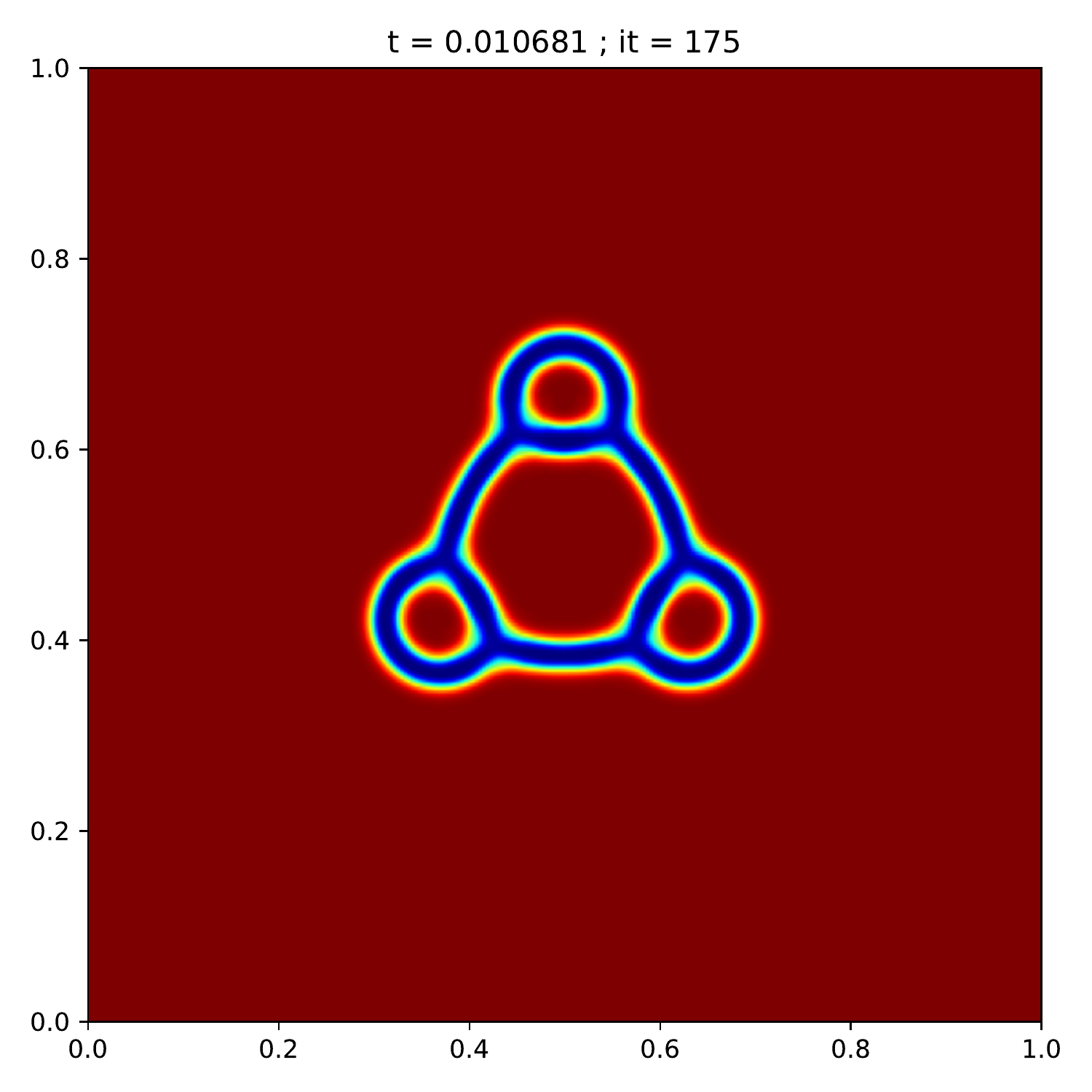}
   \includegraphics[width=0.2\textwidth]{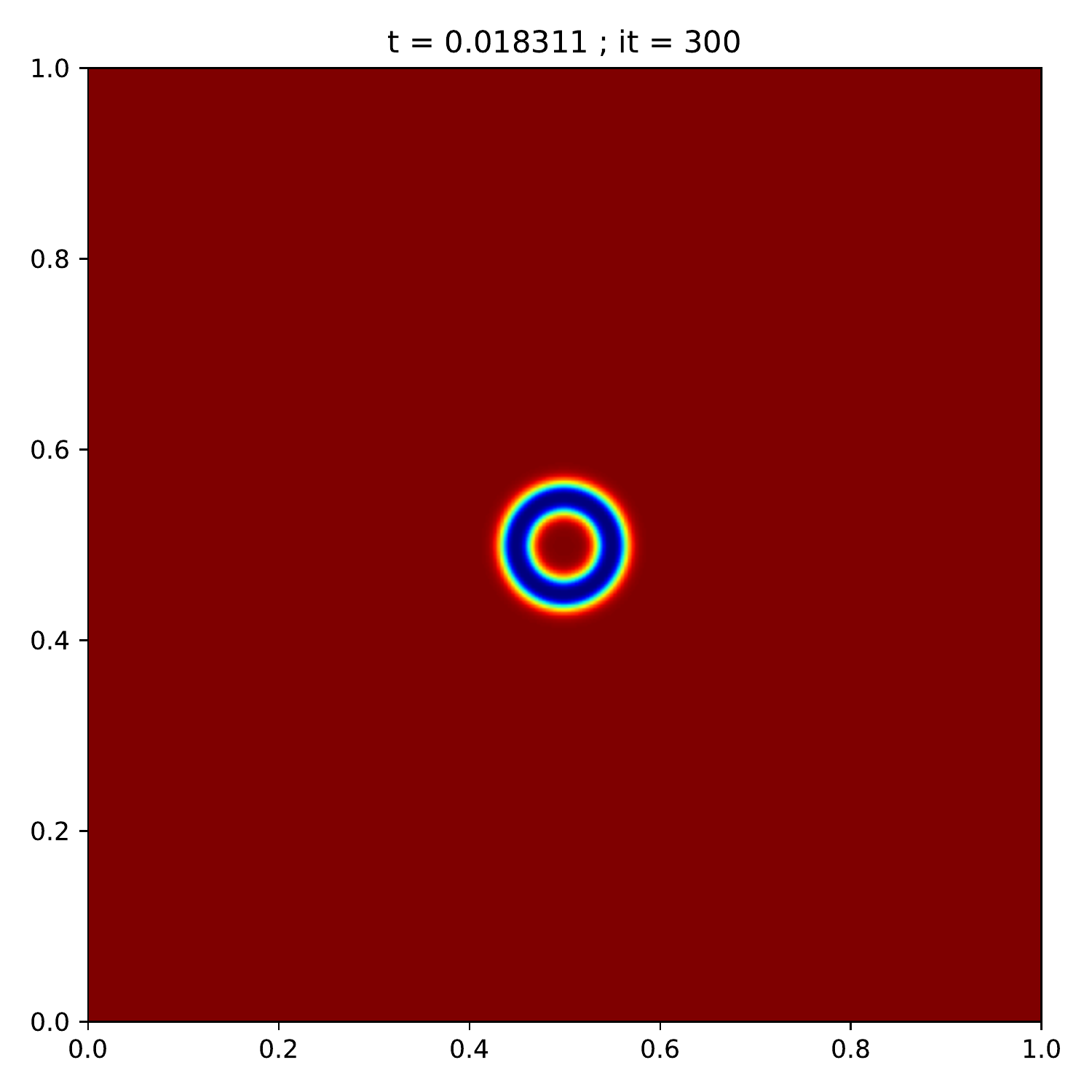}
    \caption{Approximation of the mean curvature flow of a non orientable initial set. Each line
     corresponds to a different choice of the initial set $\Gamma(0)$. We display the solution $u^{n}$ at different times $t_n$ along the iterations.}
    \label{fig:triplejunction}
\end{figure}

We also plot on figure \ref{fig:triplejunction} the numerical approximation of the mean curvature flow  obtained in the case
of non-orientable initial sets. 
In each of these experiments, we observe
the evolution of points with triple junction which seem to satisfy the Herring condition \cite{bretin2018multiphase}.   
This result is all the more surprising since our training base contained only circle evolutions and no interface with the triple point.
The presence of stable triple points is therefore very good news and shows the potential of this approach for the general case of mean curvature motion. 

\begin{figure}[htbp] 
    \centering
    \includegraphics[width=0.2\textwidth]{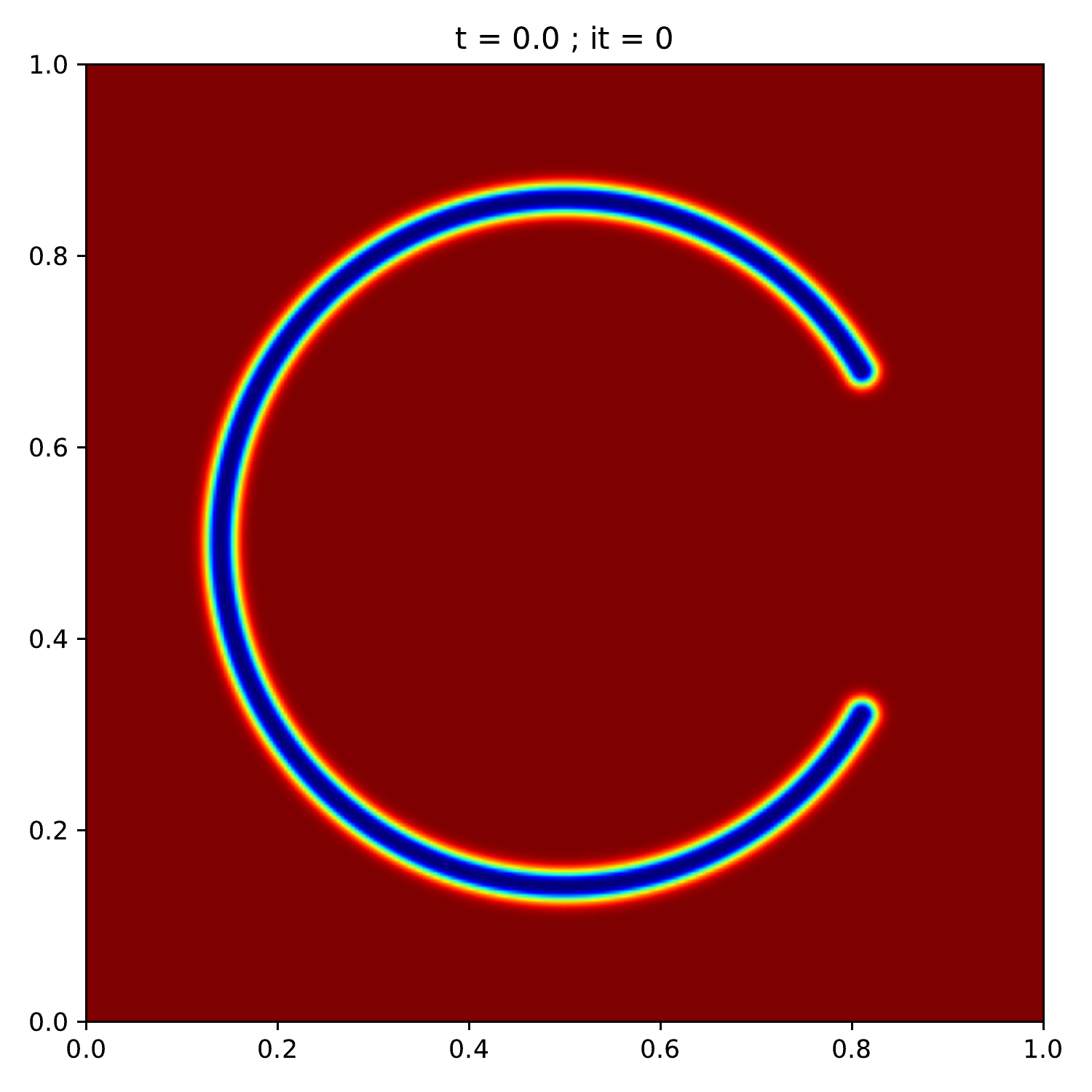}
    \includegraphics[width=0.2\textwidth]{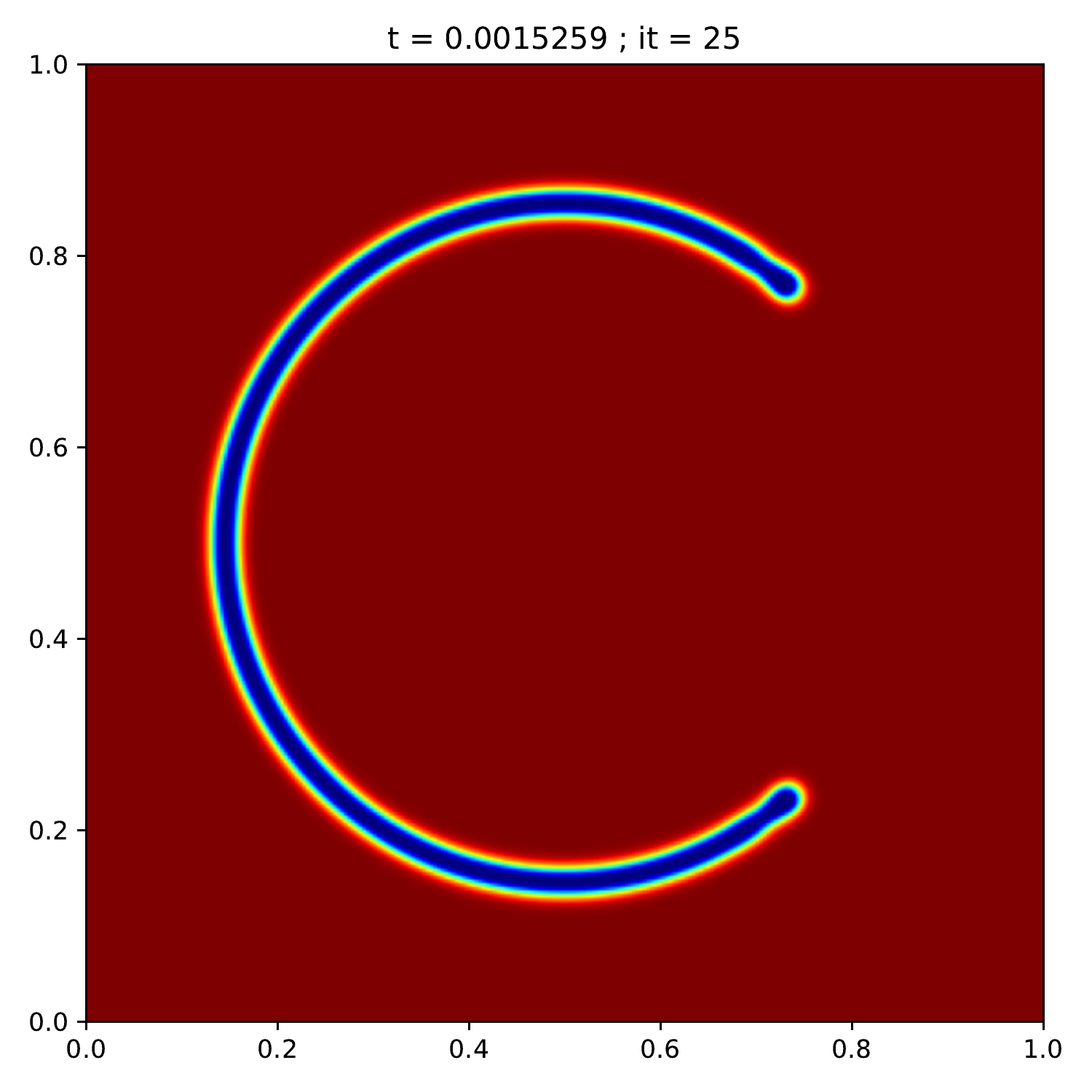}
    \includegraphics[width=0.2\textwidth]{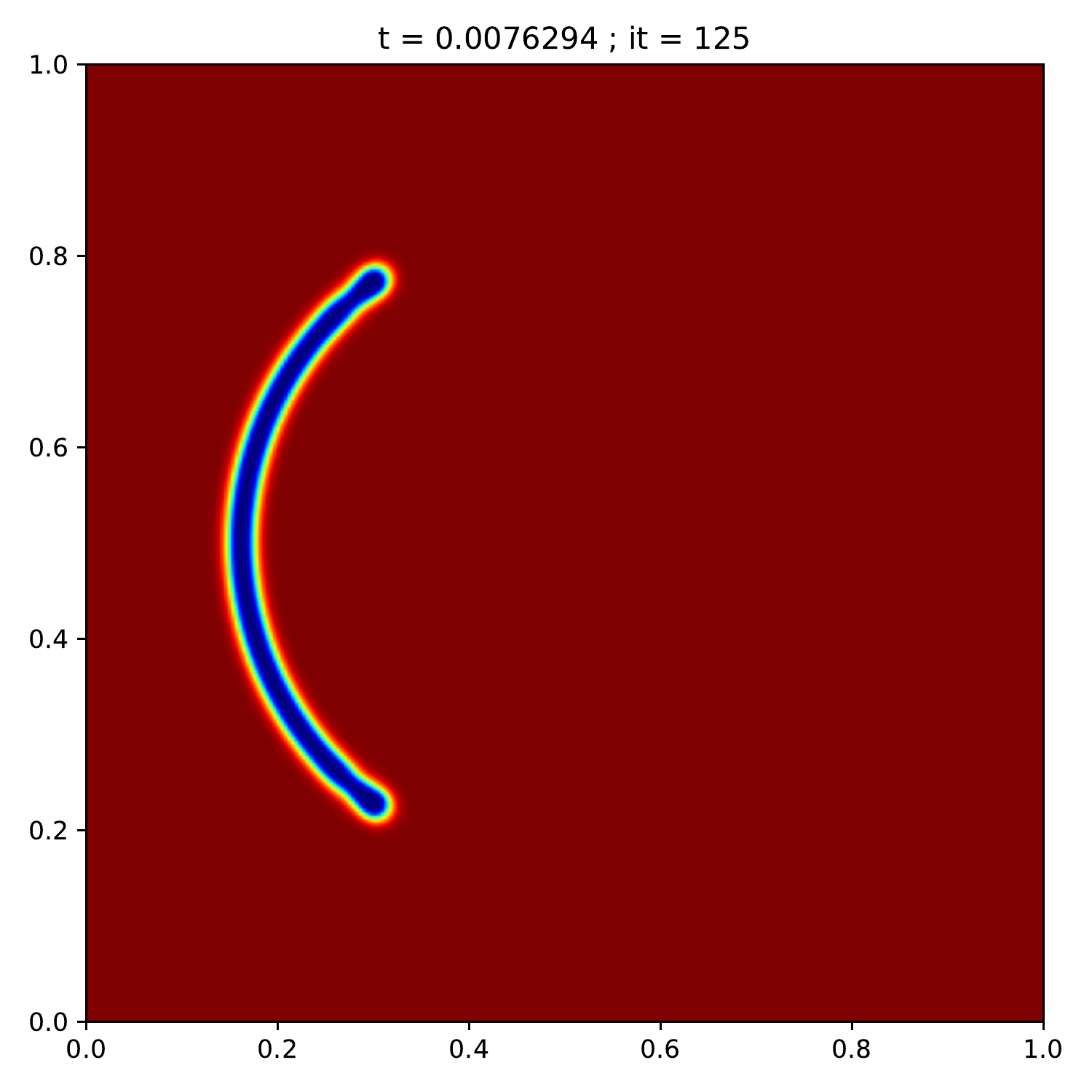}
    \includegraphics[width=0.2\textwidth]{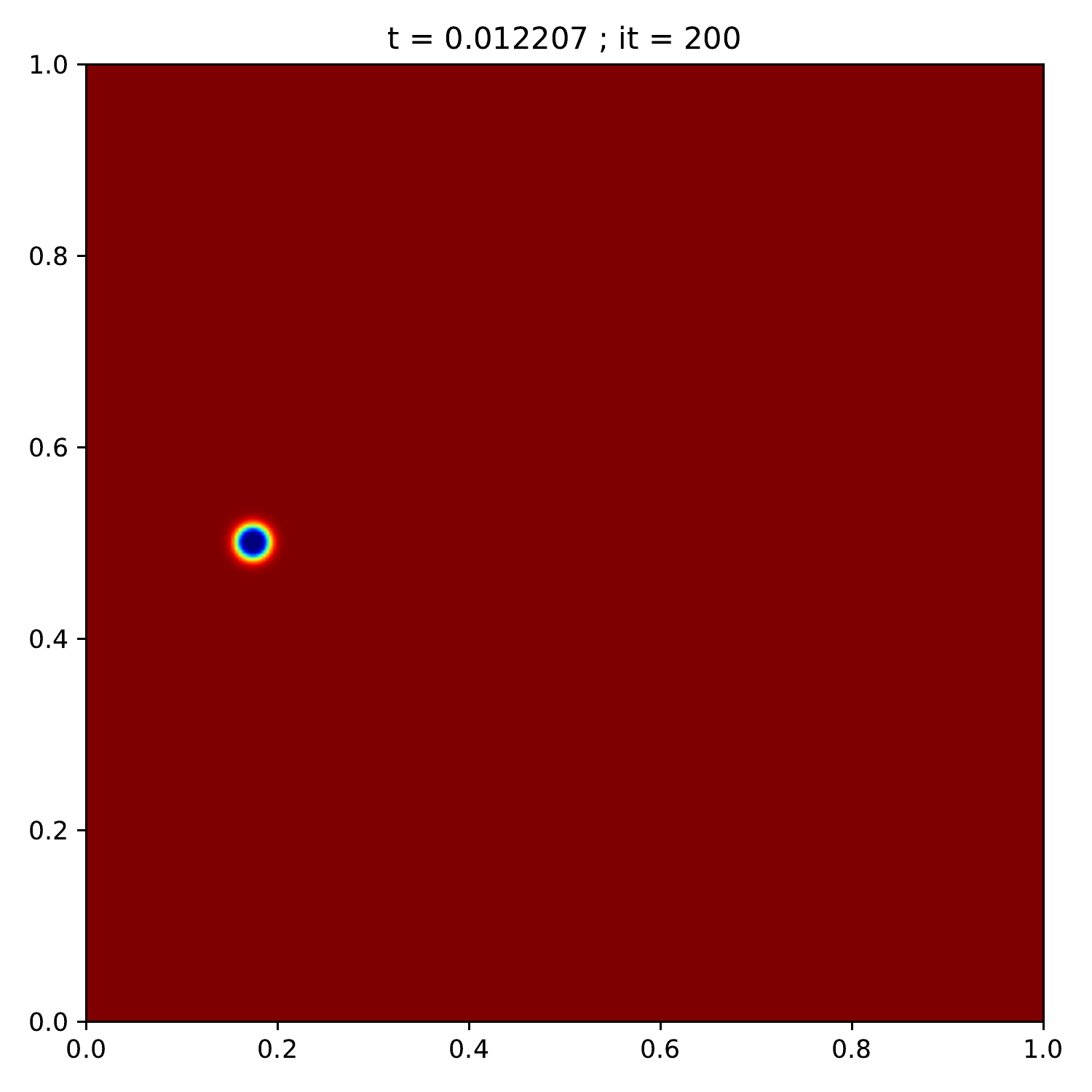} \\
    \caption{Approximation of the motion by mean curvature of an initial open curve.}
    \label{fig:Phase_field_non_closed}
\end{figure}

We give a last example in figure \ref{fig:Phase_field_non_closed} of an
approximation of mean curvature flow in the case of an initial non-closed interface.  
Here the example is quite pathological since it is very difficult to make sense of the motion by mean curvature,  
and we expect the interface to disappear in infinitesimal time. 
Surprisingly, the evolution seems relatively stable with a speed of the endpoints of the order of $1/\varepsilon$ (which is consistent with the $O(1)$ speed at regular points of the circle).

\section{Applications: multiphase mean curvature flows, Steiner trees, and minimal surfaces}
\label{sec:application}

We highlight in this section how the previous schemes derived from our trained networks $S^{NN}_{\theta,\alpha},~\alpha \in \{1,2\}$,
are sufficiently stable to be coupled with additional constraints such as volume conservation ($\int_Q u^{n} dx$ constant) or
inclusion constraints (reformulated as an inequality $ u^{n} \geq u_{in}(x) $, see~\cite{MR3738845}), and sufficiently stable to be extended to the multiphase case. 

\textcolor{black}{
Instead of retraining the networks for these applications, we simply use the same networks trained earlier on circles and spheres, 
and we couple their
predictions with an additional constraint in a post-processing step to approximate the solution. Doing so, we demonstrate the adaptability of our neural network approach.
}

{\bf Partition problem in dimension $2$.}
The first application concerns the approximation of multiphase mean curvature flow with or without volume conservation. 
We recall that the evolution of a partition with $N$ phases ${\bf \Omega(t)} = (\Omega_1(t),\Omega_2(t), \cdots \Omega_N(t))$ can be obtained
as the $L^2$-gradient flow of a multiphase perimeter 
$$ \textcolor{black}{ P({\bf \Omega}) = \frac{1}{2} \sum_{i=1}^{N} P(\Omega_i) =  \frac{1}{2}  \sum_{i=1}^{N}\int_{\partial \Omega_i} d\H^{1}},$$
with the normal velocity $V_{i,j}$ of each interface $\Gamma_{i,j} = \partial \Omega_i \cap \partial \Omega_j$ satisfying
$$ V_{i,j} = H_{i,j},$$
where $H_{i,j}$ represents the mean curvature at the interface $\Gamma_{i,j}$.\\

{\bf The Steiner problem in dimension $2$.}
The second application is the approximation of Steiner trees in dimension $2$, i.e. solutions of the Steiner problem. 
Recall that the Steiner problem consists, given a collection of points $a_1, \ldots, a_L \in \Omega$,  in finding a compact connected set $K\subset \Omega$ containing all the $a_i$'s and having minimal length. 
In other words, it is equivalent to finding the optimal solution to the following problem
	$$
	\min\left\{\H^1(K),~ K\subset \Omega, K \text{ connected, } a_i\in K, ~\forall i=1,\cdots, L \right\},
	$$
	where $\H^1(K)$ corresponds to the one-dimensional Hausdorff measure of $K$. 
As the optimal set $K$ needs not be orientable, our idea is to combine the network $S^{NN}_{\theta,2}$ trained on the non-oriented phase field database
with the inclusion constraints associated to all points $a_i$'s. \\

{\bf The Plateau problem in dimension $3$.} The last application focuses on the Plateau problem in dimension $3$ and co-dimension $1$.
Recall that it consists in finding, for 
a given closed boundary curve $\Gamma$, a compact set $E$ in $\Omega$ with minimal area and whose boundary coincides with 
$\Gamma$~\cite{david-hal-00718979}. In other words, it amounts to solving the minimization problem 
\begin{equation}
 \min \{  \H^2(E) ; E \subset \Omega, \text{connected and such as } \partial E = \Gamma \},\end{equation}
where $\H^2(E)$ stands for the two-dimensional Hausdorff measure of  $E$. 
Analogously to the method proposed for the Steiner problem, the approximation method we propose consists in 
\begin{enumerate}
 \item training a $S^{NN}_{\theta,2}$ network 
on databases built from $3D$ spheres evolving under mean curvature flow with a non-oriented phase field representation,
 \item coupling the trained network $S^{NN}_{\theta,2}$ with the inclusion constraint to force the boundary constraint 
 $\partial E = \Gamma$.
\end{enumerate}

\subsection{Evolution of a partition in dimension $2$}
As explained previously, the motivation in this first application is to approximate
multiphase mean curvature flows. Recall that the phase field approach 
consists in general~\cite{Garcke_amulti,Garcke199887,garckes-haas08,oudet2011,bretin2017new} in introducing 
a multiphase field function  ${\bf u} = (u_1,u_2,\cdots, u_N)$, solution of an Allen-Cahn system  obtained as 
the  $L^2$-gradient flow of the multiphase Cahn-Hilliard energy $P_{\varepsilon}$ defined by 
$$ P_{\varepsilon}({\bf u}) = \begin{cases}
                               \textcolor{black}{ \frac{1}{2} \sum_{k=1}^{N} \int_Q \frac{\varepsilon}{2} |\nabla u_k|^2 + \frac{1}{\varepsilon} W(u_k)dx } &\text{ if }  \sum_{k=1}^{N} u_k = 1 ,\\
                               + \infty &\text{ otherwise.}
                              \end{cases}$$

More precisely, the Allen-Cahn system  \cite{bretin_largephases,bretin2018multiphase}  reads as 
$$ \partial_t u_k = \Delta u_k - \frac{1}{\varepsilon^2} W'(u_k) + \lambda \sqrt{2 W(u_k)}, \quad k=1,\cdots, N$$
where  the Lagrange multiplier $\lambda$ is associated with the partition constraint $\sum_{k=1}^L u_k = 1$.
As explained in  \cite{bretin_largephases,bretin2018multiphase}, this PDE  can be for instance computed in two steps:
\begin{enumerate}
 \item Solve the decoupled Allen-Cahn system 
 $$u_k^{n+1/2} =  \S^{\text{AC}}_{\delta_t,\varepsilon, 1}[u_k^{n}],  \quad k=1,\cdots, N.$$
 \item Project onto the partition constraint $\sum_k u_k^{n+1} = 1$:
 $$  u_k^{n+1} = u_k^{n+1/2} + \lambda^{n+1} \sqrt{2 W(u_k^{n+1/2})},  \quad k=1,\cdots, N.$$
 Notice that the Lagrange multiplier $\lambda^{n+1}$ is computed  following the expression
 $$\lambda^{n+1} = \frac{1 -  \sum_{k=1}^N  u_k^{n+1/2}}{\sum_{k=1}^N  \sqrt{2 W(u_k^{n+1/2})}},$$   in order to satisfy exactly 
 the partition constraint.
\end{enumerate}

\begin{figure}[htbp]
    \centering
    \includegraphics[width=0.22\textwidth]{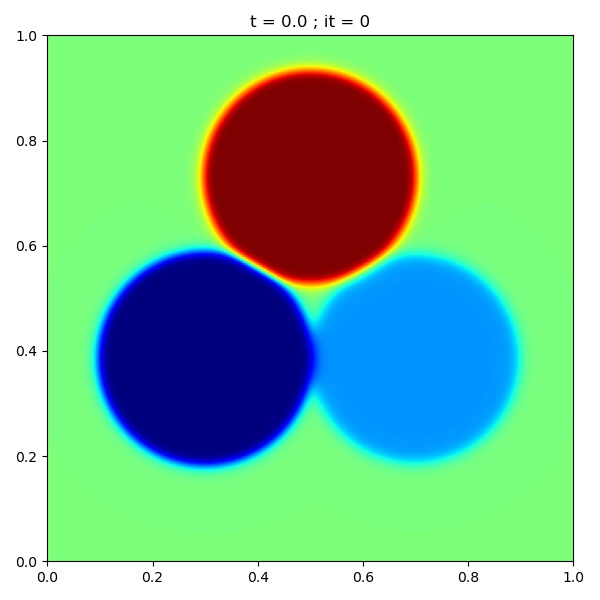}
 \includegraphics[width=0.22\textwidth]{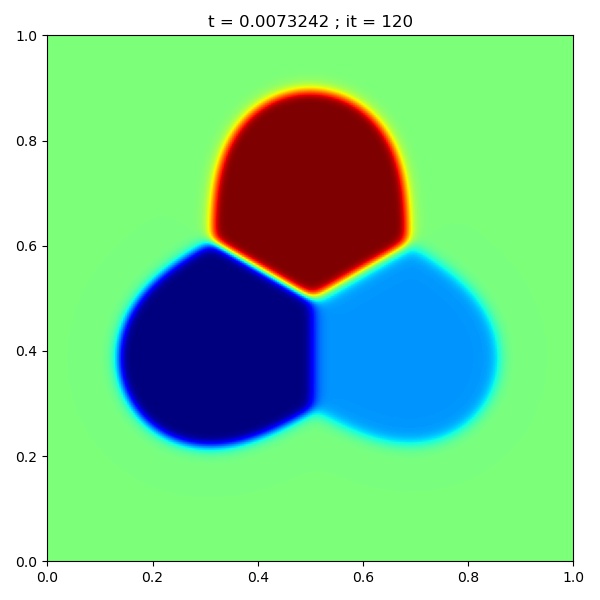}
 \includegraphics[width=0.22\textwidth]{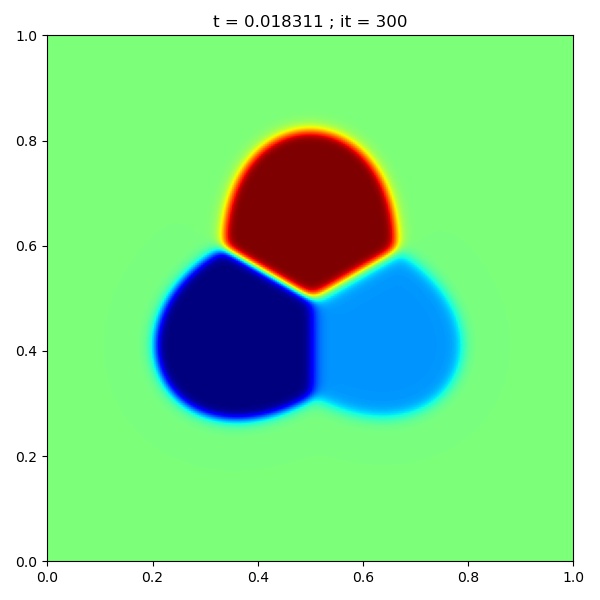}
 \includegraphics[width=0.22\textwidth]{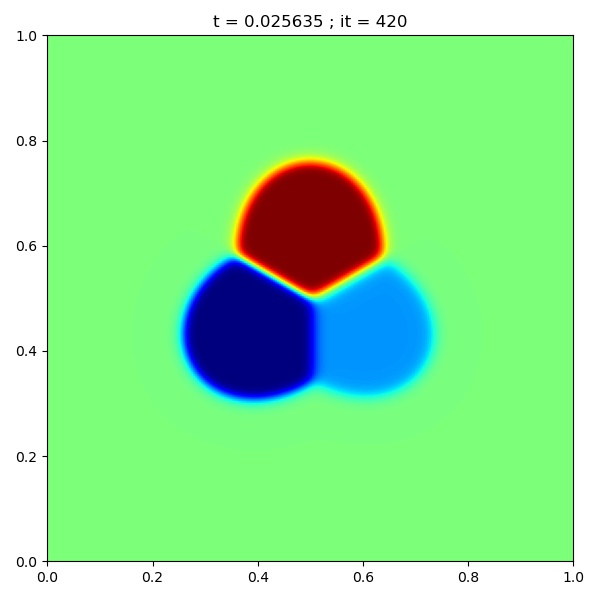} \\
     \includegraphics[width=0.22\textwidth]{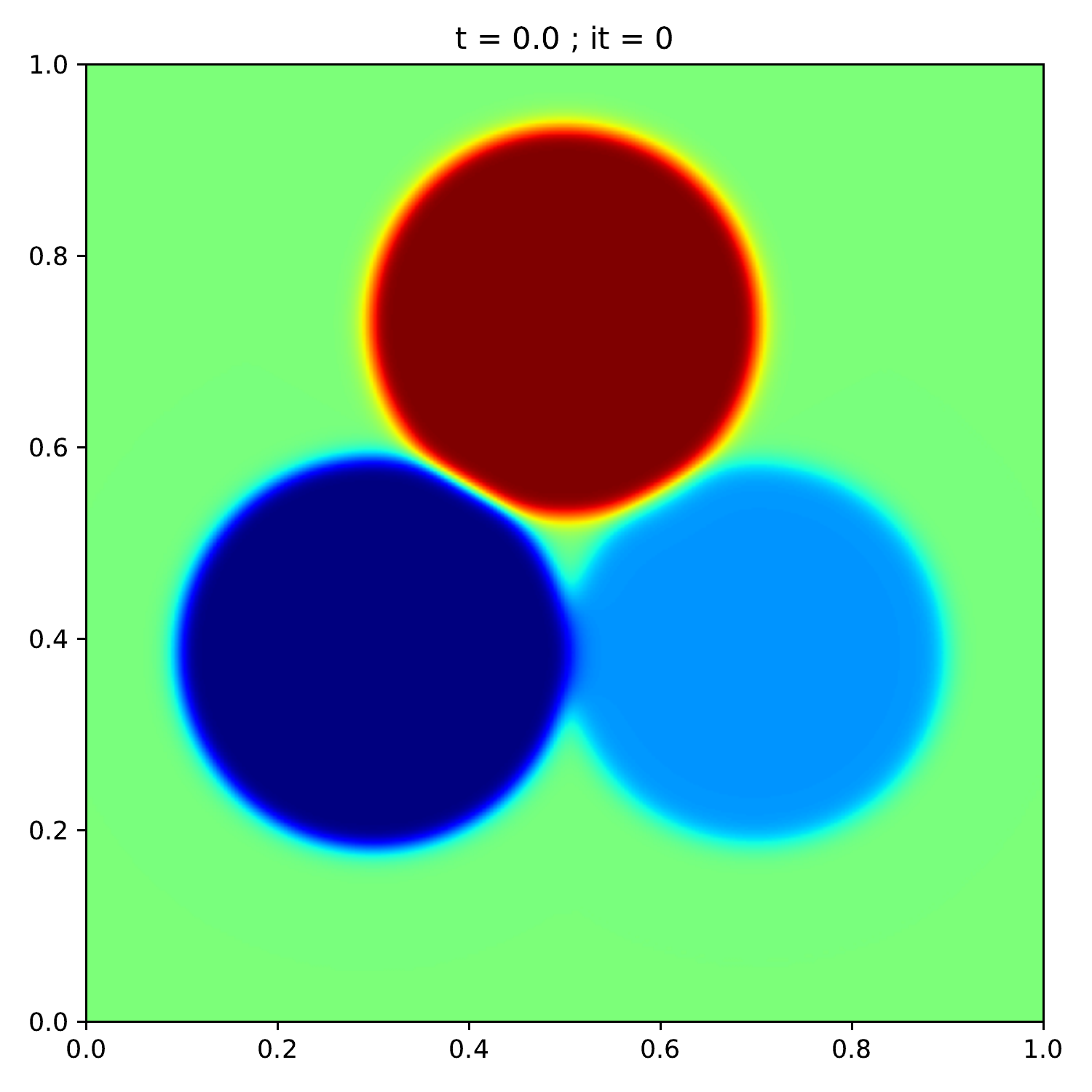}
          \includegraphics[width=0.22\textwidth]{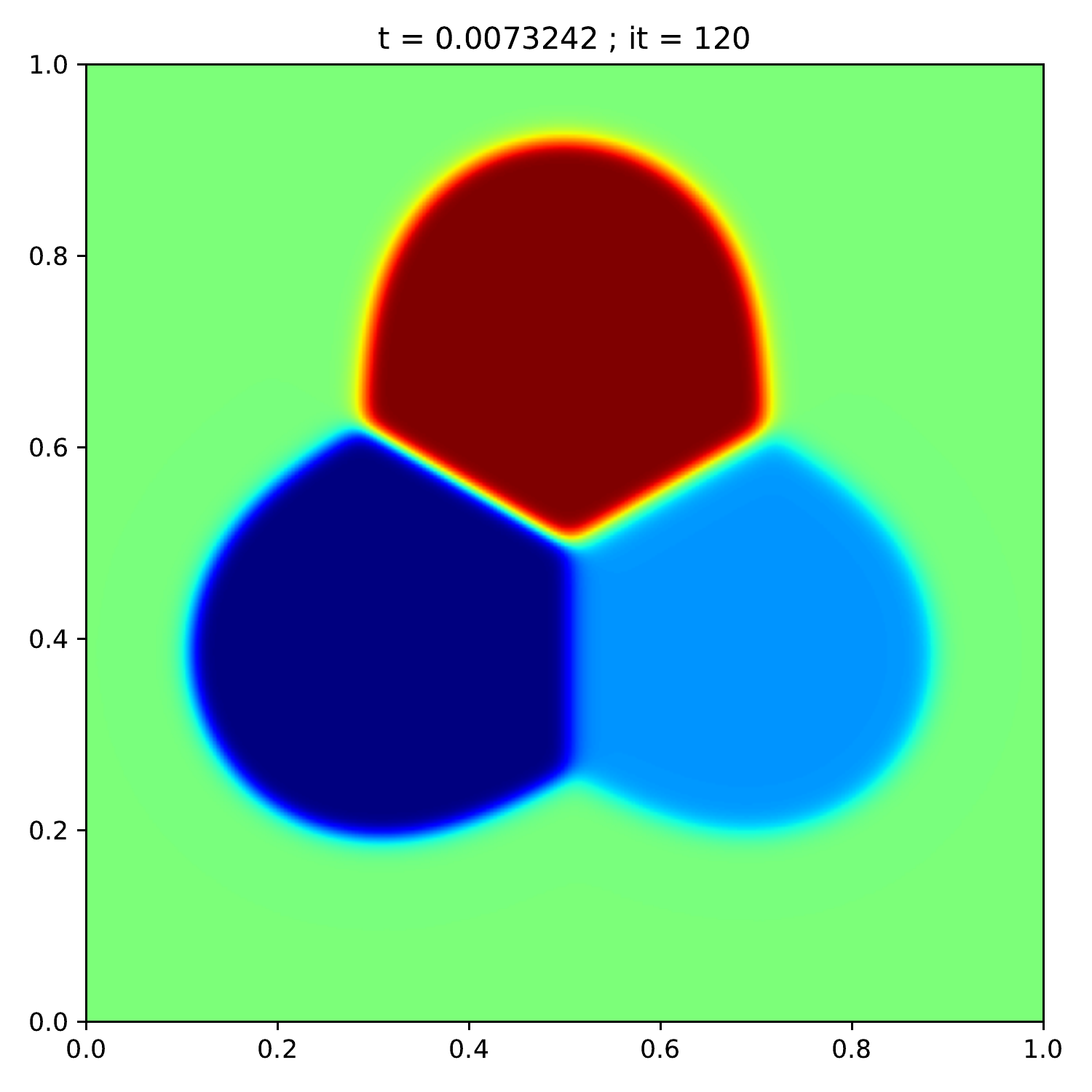}
           \includegraphics[width=0.22\textwidth]{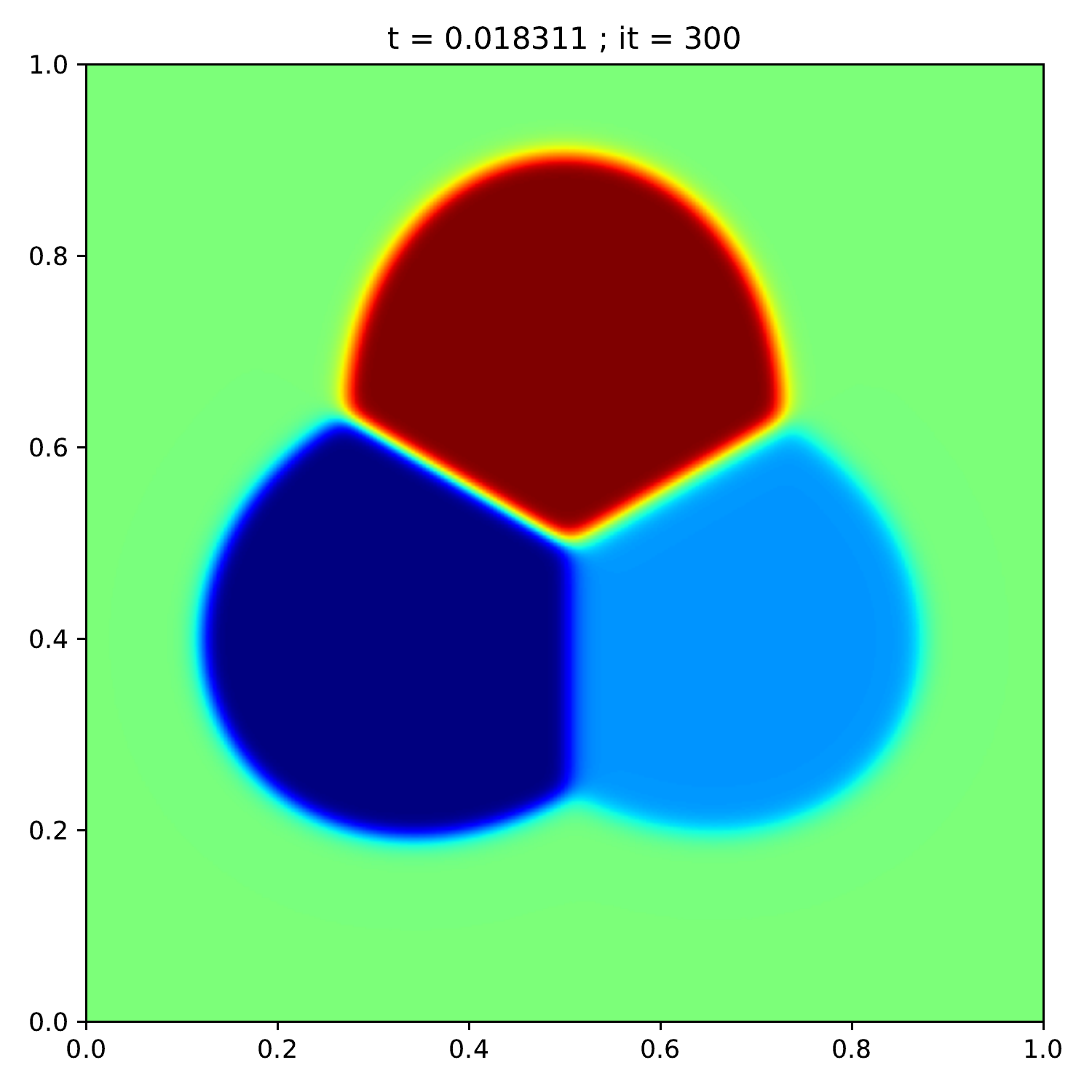}
              \includegraphics[width=0.22\textwidth]{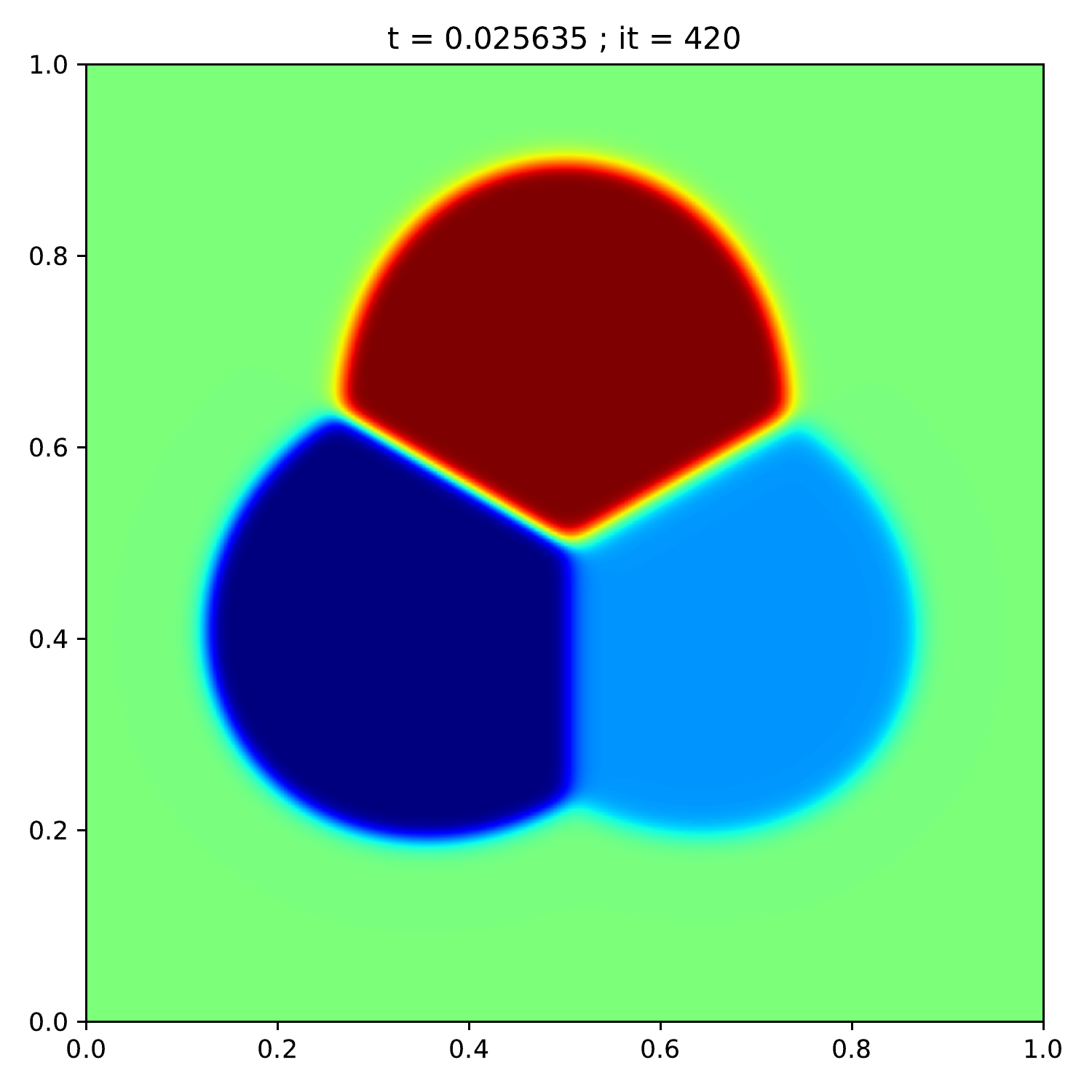}
    \caption{Mean curvature flows of three interconnecting circles without (first row) and with (second row) additional volume constraint.}
    \label{fig:mutliBubblethree}
\end{figure}

This discretization scheme associated with the Allen-Cahn system suggests to simply replace the Lie splitting Allen-Cahn operator
$\S^{\text{AC}}_{\delta_t,\varepsilon, 1}$ by our network $\S^{NN}_{\theta,1}$ trained on oriented mean curvature
flow in section \eqref{subsect:Oriented_MCF}.  
A numerical experiment obtained with such a strategy is presented on the first row of figure~\ref{fig:mutliBubblethree}.
More precisely, we consider here an evolution of a multiphase system with four phases and we observe
at least qualitatively that the flow seems to be correct with a triple junction forming angles of $2\pi/3$.
The second row of figure~\ref{fig:mutliBubblethree} gives a similar numerical experiment with additional constraint 
on the volume of each phase. Here, following the approach developed in \cite{MR2439843,PhysRevE.78.011604,BrasselBretin,bretin_largephases}, the idea is
to consider the Allen-Cahn system with volume conservation
$$ \partial_t u_k = \Delta u_k - \frac{1}{\varepsilon^2} W'(u_k) + \lambda \sqrt{2 W(u_k)} +  \mu_k \sqrt{2 W(u_k)},$$
where the new Lagrange multiplier $\mu_k$ is defined to satisfy  $\mathrm{Vol}_k = \int_Q  u_k =  \int_Q  u_0$. The previous scheme
can then be  modified by considering now a projection on the partition and the volume constraints as 
 $$  u_k^{n+1} = u_k^{n+1/2} + (\lambda^{n+1} + \mu_k^{n+1}) \sqrt{2 W(u_k^{n+1/2})},$$
 where the Lagrange multipliers $(\lambda^{n+1},\mu_1^{n+1},\mu_2^{n+1}, \cdots \mu_N^{n+1})$ are for instance given by
 $$ \mu_k^{n+1} = \frac{\mathrm{Vol}_k - \int_Q u^{n+1/2}_k}{ \int_Q \sqrt{2 W(u_k^{n+1/2}})}, \text{ and } \lambda^{n+1} = \frac{1 -  \sum_{k=1}^N \left[u_k^{n+1/2} + \mu_k^{n+1} \sqrt{2 W(u_k^{n+1})} \right]}{\sum \sqrt{2 W(u_k^{n+1/2})}},$$
which allows each of the constraints to be satisfied.  As previously, the idea consists simply to replace
 the Lie splitting Allen-Cahn operator $\S^{\text{AC}}_{\delta_t,\varepsilon_1}$ by our network $\S^{NN}_{\theta,1}$. 
 In particular, the numerical experiment plotted on the second row
of figure~\ref{fig:mutliBubblethree} shows a stable multiphase evolution where the volume of each phase
seems to be conserved.  These two numerical examples clearly show that our networks previously trained on 
mean curvature motion flow can be re-exploited in more complex situations.  
  
\subsection{Approximation of Steiner trees in $2D$}~\\
Variational phase field models have been recently introduced 
\cite{MR3337998,MR4011534,BonafiniOrlandiOudet,BonafiniOudet,ChambolleFerrariMerlet2019-1} 
to tackle the Steiner problem. For example, the model introduced in \cite{MR3337998} proposes the Ambrosio-Tortorelli functional  
$$ F_{\varepsilon}(u) =  \int_{\textcolor{black}{\Omega}} \varepsilon  |\nabla u|^2 + \frac{1}{\varepsilon} (1-u)^2 dx,$$
coupled with an additional penalization term $G_{\varepsilon}(u)$ which forces the connectedness of the set $K$ 
and which is defined by
$$
G_{\varepsilon}(u) = \frac{1}{\varepsilon} \sum_{i=1}^{N} {\bf D}(u^2; a_0 , a_i),  \text{ with } {\bf D}(w;a,b):=\inf_{\Gamma: a \leadsto b}\int_{\Gamma} w\,d\mathcal{H}^1 \in [0,+\infty],
$$
where the notation $\Gamma: a \leadsto b$ means that $\Gamma$ is a rectifiable curve in $\Omega$ connecting $a$ and $b$.

Intuitively, the phase field function $u_{\varepsilon}$ is expected to be of the form $u_{\varepsilon} = \varphi(\operatorname{dist}(x,K)/ \varepsilon)$, the Ambrosio-Torterelli term approximates the length of $K$, and the geodesic term $G_{\varepsilon}$
forces the solution $u_\varepsilon$ to vanish on $K$, and $K$ to connect the different points $a_i$.
An analysis of this  model and, more precisely,
 a $\Gamma$-convergence result established in \cite{MR3337998} suggest that the minimization of this phase field model 
 should give an approximation of a Steiner solution, i.e. a Steiner tree.  However, numerical experiments as proposed in  \cite{MR3337998} and in  \cite{MR4011534} show the ability of this coupling to approximate solutions of the Steiner problem in dimensions $2$ and $3$
but show also all the difficulties to minimize efficiently  the geodesic term $G_{\varepsilon}$ to preserve the connectedness of the set $K$.
The conclusions are quite similar for the approaches developed in \cite{BonafiniOrlandiOudet,BonafiniOudet,ChambolleFerrariMerlet2019-1} 
where the idea is rather to use the measure-theoretic notion of current and the connectedness of the set $K$ 
is ensured by adding a divergence constraint of the form  $ div(\tau) = \sum \alpha_i \delta_{a_i}$. \\

We propose a completely different approach by considering a non oriented mean curvature flow 
$t\mapsto \Gamma(t)$ of an initial connected set $\Gamma(0)$ containing all the points $a_i$ coupled with the additional inclusion
constraint $\{a_i,~i=1,\cdots,N\} \subset \Gamma(t)$. Indeed, we expect that the stationary state of such an evolution
is at least a local minimum
of the Steiner problem. From a phase field point of view, the strategy is to consider
the non-oriented phase field representation $q'(\operatorname{dist}(x,\Gamma(t)/\varepsilon))$, to use the non oriented trained network 
$\S^{NN}_{\theta,2}$ to let the interface evolve by mean curvature, and to incorporate  the inclusion constraint
by using an additional inequality constraint 
$$ u \leq u_{\text{in}}(x) = \sum_{i=1}^{N}  q'(\operatorname{dist}(a_i,x)/\varepsilon),$$
in the spirit of \cite{MR2996338,MR3738845}.
Finally, the scheme reads as follows:
\begin{enumerate}\label{scheme:Non-oriented+inclusion}
 \item Approximation of a non-oriented mean curvature flow step  
 $$u^{n+1/2} =  \S^{\text{NN}}_{\theta,2}[u^{n}].$$
 \item Projection on the inclusion constraint $ u \leq u_{\text{in}}(x)$:
 $$  {u^{n+1}} = \text{min}(u_{\text{in}},u^{n+1/2}).$$
\end{enumerate}

\begin{figure}[htbp]
    \centering
    
   \includegraphics[width=0.22\textwidth]{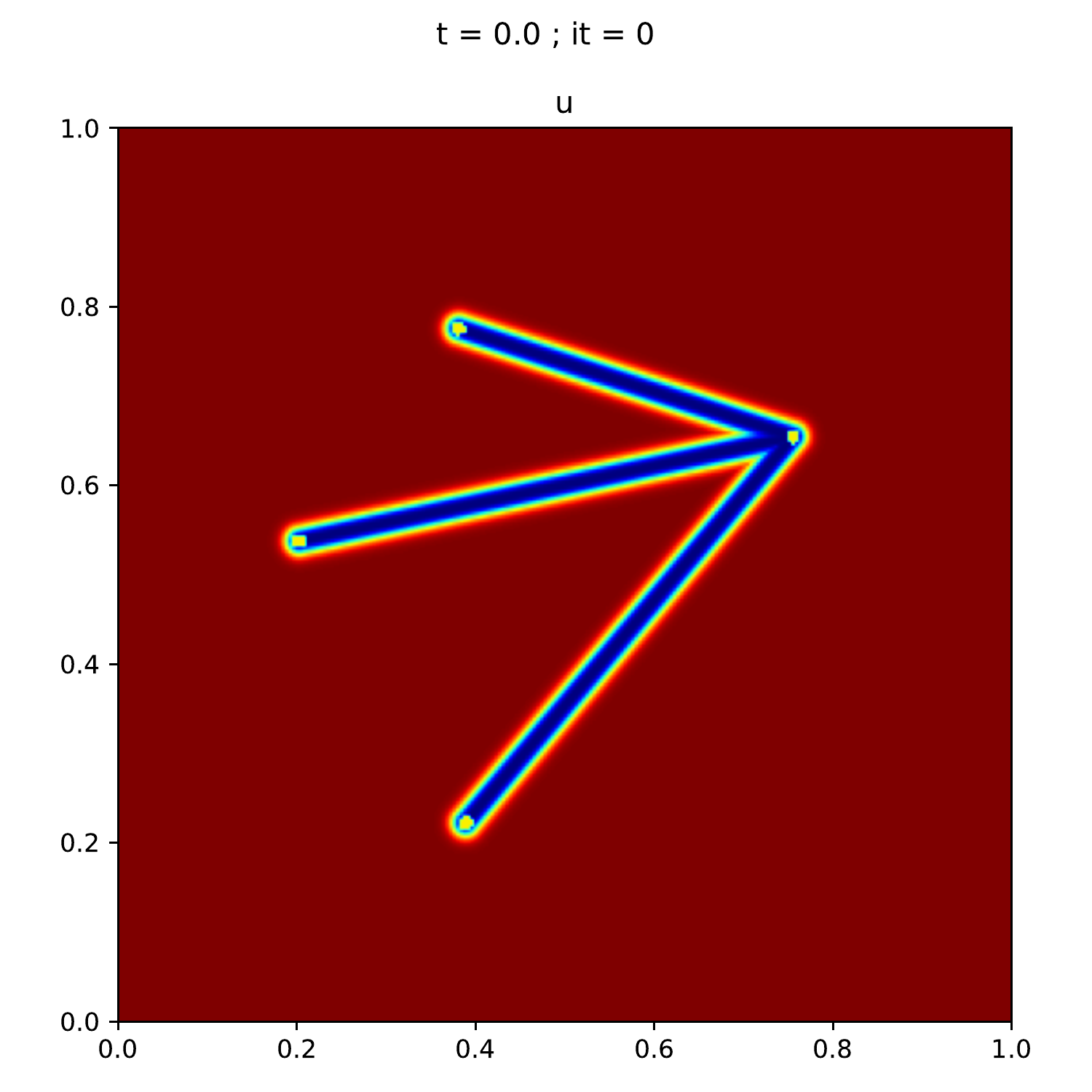}
   \includegraphics[width=0.22\textwidth]{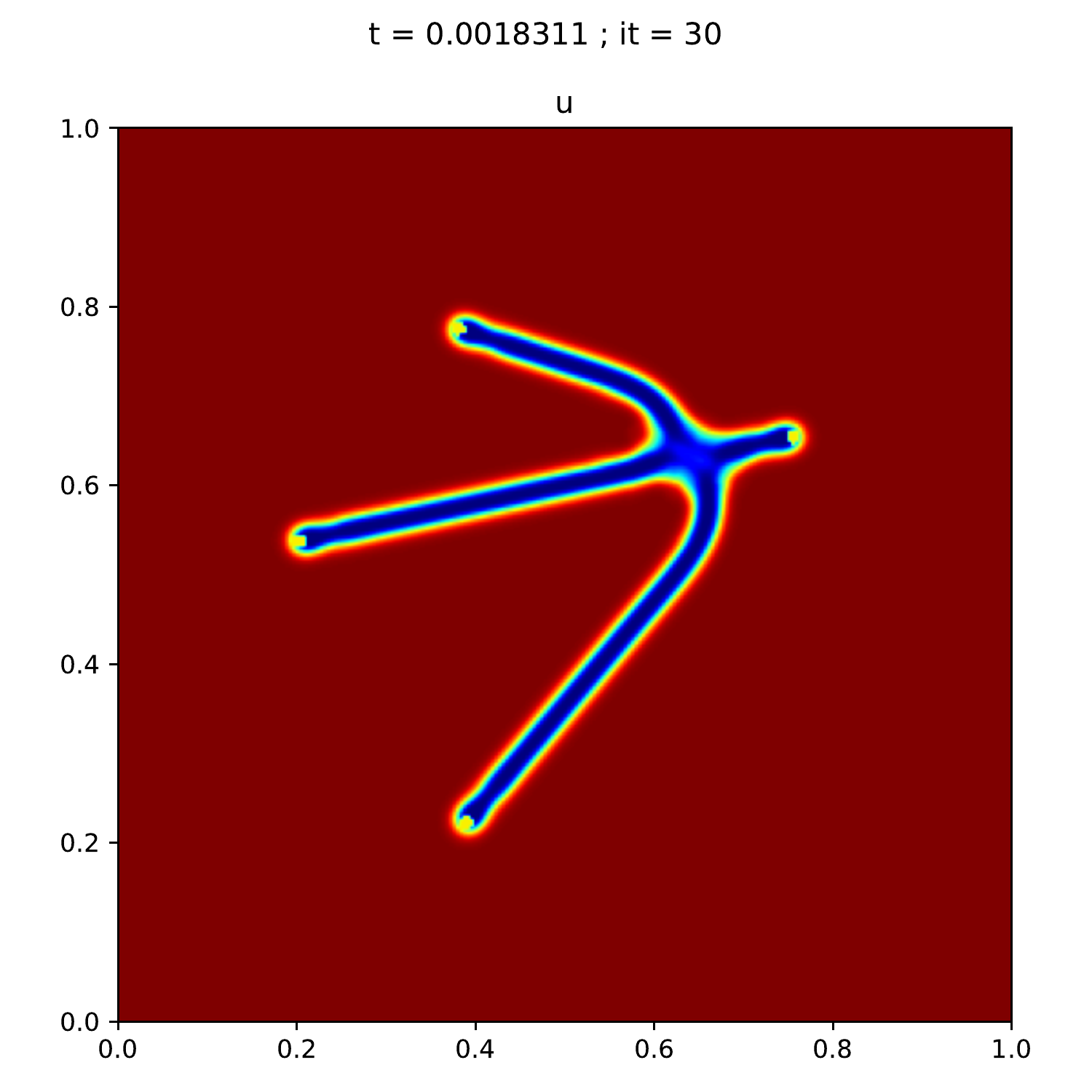}
   \includegraphics[width=0.22\textwidth]{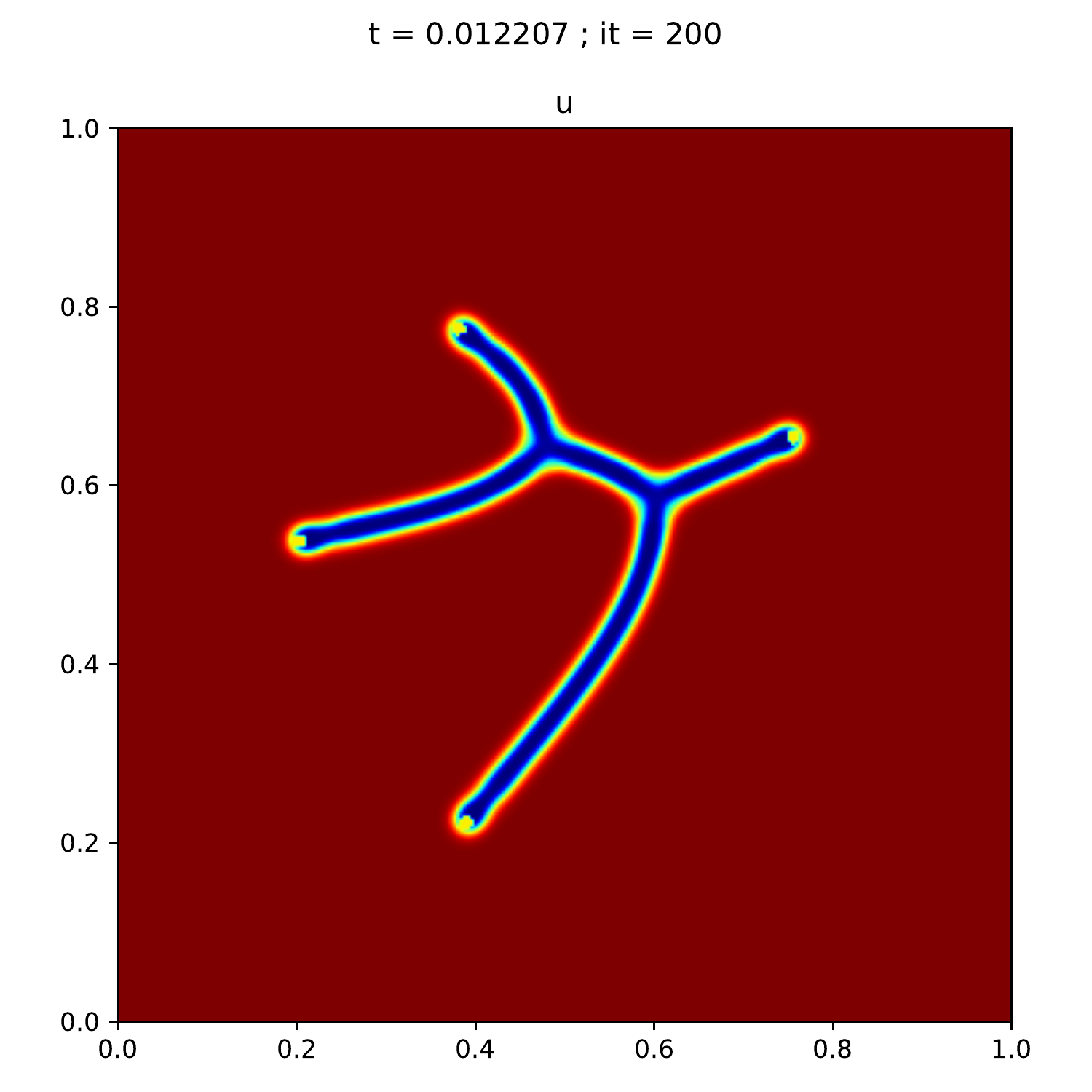}
   \includegraphics[width=0.22\textwidth]{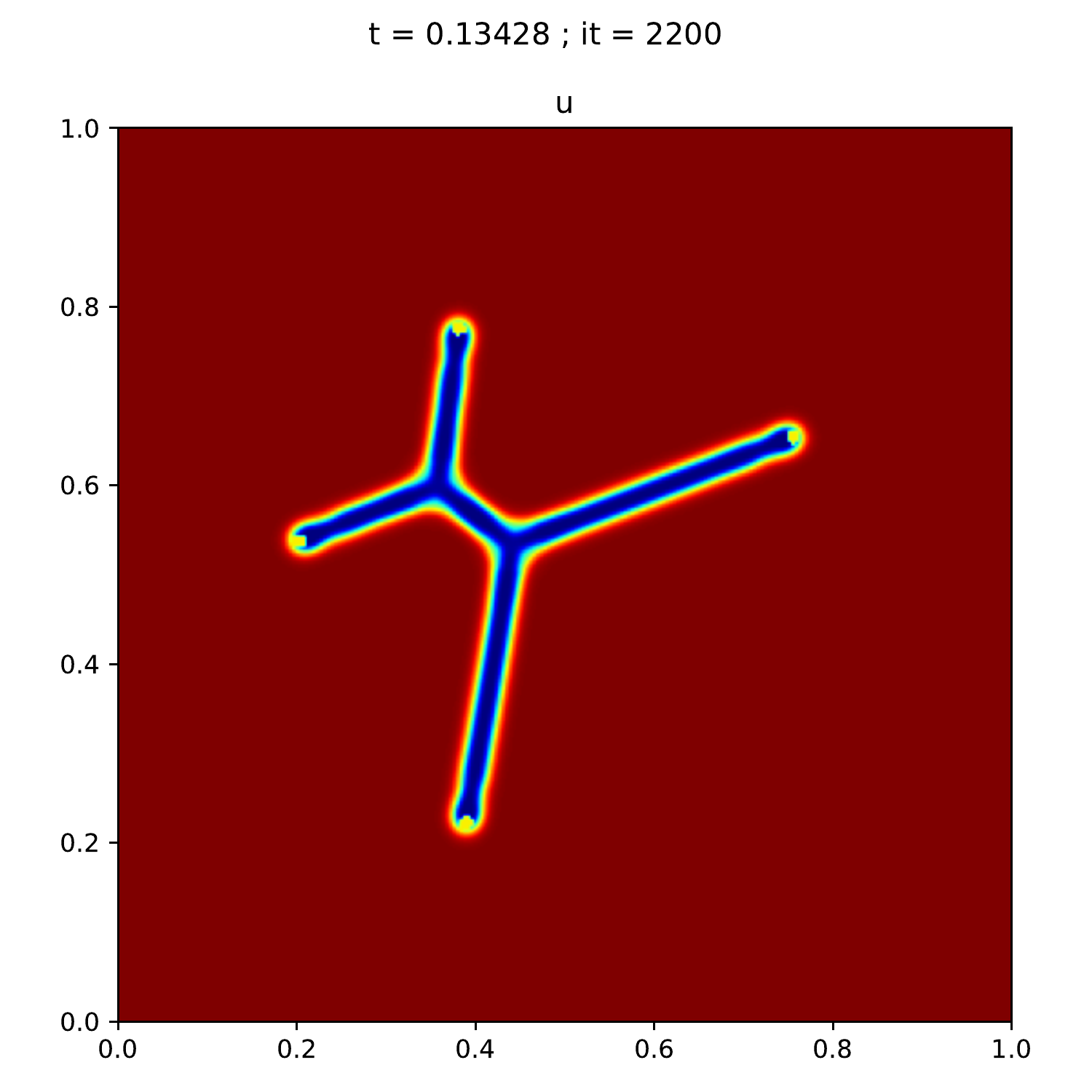}

   \includegraphics[width=0.22\textwidth]{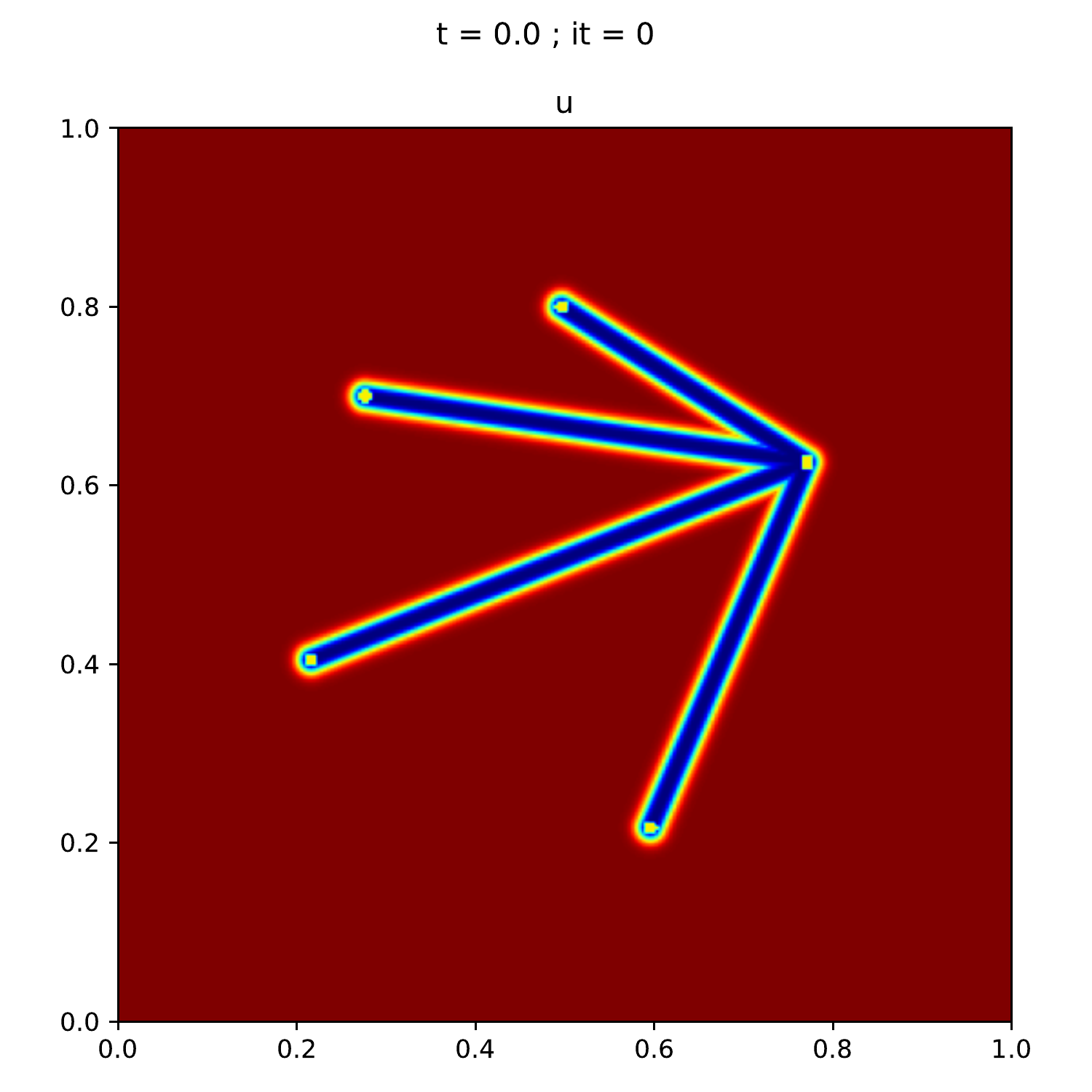}
   \includegraphics[width=0.22\textwidth]{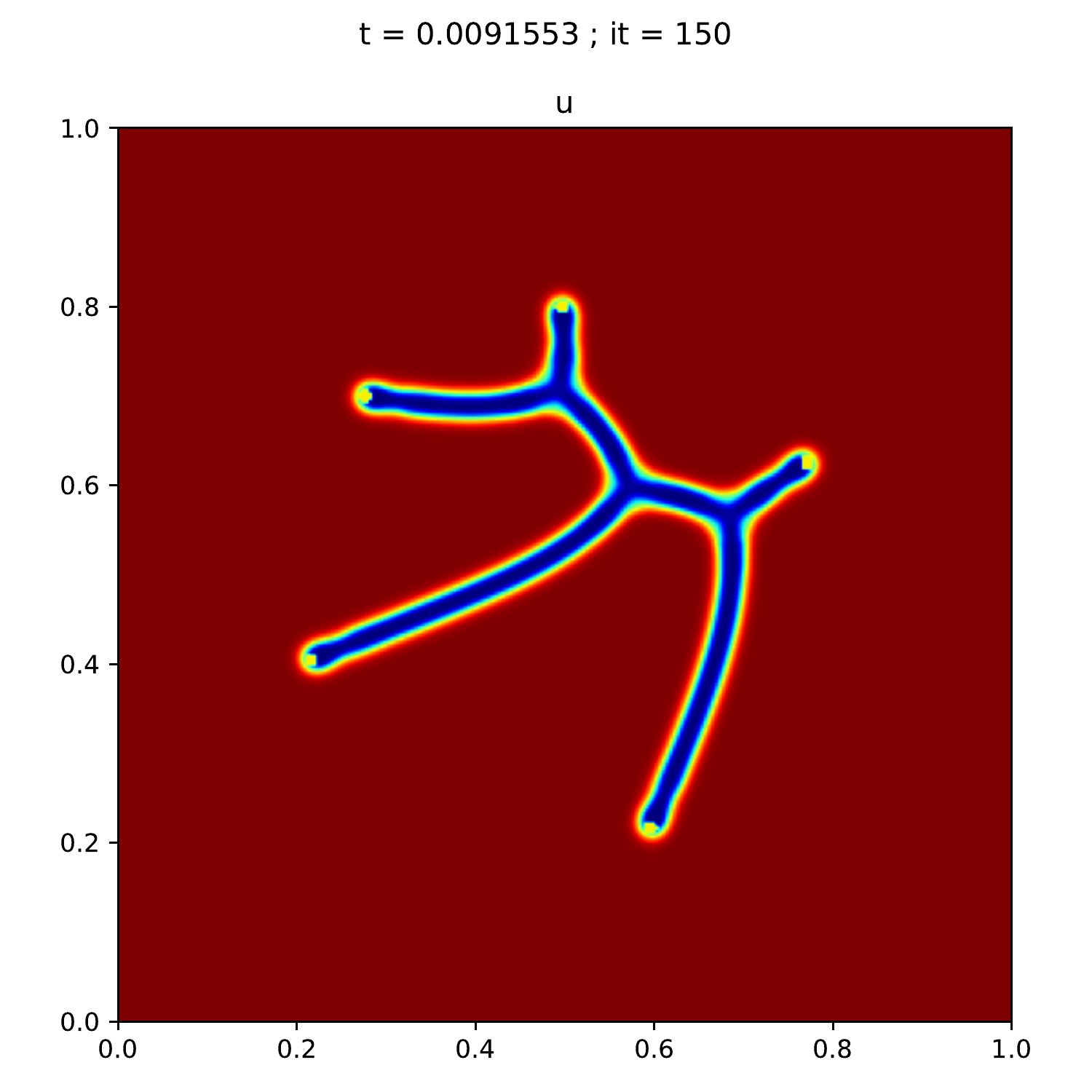}
   \includegraphics[width=0.22\textwidth]{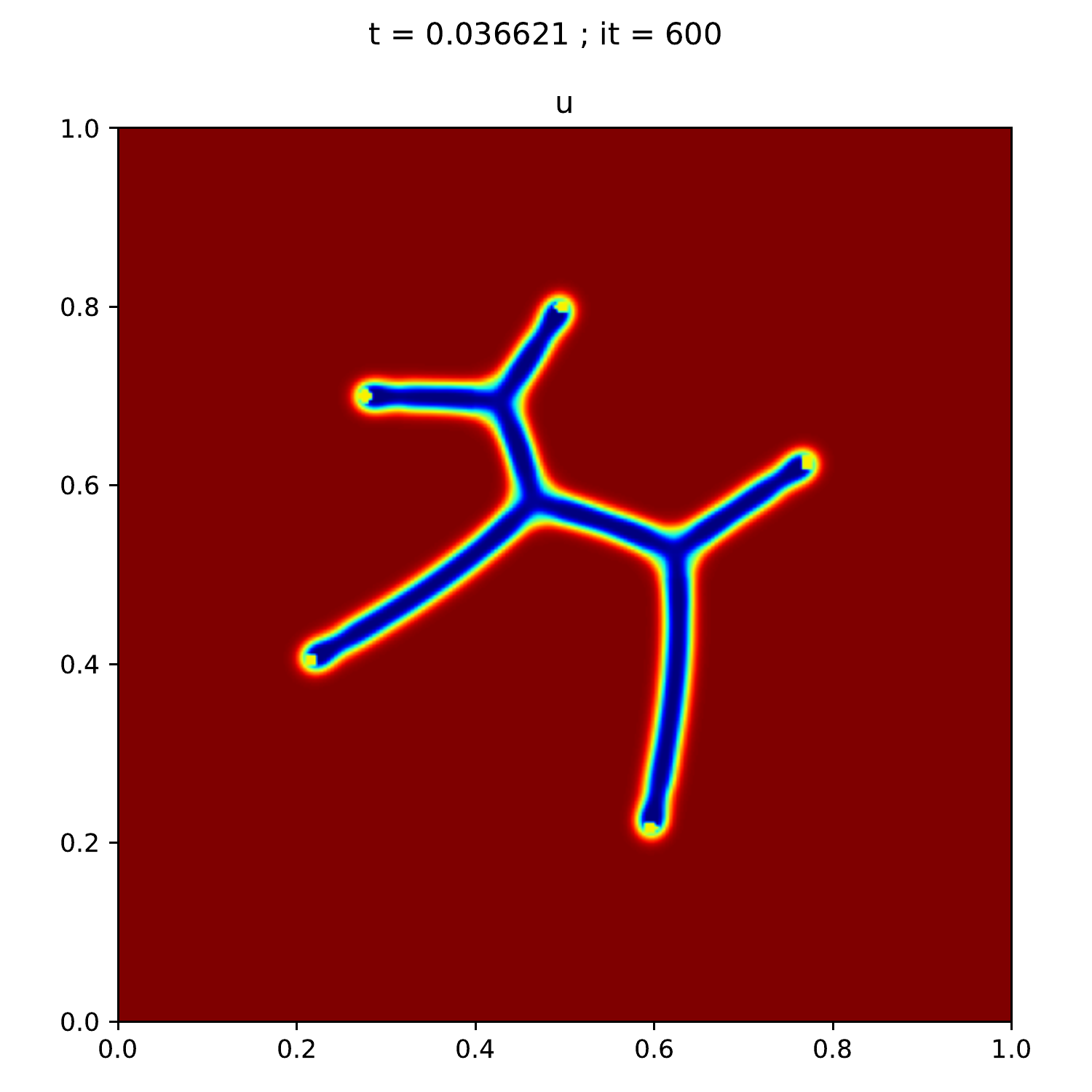}
   \includegraphics[width=0.22\textwidth]{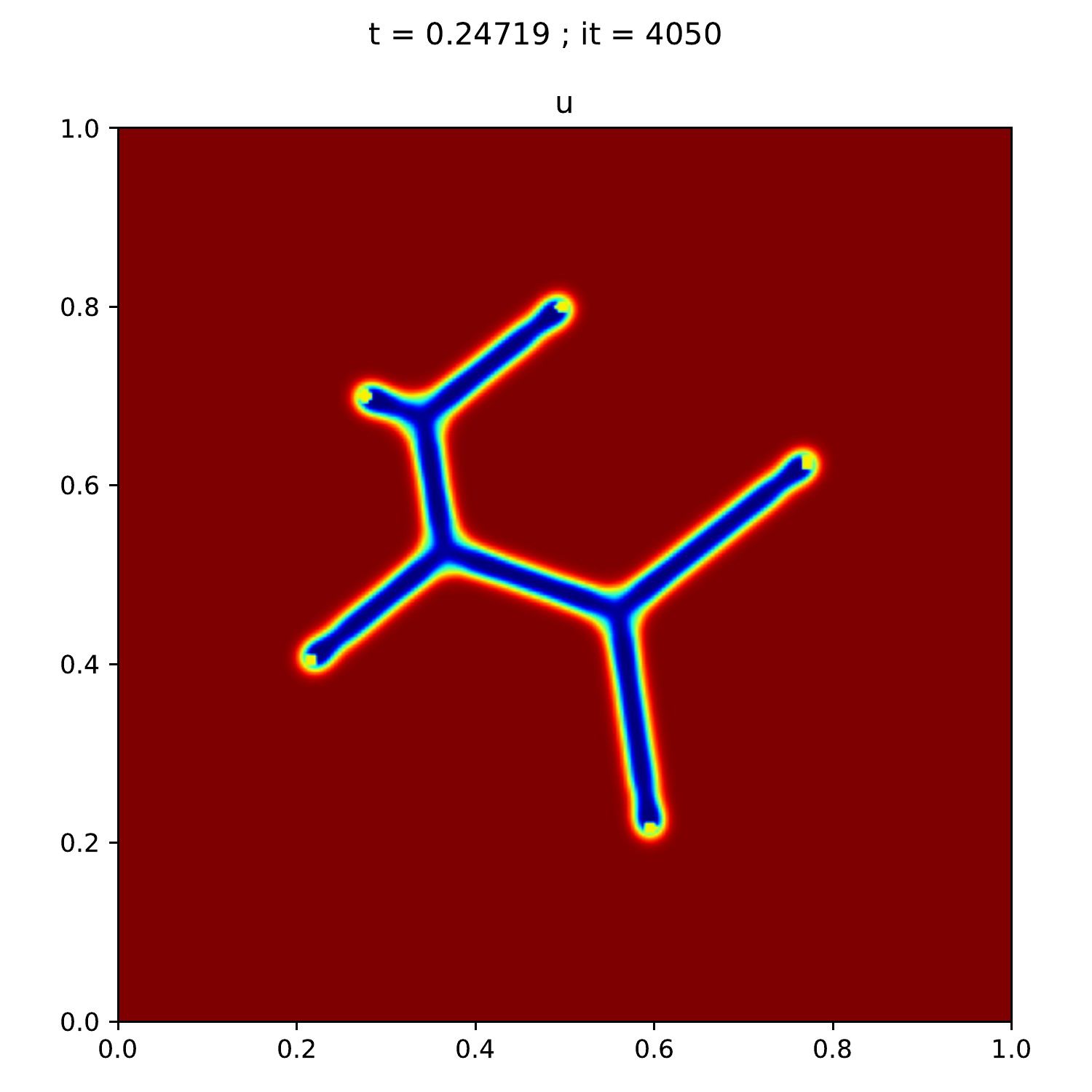}
   
   \includegraphics[width=0.22\textwidth]{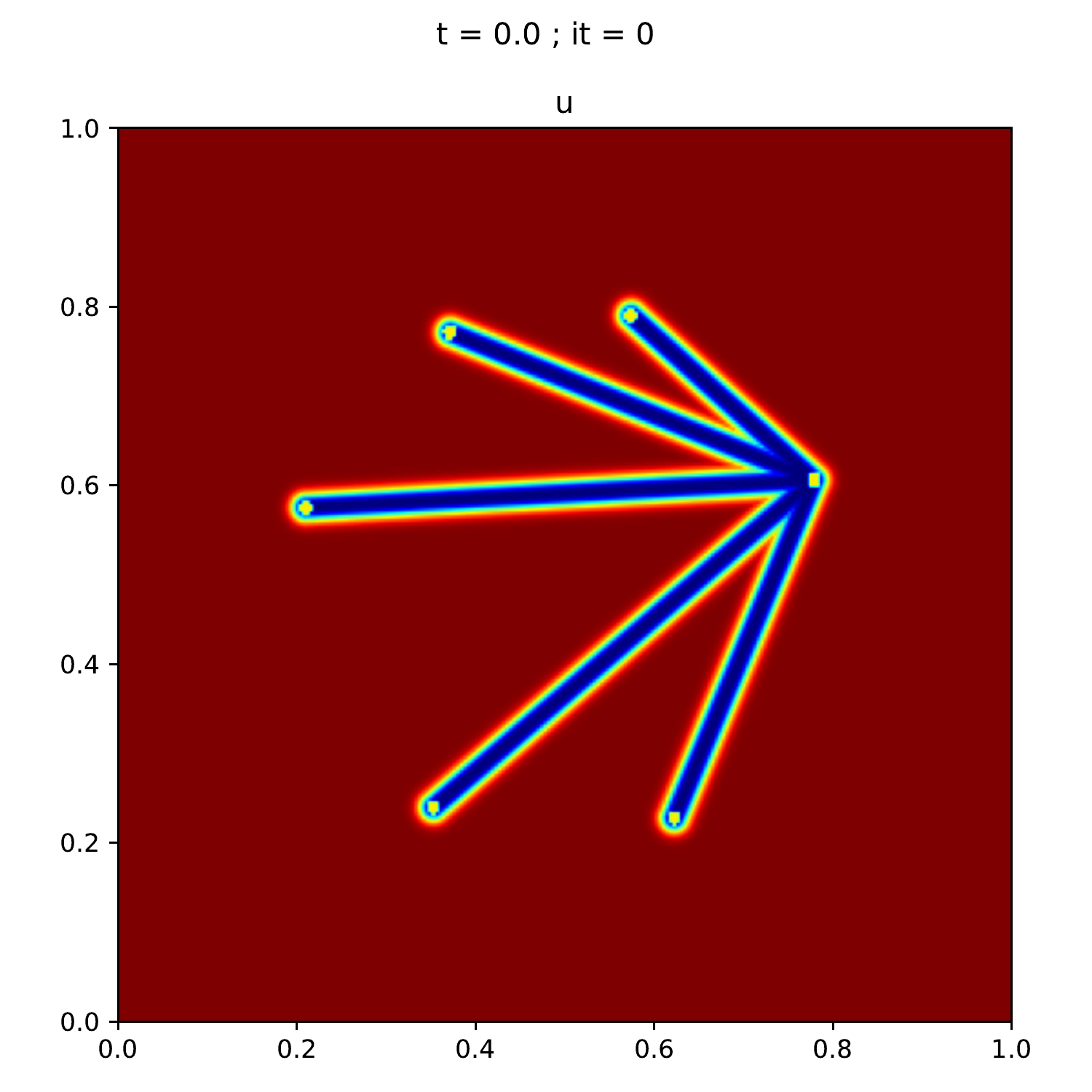}
   \includegraphics[width=0.22\textwidth]{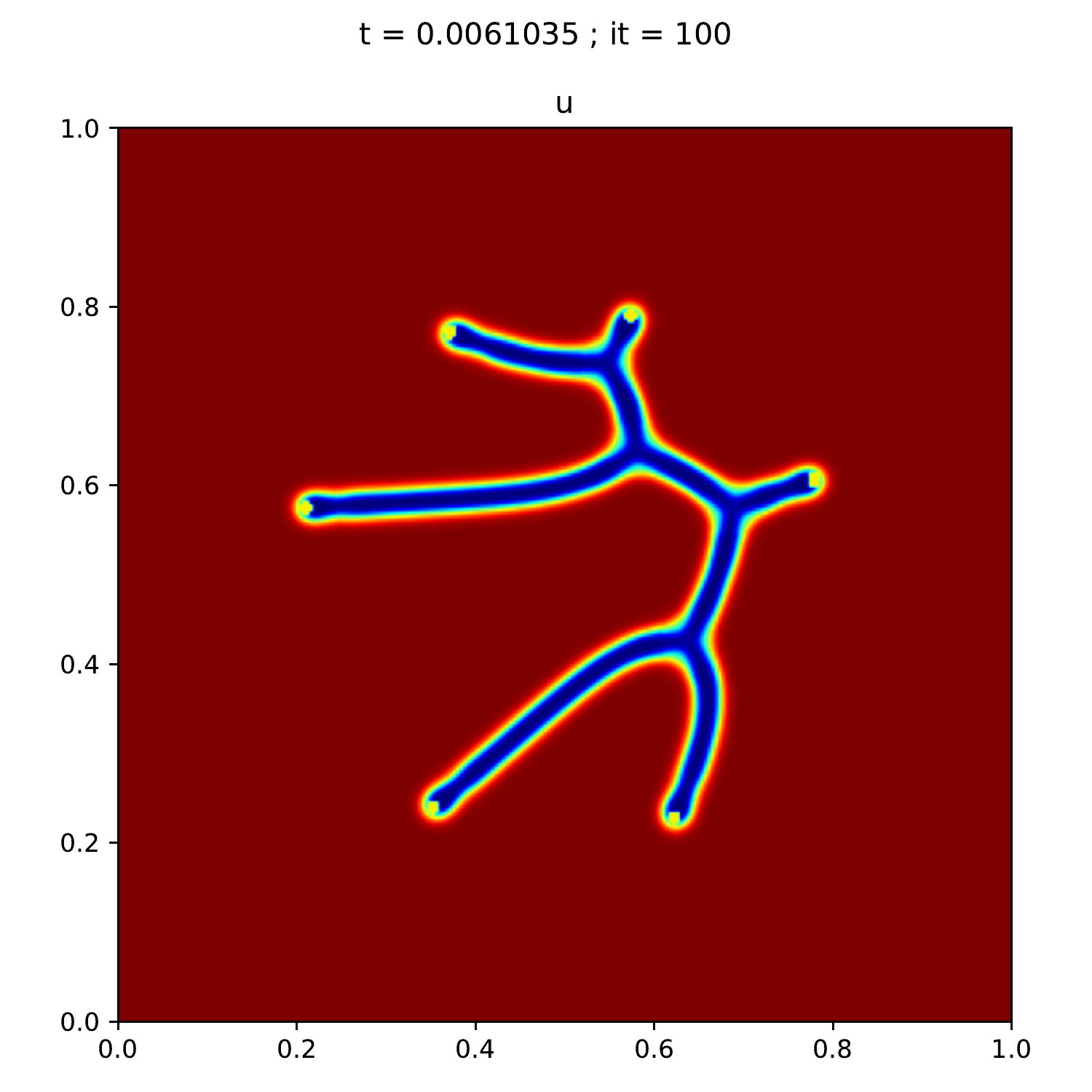}
   \includegraphics[width=0.22\textwidth]{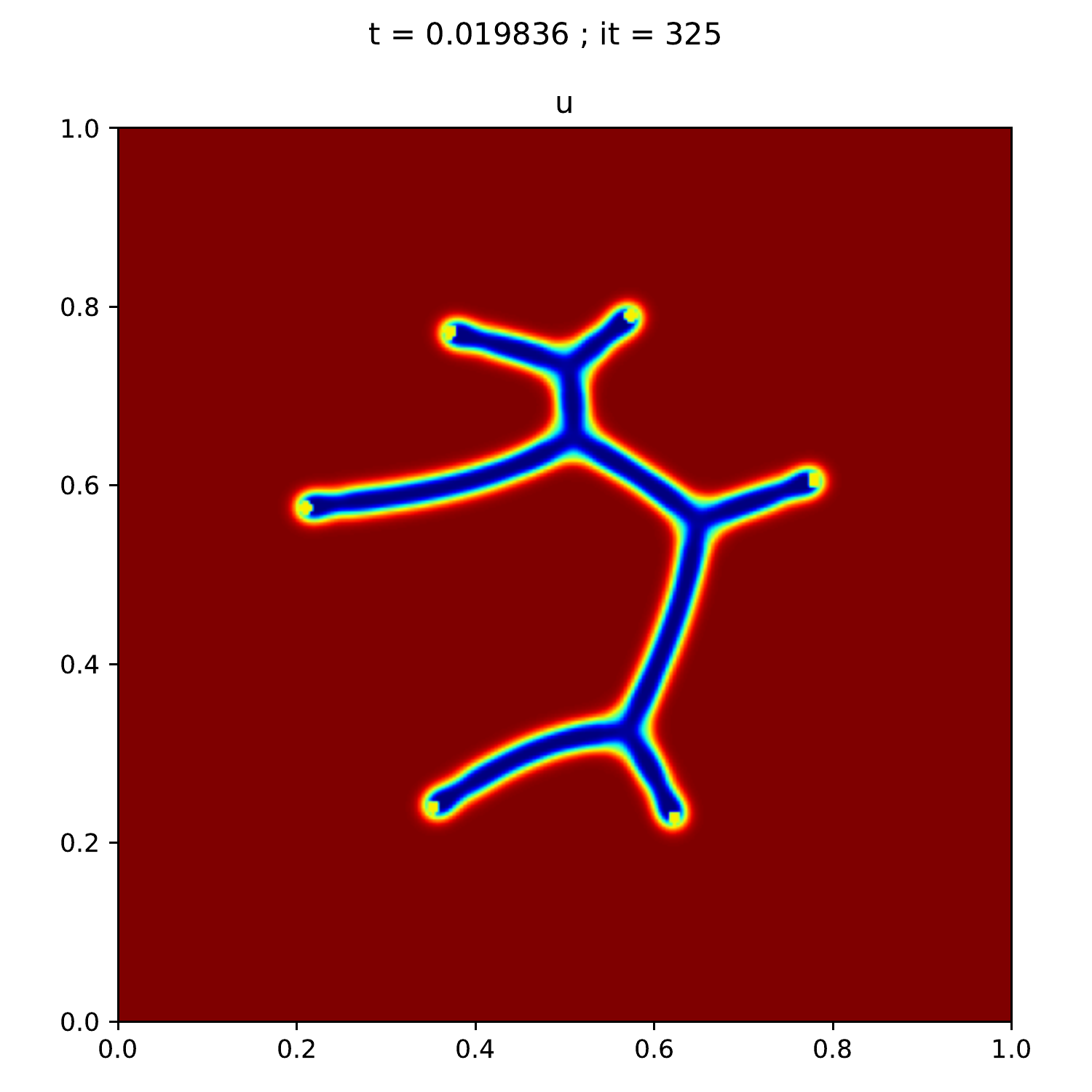}
   \includegraphics[width=0.22\textwidth]{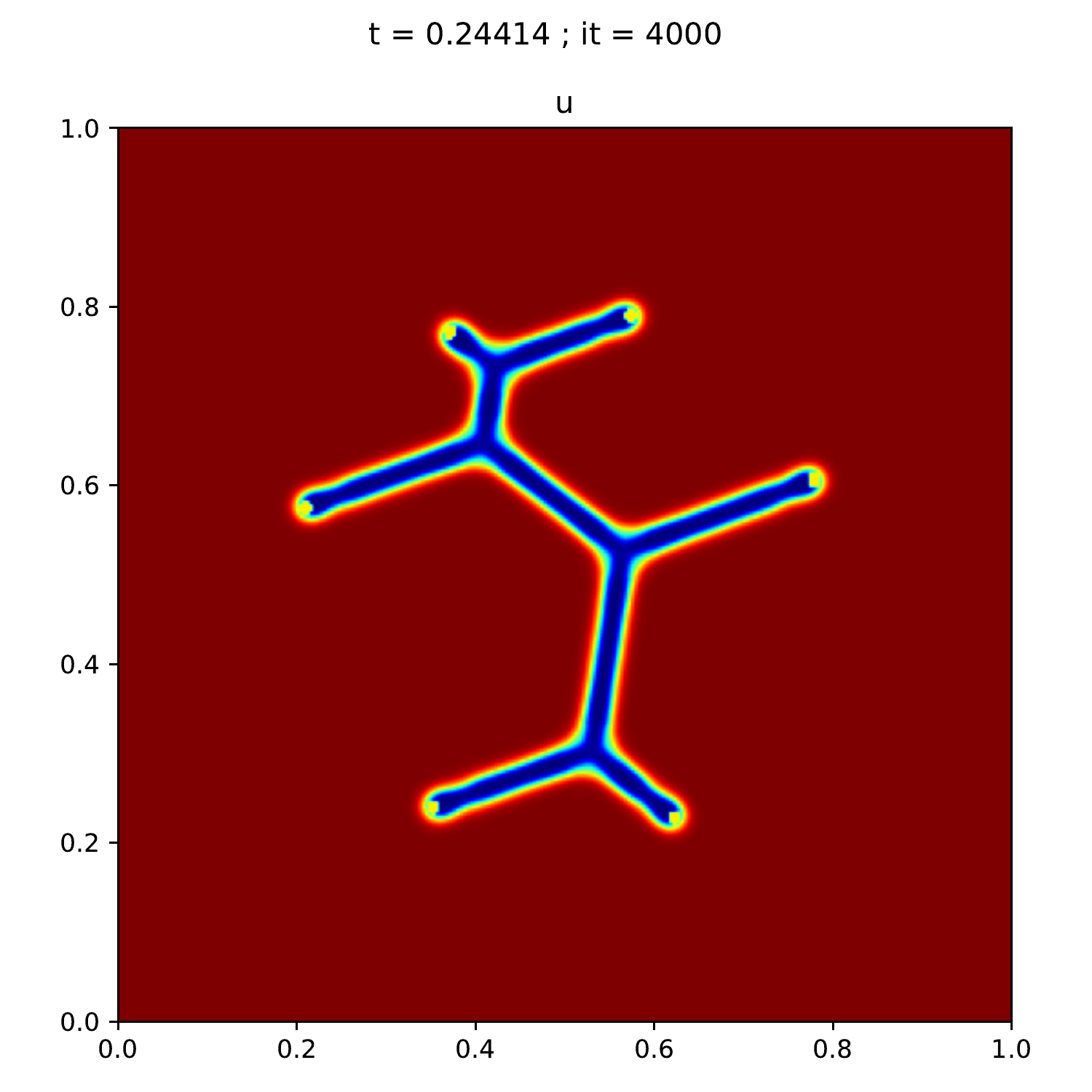}

    \caption{Approximation of Steiner trees in $2D$ using a non oriented mean curvature flow coupled with
    inclusion constraints according to the scheme \eqref{scheme:Non-oriented+inclusion}; Each line shows an 
    evolution of the numerical solution $u^n$ along the iterations
    with, respectively, $4$, $5$, and $6$ points $a_i$.}
    \label{fig:Steiner2D}
\end{figure}

We present in figure~\ref{fig:Steiner2D} three numerical experiments using, respectively, $4$, $5$, and $6$ points $a_i$ arranged non-uniformly on a circle. The first picture of each line corresponds to the initial set $\Gamma(0)$ 
constructed in such a way to connect the point $a_1$ in a linear way to all other points $a_i$, ${i \in \{ 2 \cdots N\}}$.
We then clearly observe an evolution of the set $\Gamma(t)$ which seems to converge to a Steiner tree connecting
all points $a_i$. In particular, the triple points seem to be handled accurately enough. In the end, the method seems extremely effective, simple to implement, and the computation of 
Steiner solutions is really fast
compared to other methods proposed in the literature. These results illustrate that our network trained on non-oriented mean curvature flows has sufficient numerical stability properties
to be coupled with inclusion constraints.     

\subsection{Approximation of minimal surfaces in $3D$} ~\\	   
The Plateau problem was formulated by Lagrange in 1760 and consists in showing the existence of a minimal surface in $\R^3$
with prescribed boundary $\Sigma$. Existence and regularity of solutions have been studied  in different contexts, 
see for instance  \cite{MR1501590} for smooth and orientable solutions and  \cite{MR117614}
for existence and uniqueness  of soap films  including orientable and non orientable surfaces and possibly multiple junctions. \\

From the numerical point of view, there exist many methods to compute minimal surfaces,
as for instance the first one \cite{MR1502829} and others  \cite{MR0458923,MR942778}   
which use a parametric representation of the surface. In the same spirit, the papers \cite{MR1613695,MR1613699} exploit
a finite element method to compute  numerical approximations of  minimal surfaces. 
Finally, note also some numerical approaches using an implicit representation of the minimal surface in \cite{MR1214016,MR2143330}
where a level set method is used. As for phase field approaches, a current-based method was proposed and analyzed in 
\cite{ChambolleFerrariMerlet2019-2}, and provides numerical approximations in the case of oriented surfaces. \\

We propose in this section to extend to $3D$ the phase field approach used for the Steiner problem. The idea is very similar in the sense that
\begin{enumerate}
 \item We  train a new network $\S^{NN}_{\theta,2}$ to approximate the semigroup $S_{\delta_t,\varepsilon,q'}$ in dimension $3$.
 Here we consider the case of the mean curvature flow of co-dimension $1$. The training of the network is performed on a
 database consisting of spheres involving under
 motion by mean curvature which is also explicit in this case. 
 \item We apply the previous scheme where the function $u_{in}$ is now defined by
 $$ u_{\text{in}}(x) = \sum_{i=1}^{N}  q'(\operatorname{dist}(\Sigma,x)/\varepsilon).$$
\end{enumerate}

\begin{figure}[htbp]
    \centering
     \includegraphics[width=0.2\textwidth]{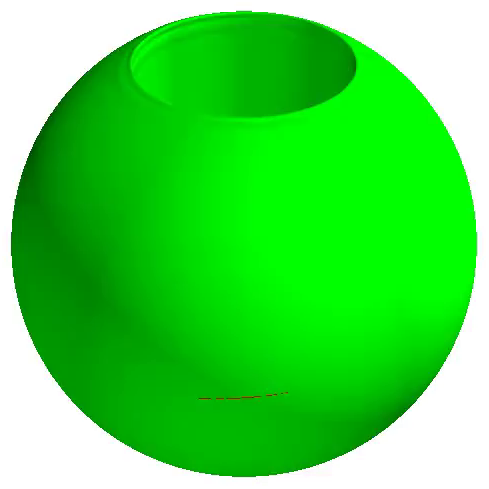}
     \includegraphics[width=0.2\textwidth]{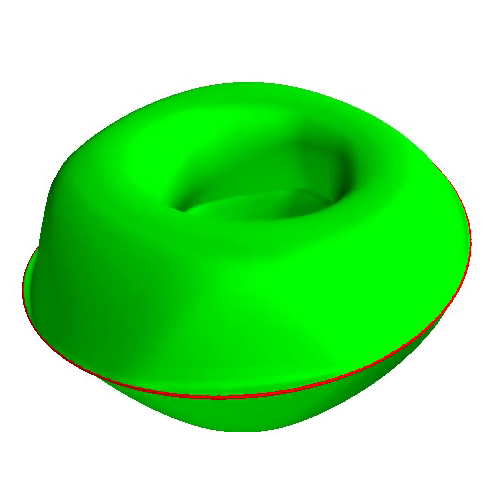}
     \includegraphics[width=0.2\textwidth]{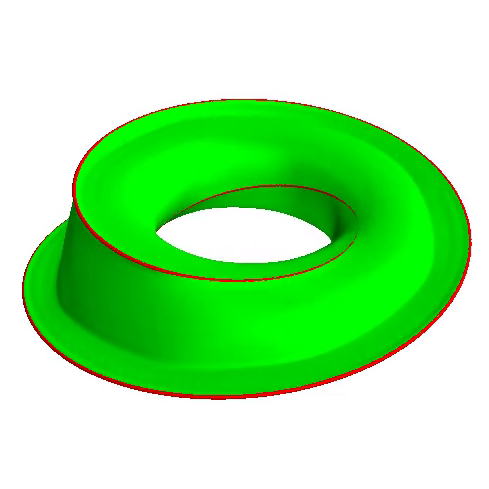}
    \includegraphics[width=0.2\textwidth]{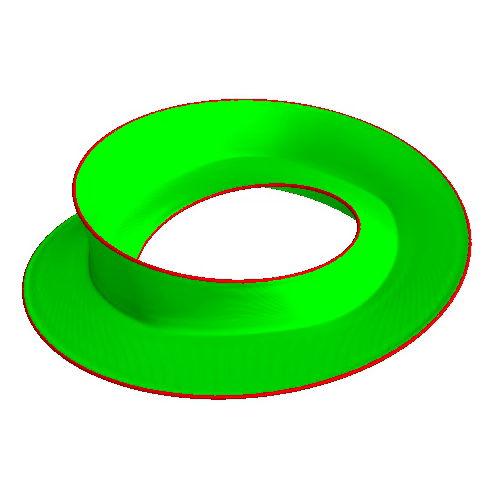}
    \caption{Approximation of a minimal surface using a non-oriented mean curvature flow coupled with an
    inclusion constraint. The images represent the numerical solution $u^n$ at four different times. 
    The prescribed boundary $\Sigma$  and  the boundary of $ \{ x \in Q ; u_n(x) \geq -0.2 \}$ are plotted in red and green, respectively .
    }
    \label{fig:Plateau3D_1}
\end{figure}

\begin{figure}[htbp]
    \centering
     \includegraphics[width=0.3\textwidth]{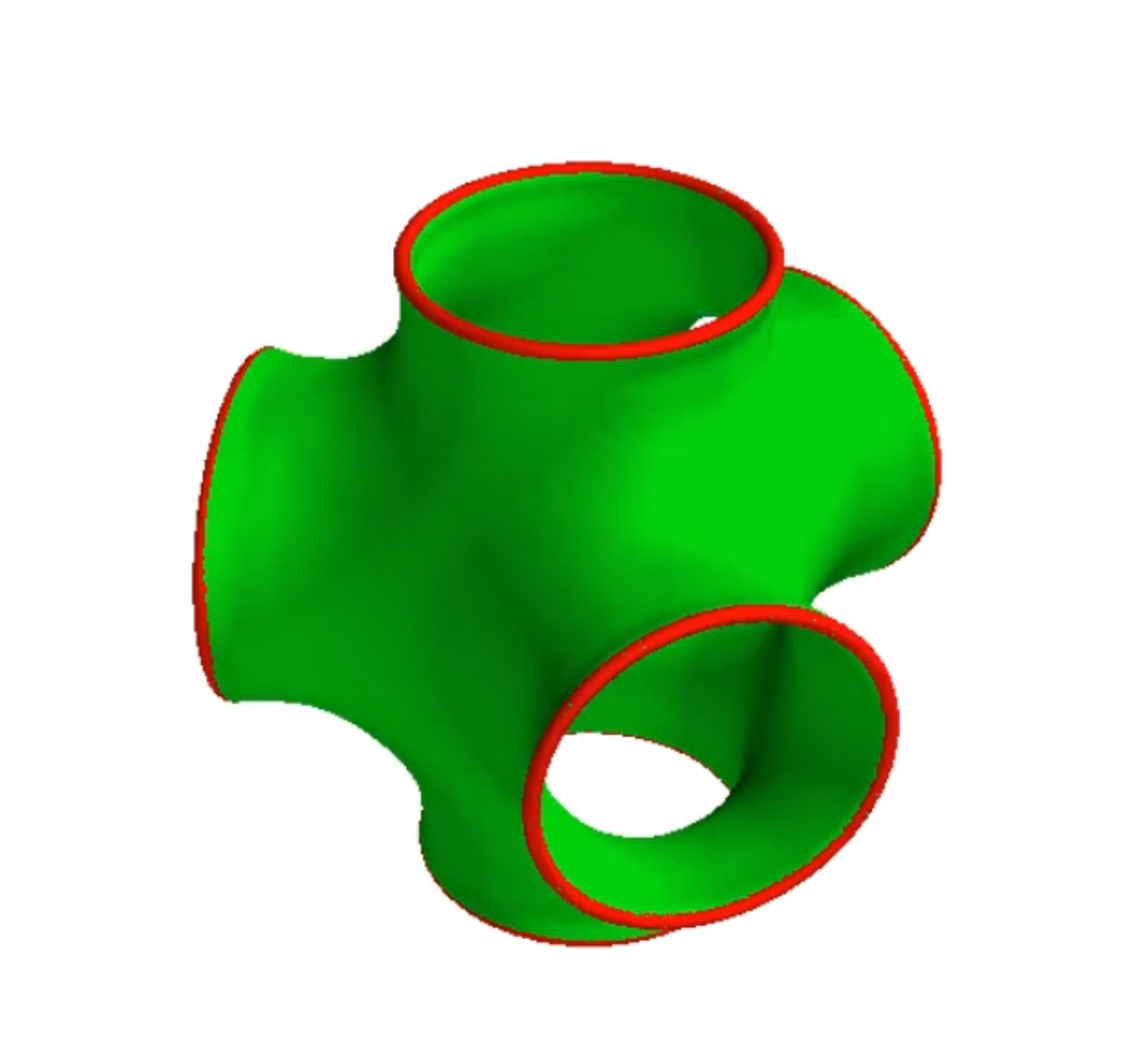}
    \includegraphics[width=0.3\textwidth]{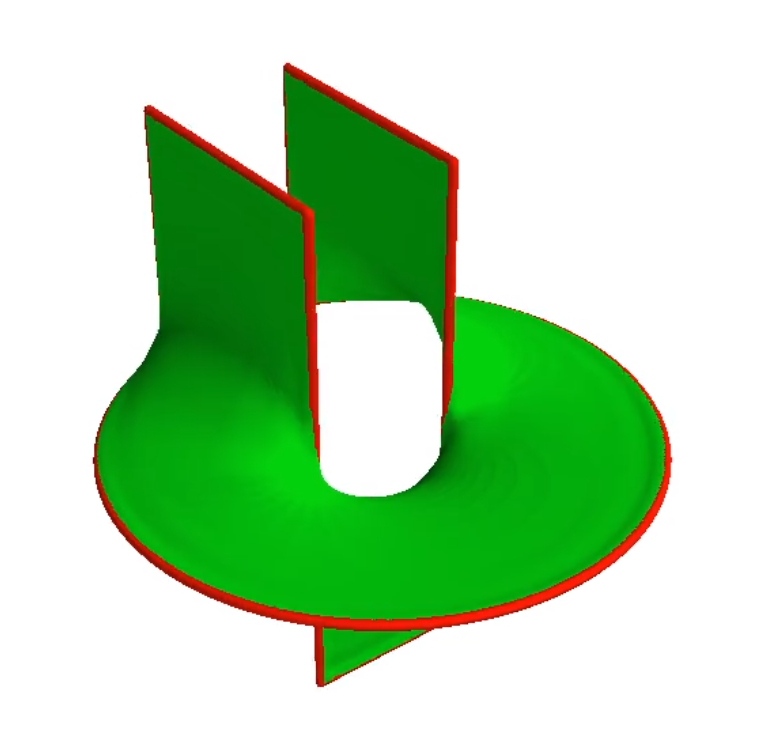}
     \includegraphics[width=0.3\textwidth]{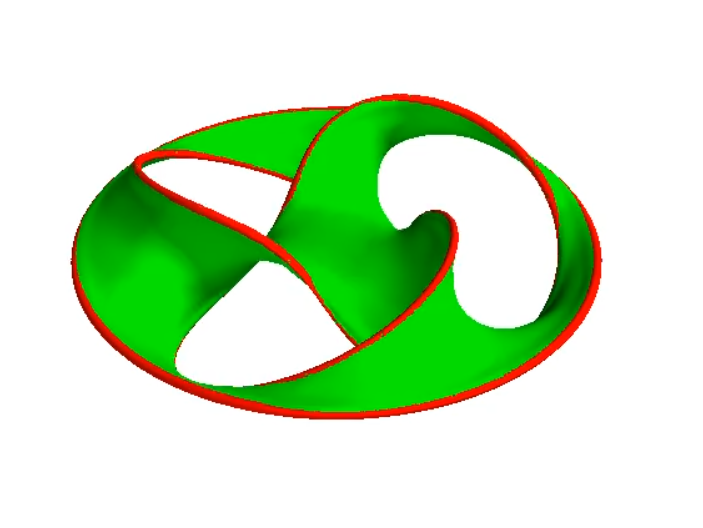}
    \caption{Approximation of minimal surfaces using a non-oriented mean curvature flow coupled with an
    inclusion constraint; Each picture shows the numerical stationary solution $u^{n}$ in green and the prescribed boundary $\Sigma$ in red.}
    \label{fig:Plateau3D_2}
\end{figure}

We present in figures~\ref{fig:Plateau3D_1} and~\ref{fig:Plateau3D_2} four numerical experiments using respectively
the boundary of a Mobius strip, a trefunknot, a union of rings, and a pit.
The prescribed boundary $\Sigma$  and  the boundary of $ \{ x \in Q ; u_n(x) \geq -0.2 \}$ are respectively plotted in red and green along
the iterations. We display in figure~\ref{fig:Plateau3D_1} the evolution of an interface along the iterations toward the stationary numerical solution. We show in figure~\ref{fig:Plateau3D_2}  only the stationary solutions, i.e. the minimal surfaces associated with the respective prescribed boundaries.

These numerical results are very encouraging because they are examples of non-orientable solutions with more 
or less complex topologies, and possibly with triple line singularities. 

\textcolor{black}{Note that the results obtained with our neural networks are perfectly consistent with those obtained using other methods, see for instance~\cite{MR4011534} for Steiner 2D trees and \cite{huang:tel-03584255} for minimal surfaces.}

\section{Conclusion}
We have introduced in this work neural networks that provide simple and very effective numerical schemes to accurately approximate 
motions by mean curvature in the classical case of oriented interfaces, as well as in the case of non orientable interfaces for 
which classical phase field approaches are not suitable.
The structures of our networks are inspired by discretization schemes of the Allen-Cahn equation that alternate linear diffusion and local nonlinear reaction. 
Our networks are very effective despite their low number of degrees of freedom (less than 1000) compared to standard networks such as U-net which has millions of 
parameters. 
This low complexity facilitates faster training and enforces robustness. 
This is best illustrated by the ability of our networks to generalize well despite a limited and simplistic training database.
The first numerical results presented in this paper are very encouraging.
In particular, the application to the Plateau problem of our numerical method based on neural networks shows that it can well approximate non-orientable minimal surfaces with triple junction lines. 
It is quite reasonable to expect that our approach can 
be extended to many other interface geometric flows, e.g., anisotropic mean curvature motion, 
surface diffusion, the Willmore flow, etc., possibly with various constraints. Furthermore, our approach seems quite generic and possible extension is expected 
to more general diffusion-reaction equations and the rich world of associated applications.
Lastly, since the network structures we propose can be easily coupled with constraints, it would be interesting to explore further in the near future the coupling with constraints such as bounds on solutions or the positivity of convolution kernels to ensure the monotonicity 
of the associated numerical schemes.

\section*{Acknowledgments}
The authors acknowledge support from the French National Research Agency (ANR) under grants ANR-18-CE05-0017 (project BEEP) and ANR-19-CE01-0009-01 (project MIMESIS-3D). Part of this work was also supported by the LABEX MILYON (ANR-10-LABX-0070) of Universit\'e de Lyon, within the program "Investissements d'Avenir" (ANR-11-IDEX- 0007) operated by the French National Research Agency (ANR), and by the European Union Horizon 2020 research and innovation programme under the Marie Sklodowska-Curie
grant agreement No 777826 (NoMADS).

\end{document}